\documentclass{amsart}
\usepackage{amsmath}
\usepackage{graphicx}
\usepackage{verbatim} 

\newtheorem{theorem}{Theorem}[section] 
\newtheorem{lemma}{Lemma}[section] 

\newtheorem{definition}{Definition}[section]

\newcommand{\bed}{\begin{definition}}
\newcommand{\eed}{\end{definition}}

\newcommand{\eps}{\epsilon}

\newcommand{\bitem}{\begin{itemize}}
\newcommand{\eitem}{\end{itemize}}

\newcommand{\goto}{\rightarrow}

\newcommand{\beqn}{\begin{equation}}
\newcommand{\eeqn}{\end{equation}}
\newcommand{\balign}{\begin{align}}
\newcommand{\ealign}{\end{align}}

\newcommand{\bR}{\bf R}

\def\cal{\mathcal}

\newcommand{\lp}{(LP)}
\newcommand{\p}{(P1)}

\usepackage{amssymb}
\begin{document}

\title[Observed Universality of Phase Transitions]{Observed Universality of Phase Transitions
in High-Dimensional Geometry, with Implications for
Modern Data Analysis and Signal Processing}

\author{David L. Donoho}
\address{Department of Statistics, Stanford University}
\curraddr{Department of Statistics, Stanford University}
\email{donoho@stanford.edu}
\author{Jared Tanner}
\address{School of Mathematics, University of Edinburgh}
\curraddr{School of Mathematics, University of Edinburgh}
\email{jared.tanner@ed.ac.uk}
\thanks{The authors thank
the Isaac Newton Mathematical Institute 
 for hospitality during the programme
{\it Statistical Theory and Methods for Complex, High-Dimensional Data},  
and for a Rothschild Visiting Professorship held by DLD.  The authors 
thank Erling Andersen for donating licenses 
for the Mosek software package which saved a great deal of computer time
in the studies described here.
This work has made use of the resources provided by the Edinburgh
Compute and Data Facility (ECDF). 
DLD was partially 
supported by NSF DMS 05-05303,  and  JT was partially supported by Sloan
and Leverhulme Fellowships.}

\date{14 May 2009}

\begin{abstract}

We review connections
between phase transitions in high-dimensional
combinatorial geometry and
phase transitions occurring in modern high-dimensional
data analysis and signal processing.
In data analysis, such transitions arise
as abrupt breakdown of linear model selection,
robust data fitting or compressed sensing reconstructions,
when the complexity of the model or the number of outliers
increases beyond a threshold.  
In combinatorial geometry these transitions appear
as abrupt changes in the properties of face
counts of convex polytopes when the dimensions
are varied.  The thresholds in
these very different problems appear in the same critical
locations after appropriate
calibration of variables.

These thresholds are important in each
subject area: for linear modelling, they place hard limits
on the degree to which the now-ubiquitous high-throughput data analysis
can be successful; for robustness, they place hard limits
on the degree to which standard robust fitting methods
can tolerate outliers before breaking down; for
compressed sensing, they define the sharp boundary
of the undersampling/sparsity tradeoff curve in
undersampling theorems.

Existing derivations of phase transitions in combinatorial geometry
assume the underlying matrices have independent and 
identically distributed (iid) Gaussian elements.  In
applications, however, it often seems that Gaussianity is not required.
We conducted an extensive
computational experiment  and formal
inferential analysis to test the hypothesis
that these phase transitions are {\it universal}
across a range of underlying matrix ensembles.
We ran millions of linear programs using random
matrices spanning several matrix ensembles and problem sizes; to the naked
eye, the empirical phase transitions do not depend on the
ensemble, and they agree extremely
well with the asymptotic theory assuming Gaussianity.
Careful statistical analysis reveals discrepancies
which can be explained as transient terms, decaying
with problem size. The experimental results
are thus consistent with an asymptotic large-$n$
universality across
matrix ensembles;  finite-sample universality can
be rejected.

\end{abstract}

\maketitle



\newcommand{\ntrain}{n}
\section{Introduction}
\setcounter{equation}{0} 

Recent work has exposed a phenomenon
of  abrupt {\it phase transitions} in high-dimensional
geometry. The phase transitions
amount to rapid shifts in the likelihood of a 
property's occurrence when a dimension parameter
crosses a critical level (a {\it threshold}).
We start with  a concrete example, and then
identify surprising parallels in data analysis and signal processing.

\subsection{Convex hulls of Gaussian point clouds}
\label{ssec-ConvGauss}

Suppose we have a sample $X_1$, \dots , $X_n$
of independent standard normal random variables
in dimension $d$, forming a point cloud 
of $n$ points in $\bR^d$.  Our intuition suggests that a
few of the points will lie on the `surface' of the dataset,
that is, the boundary of the convex hull; the rest will lie  `inside', i.e. interior to the hull.
However, if $d$ is a fixed fraction of $n$ and both are large,
our intuition is completely violated. Instead, {\it all}
of the points are on the boundary of the convex hull -- {\it none} is interior.
Moreover, the line segment connecting the typical
pair of points {\it does not} intersect the interior;
in complete defiance of expectation, it
stays on the boundary.
Even more, for $k$ in some appreciable range,
{\sl the typical $k$-tuple spans a convex hull which does not intersect the
interior}! For humans stuck all their lives in three-dimensional space, 
such a situation is hard to visualize.

The phenomenon of phase transition appears as follows: such seemingly strange behavior
continues for quite large $k$, up to a {\it predictable} threshold
given by a formula
$k^* = d \cdot \rho(d/n;T)$, where $\rho()$ is 
defined in \S \ref{sec:GeoComb} below.
Below this threshold (i.e. $k $ a bit smaller than $k^*$),
the strange behavior is observed;
but suddenly, above this threshold (i.e. for $k$ a bit larger) our normal
low-dimension intuition works again -- convex hulls of $k$-tuples
of points indeed intersect the interior.

This curious phenomenon in high-dimensional geometric probability
is one of a small number of fundamental
such phase transitions.  We claim they have consequences
in several applied fields:
\bitem
\item  in selecting models for statistical data analysis of large datasets, 
\item  in coping with outlying measurements in designed experiments,
\item  in determining how many samples we need to take in
          designing imaging devices.    
\eitem 
The consequences  can be both profound and important. They range 
from negative-philosophical -- if your database has too many 'junk' variables in it
nothing can be learned from it --
to positive-practical -- it isn't really necessary to sit cooped up for an hour in a medical
MRI scanner: with the right software, the necessary  data
could be collected in a fraction of the time commonly used today.

Our paper will help the reader understand
more precisely what these phase transitions are
and where they may occur in science and technology;
it will then discuss our
\begin{quotation}
{\bf Main contribution.}
We have observed a {\it universality}
of threshold locations across
a range of underlying probability distributions. 
We are able to change the underlying distribution from
Gaussian to any one of a variety of non-Gaussian choices, and
we still observe phase transitions at the
same locations.
\end{quotation}
We compiled evidence based on
millions of random trials and 
observed the same phase transitions  even
for several  highly non-iid ensembles.
We here formally state and test the universality hypothesis.

Our research leads to an intriguing challenge
for high-dimensional geometric probability:
\begin{quotation}
{\bf Open problem.}
Characterize the {\it universality class} containing the standard Gaussian:
i.e. the class of matrix ensembles leading to
phase transitions matching those for Gaussian polytopes.
\end{quotation}

Evidently this class is fairly broad.
In view of the significance of these phase transitions
in applications, this is quite an attractive challenge.
We begin by illustrating three surprising appearances
of these phase transitions.

\subsection{First surprise: model selection with large databases}
\label{ssec-Surp1}

A  characteristic feature of today's {\it data deluge}  is the tendency 
in each field to  collect ever
more and more measurements on each observed entity,
whether it be a pixel of sky, a sample of blood or a sick patient.
Technology continually puts in our hands
{\it high-throughput} measurement equipment making ever more
varied and ever more detailed numerical measurements on the spectrum of light,
the protein expression in whole blood or fluctuations in neural
or muscular activity.

As a result, observed entities
are represented by ever higher-dimensional feature
vectors.  In fact the transition between the 20th and
21st centuries marked a sudden increase in the dimensionality
of typical datasets that scientists studied, so that
it became unremarkable for each observational unit to
be represented as a data point in a $p$-dimensional
space with $p$ very large -- in the hundreds, thousand or millions.

The modern trend to high-throughput measurement
devices often does not address the fundamental difficulty of obtaining
good observational units.   Scientists face the same troubles
they always have faced when  searching for
 subjects affected by a rare disease, 
or observing rare events in distant galaxies.  Hence, in many fields 
the number of observational units
stays small, perhaps in the  hundreds (or even dozens),
but each of those few units can now be routinely subjected 
to unprecedented density of
numerical description.

Orthodox statistical methods assumed a quite
different set of conditions:  an abundance of observational units
and a very limited set of measured characteristics on each unit.
Modern statistical research
is intensively developing new tools
and theory to address the new unorthodox setting; such
research comprised much of the activity  in the 
recent 6-month Newton Institute programme 
{\it Statistical Theory and Methods for Complex,   High-Dimensional Data}.

Consider a linear modelling scenario going back to
Legendre, Gauss and perhaps even before.
We have available a response variable $Y$
which we intend to model as a linear function of 
up to $p$ numerical predictor variables
$X_1$, \dots, $X_p$.   We contemplate
an utterly standard multvariate linear model,
$ Y = \alpha + \sum_j \beta_j X_j + Z$ where the $\beta_j$ are
regression coefficients and $Z$ is a standard normal 
measurement error.  In words, the expected value
of $Y$ given $\{X_j \}$  is a linear combination
 with coefficients $\beta_j$.
 
Suppose we have a collection of measurements
$(Y_{i},X_{i,j}, j=1,\dots,p)$, one for each observational
unit $i$.  We will use these data to
estimate the $\beta_j$'s, allowing future
predictions of $Y$ given the $X_j$'s.

In the `21st Century Setting' described above,
we have {\it more predictor variables than observations},
meaning $p > n$.  While Legendre and Gauss may have
understood the $p < n$ case, they would have been
very troubled by the $p > n $ case: there are more
unknowns than equations, and there is noise to boot!

A key feature of high-throughput analysis is that
batches of potential predictors are automatically measured
but one does not know in advance which, if any, may be useful
in a particular project.
Researchers in applied sciences
where high throughput studies are
popular (e.g. genomics, proteomics, metabolonomics) 
believe that some small fraction of the measured
features  are useful, among many useless ones.
Unfortunately,  high-throughput techniques
give us everything, useful and useless, all
mixed together in one batch.

In this setting, a reasonable response is {\it  forward stepwise linear regression}.
We proceed in stages, starting with
the simple model $Y = \alpha + Z$ (i.e. no dependence on
$X$'s) and at each stage expand the model by
adding the single variable  $X_j$ offering the strongest improvement in
prediction.  

\begin{figure}
\begin{center}
\includegraphics[width=3in,height=2.1in]{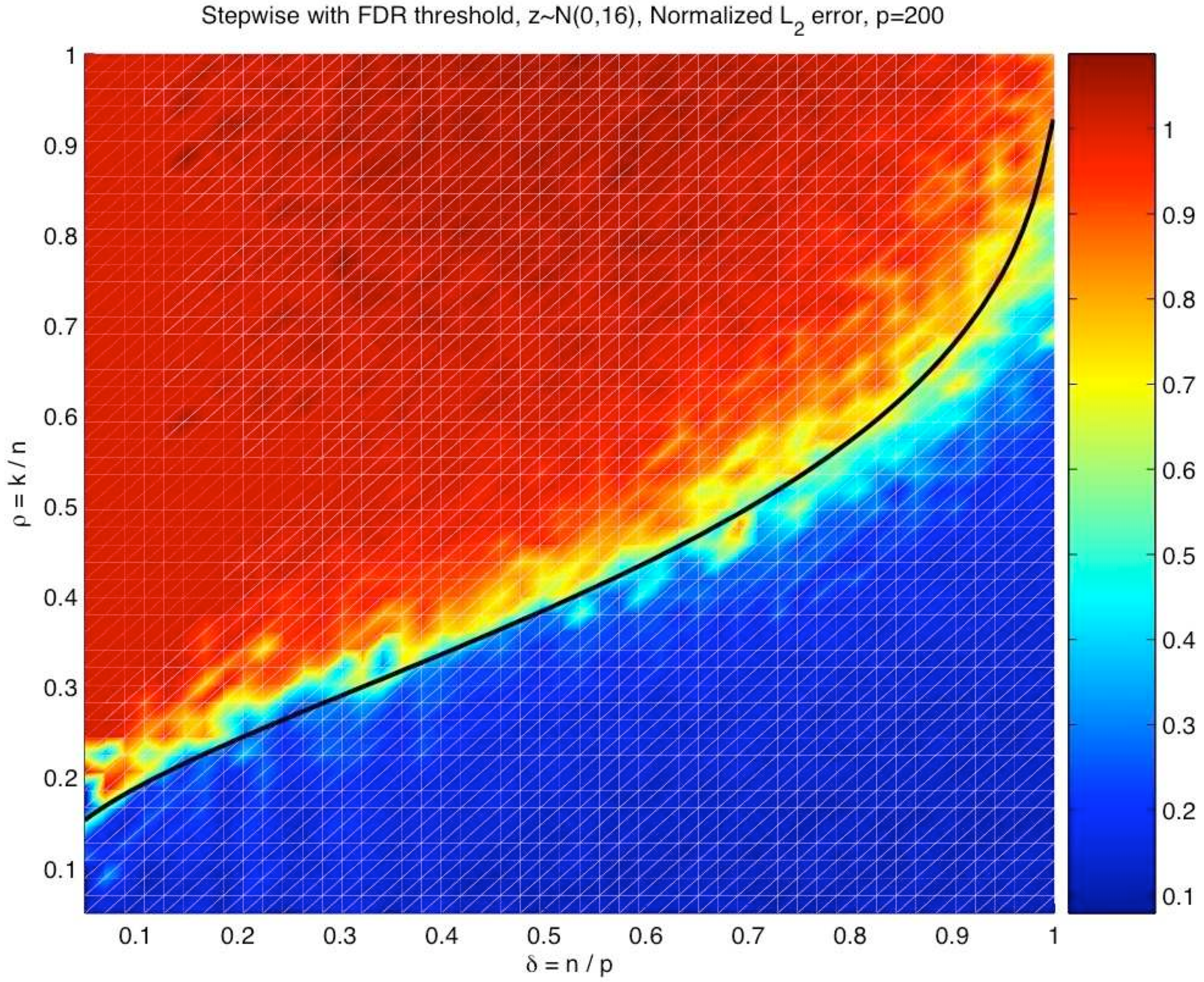}
\end{center}
\caption{{\it Phase Transition in Stepwise Regression Performance}.
  Average fractional error of the regression coefficients from Forward
  Stepwise Regression with a False Discovery Threshold,
  $\|\hat{\beta}-\beta\|/\|\beta\|$ (Donoho \& Stodden 2006). 
  Horizontal axis: $\delta = n/p$ (number of observations/number of variables). 
  Vertical axis $\rho = k/n$ (number of useful variables / number of
observations).  Solid curve $\rho(\delta; C)$ from combinatorial geometry. 
The rapid deterioration in performance roughly coincides with this curve} \label{fig:stodden}
\end{figure}

Donoho \& Stodden (2006) conducted a simulation experiment
using forward stepwise regression with
False-Discovery-Rate control stopping rule. Their experiment
chose $p=200$ and explored a range of $n < p$ cases.
Letting $k$ denote the number of useful
predictors among the $p$ potential predictors,
they set up true underlying regression models
reflecting the choice $(k,n,p)$, ran the stepwise regression
routine, and recorded the mean squared prediction error
of the resulting estimate.  Figure \ref{fig:stodden} displays
results: the coloured attribute gives the relative mean-squared
error of the estimate; the axes present the ratio $\delta = n/p$
of observations to variables, with $0 < \delta < 1$ in this brave new
world, and $\rho =k/n$ of useful variables to observational units.
Evidently there is an abrupt change in performance: one
can suddenly `fall off a cliff' by slightly increasing the number of
useless variables per useful variable.
 The {\it surprise} is that
the `cliff'  is roughly at the same position as 
the overlaid curve.  That curve, denoted by $\rho( \delta ; C)$
and defined fully below,  derives from combinatorial
geometry\footnote{The auxilary parameter $C$ in $\rho(\delta; C)$ 
is used to indicate the connection of this curve with the standard 
cross-polytope $C$, defined in \eqref{eq:crosspolytope}, from which 
it is derived.}
notions similar to those in \S \ref{ssec-ConvGauss}.

Interpretation: 

\bitem
\item A standard scientific data analysis approach in the `21st century setting'
`falls off a cliff', failing 
abruptly when the model becomes too complex.
\item The location of this failure (ratio of model variables to
observations)  matches a curve derived
from the field of geometric combinatorics!
\eitem

\subsection{Surprise 2:  robustness in designed experiments}
\label{ssec-Surp2}

We now consider a problem in robust statistics.
Suppose that  a response variable $Y$
is thought to depend on $p$ independent variables
$X_1$, \dots, $X_p$.
Unlike the data-drenched
high-throughput observational studies
of \S \ref{ssec-Surp1} we are in a classical
designed experiment, with $p < n$.  The dependence
is linear, so we again have $Y = \sum_j \beta_j X_j + Z + W$.
The error $Z$ is again normal,
but $W$ is a `wild' variable containing occasional
very large  outliers.  

Most scientists realize that such outliers could upset the usual least-squares
procedure for estimating $\beta$, and many know that $\ell_1$ minimization,
\begin{equation}\label{eq:l1Robust}
     \min_b \| Y - X'b \|_1,
\end{equation}
is purported to be  `robust', particularly in designed experiments
where wild $X$'s do not occur.
Let us focus on a specific designed experiment, 
where $X$ is an $n$ by $p$ partial Hadamard matrix,
i.e. $p$ columns chosen at random from an $n \times n$ Hadamard matrix.
This design chooses $X_{i,j}$ either  $+1$ or $-1$ 
in a very specific way; there are no wild $X$'s.
We let the outlier generators
$W$ have most entries $0$, but, in our study,
a small fraction $\eps = k/n$ -- at randomly-chosen
sites -- will have very large values; in any one realization
they can either be all large positive or all large negative. 

\begin{figure}
\begin{center}
\includegraphics[bb=76 215 543 691,width=2.3in,height=2.1in]{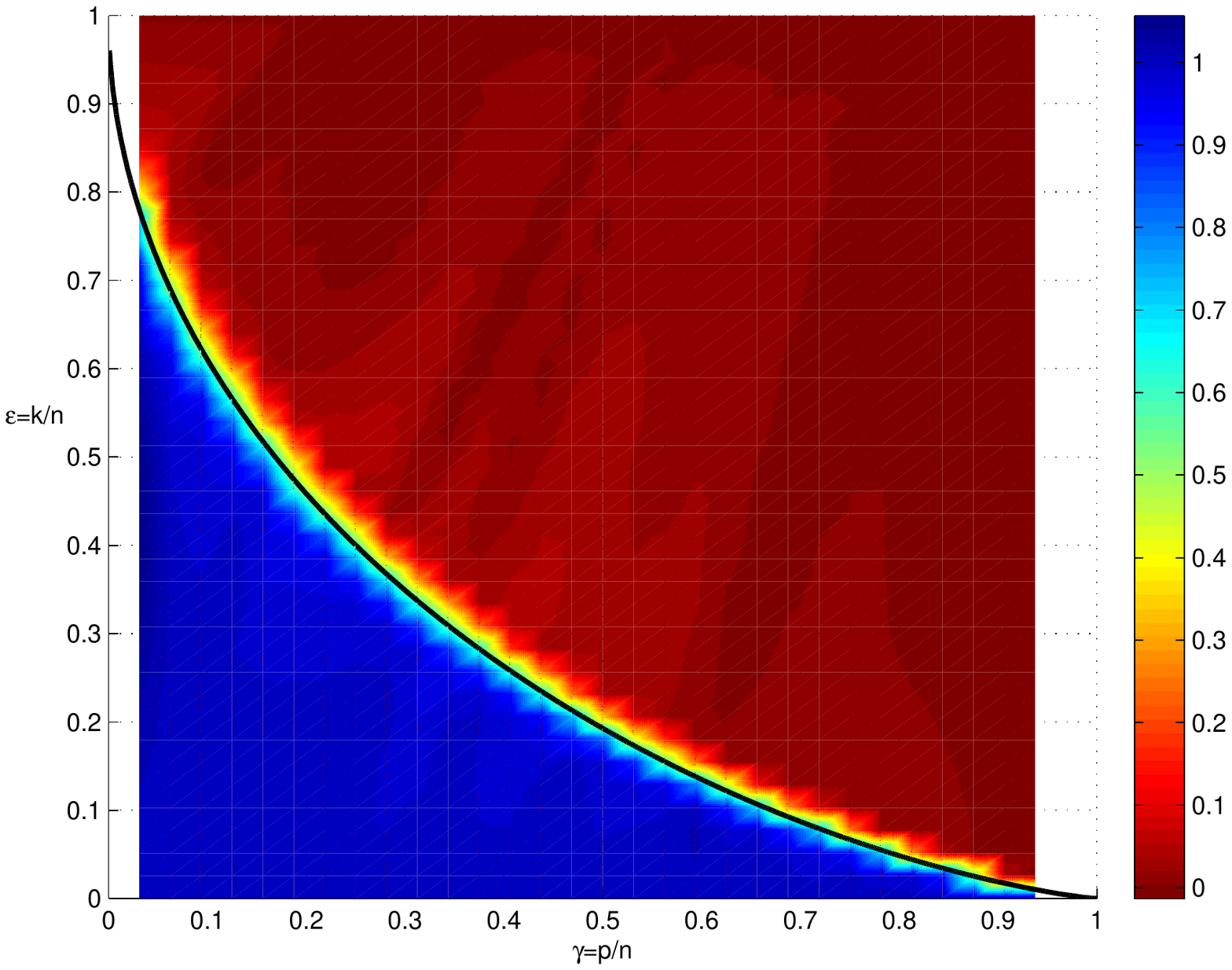} 
\includegraphics[bb=76 215 543 691,width=2.3in,height=2.1in]{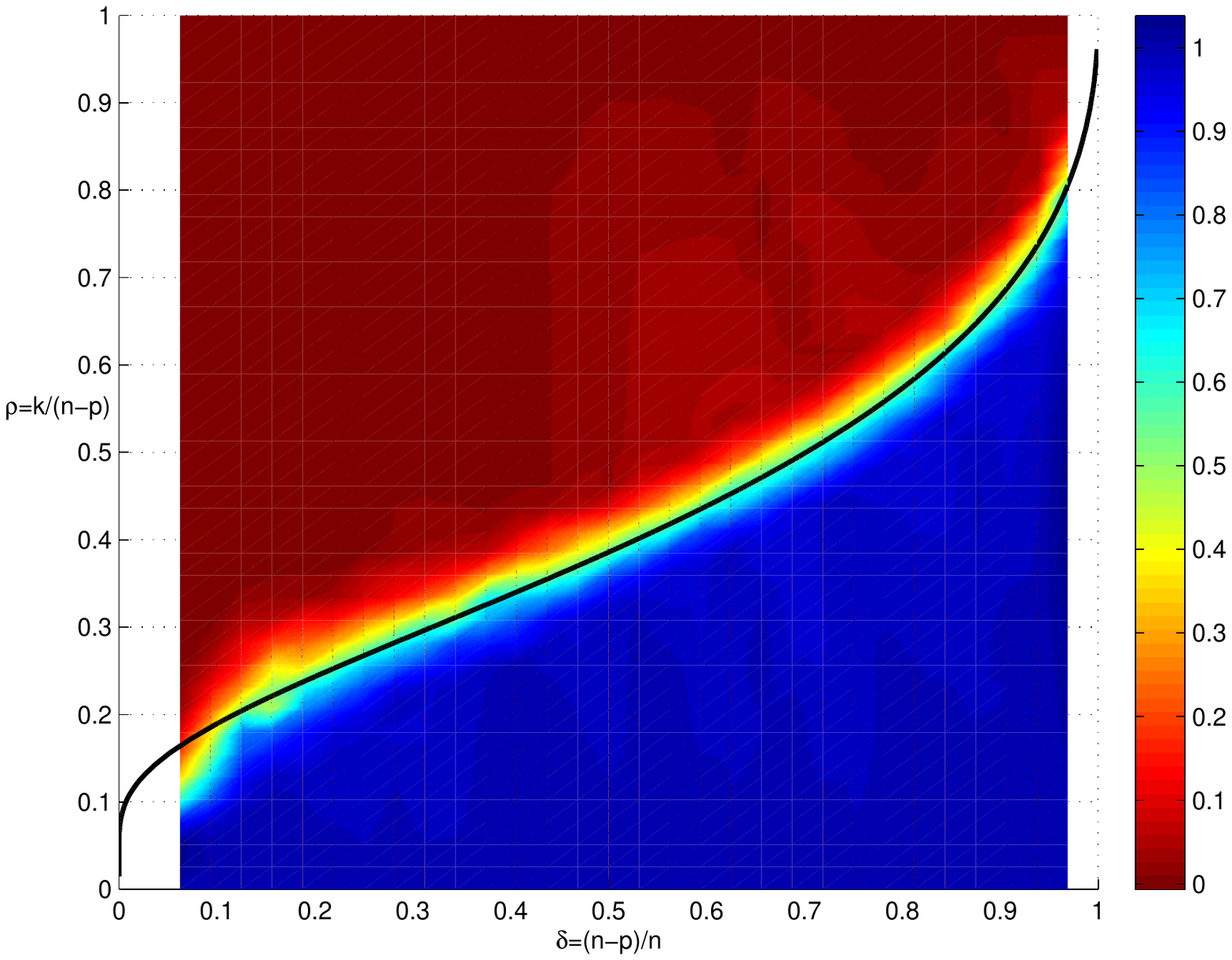} 
\end{center}
\caption{{\it Breakdown Point of $\ell^1$ Fitting}.  
Shaded attribute: fraction of realizations in which regression coefficients from 
\eqref{eq:l1Robust} are accurate to within six digits. 
Panel (a) Horizontal axis: $\gamma=p/n$. Vertical axis: $\epsilon=k/n$.
Curve: $\eps^*(\gamma)$.
Panel (b) Horizontal axis: $\delta=(n-p)/n$. Vertical axis: $k/p$.
Curve: $\rho(\delta; C)$.  Same data underly both panels. The left panel
uses variables that are natural for robustness experts; the right panel
uses variables showing the agreement with the pattern already seen in figure 1.
} \label{fig:hadamard}
\end{figure}

To clarify better the relationships we are trying to make 
in this paper we set the standard error $Z$ to zero, and 
consider only the 'wild' outliers in $W$.  
We quantify breakdown properties of the
$\ell_1$ estimator by a large computational experiment. We  vary
parameters $(k,n,p)$ creating a range of
situations. At each one, we
solve the $\ell_1$ fitting problem (1.1) and 
measure to how many digits $\hat{\beta}$ agrees with $\beta$; we
record $\hat{\beta}$ as {\it breaking down} due to outliers
when fewer than six digits agree\footnote{
In order to save computer time, the actual experiment
conducted used an asymptotic approximation,
asymptotic in the size of the amplitude of the Wild component.
In our experiment, we set Z to zero, and modified the definition
of breakdown; instead of declaring breakdown when the 
estimated beta was wrong by more than 5 standard errors,
we declared breakdown when the estimated beta was wrong
in the sixth digit.  
Because of a scale invariance and continuity enjoyed by $\ell^1$, 
the experiments here can be viewed as the limit of standard 
robust statistics with ordinary noise, as the size of the wild component 
increases relative to the size of the ordinary noise $Z$.}.
Panel (a) in figure \ref{fig:hadamard} shows the results of this experiment,
depicting the breakdown fraction.  

Evidently, there is
an abrupt change in behaviour at a certain critical fraction
$\eps^*$; this depends on $\gamma= p/n$.
Now $0 < \gamma < 1$, so there are here more observational units
than predictors.  
When there are many observation units per predictor,
i.e. $\eps^* \approx 1$, $\ell_1$ fitting can resist a large
fraction of outliers. When the model is almost saturated,
i.e. nearly one predictor per observation, $\eps^* \approx 0$,
it takes very little contamination to break down the $\ell_1$ estimator.
On figure \ref{fig:hadamard}(a) we overlay a theoretical curve 
$ \eps^*(\gamma; C) = (1-\gamma) \rho((n-p)/n; C)$
where $\rho(\cdot;Q)$ is derived from geometric combinatorics.  
Evidently the curve coincides
with the observed breakdown point of the $\ell_1$ estimator
in a designed experiment. Figure \ref{fig:hadamard}(b) depicts a transformation
of panel (a), into new axes, with variables $\rho = k/(n-p)$ and $\delta = (n-p)/n$.
The display looks now similar to Stodden's figure 1.

Interpretations:
\bitem
\item  Standard $\ell_1$ fitting in a standard designed
experiment (but with large $n$ and $p$, this time with $p < n$)
turns out to be  robust below a certain critical fraction of outliers,
at which point it breaks down.
\item The simulation results, properly calibrated,
closely match seemingly unrelated
phenomena in sparse linear modelling in the $p > n$ case.
\item This critical fraction matches a known phase transition in
geometric combinatorics.
\eitem

\subsection{Surprise 3: compressed sensing}

We now leave the field of data analysis
for the field of signal processing.

Since the days of Shannon, Nyquist, Whittaker and Kotelnikov,
the `sampling theorem' has helped engineers 
decide how much data need to be acquired 
in design of measurement equipment.  Consider, therefore,
the following imaging problem.   We wish to acquire a signal $x_0$
having $N$ entries.
Now suppose that only $k \ll N$ of those pixels are
actually nonzero -- we do not know which ones are nonzero,
or even that this is true.  There are $N$ degrees
of freedom here, since any of the $N$ pixel values
vary.

Consider making $n \ll N$ measurements
of a special kind.  We simply observe $n$ random
Fourier coefficents of $x$.
Here $n \ll N$ so that, although the image has
$N$ degrees of freedom, we make far fewer measurements.
Let $A$ be the linear operator that delivers the selected $n$
Fourier coefficients and let $y$ be the resulting measured coefficients.
We attempt to reconstruct by solving for the object $x_1$ with
smallest $\ell_1$ norm subject to agreeing with the
measurements $y$:
\begin{equation} \label{PeeOne}
\min \| x \|_1 \mbox{ subject to } y = A x.
\end{equation}

\begin{figure}
\begin{center}
\includegraphics[bb=76 215 543 691,width=3in,height=2.1in]{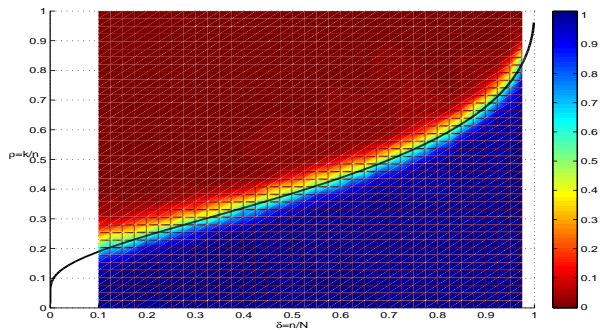} 
\end{center}
\caption{{\it Compressed Sensing from random Fourier measurements}.  
Shaded attribute: fraction of realizations in which 
$\ell_1$minimization  \eqref{PeeOne} reconstructs an image 
accurate to within six digits.
Horizontal axis:  undersampling fraction
$\delta = n/N$. Vertical axis: sparsity fraction  
$\rho = k/n$.} \label{fig:fourier}
\end{figure}

Figure \ref{fig:fourier} shows results of  computational experiments
conducted by the authors.  In those experiments we chose $N=200$
and varying levels of $k$ and $n$.  The horizontal axis $\delta = n/N$
measures the undersampling ratio -- how many fewer measurements
we are making than the customary $N$.  The vertical axis $\rho = k/n$
measures to what extent the effective number of degrees of freedom
$k$ is smaller than the number of measurements.  
The contours indicate the success rate.
The superimposed curve $\rho(\delta;C)$ roughly coincides with the 
empirical 50\% success rate curve.

Interpretation: We can violate the
usual `sampling theorem' ($n \geq N$) with impunity!
 The true limit is
$n \gtrsim k/\rho(n/N;C)$, where $\rho(\delta;C)$ is 
the curve decorating 
the display\footnote{The symbol $\gtrsim$ denotes
an asymptotic relatonship; for precise conditions see Donoho \& Tanner (2009$a$).}. 
We have seen this curve
twice before already;
it arises in a superficially unrelated problem in high-dimensional
geometric combinatorics.

\subsection{The connection to high-dimensional geometric combinatorics }

Recall the problem in geometric probability
we discussed in \S \ref{ssec-ConvGauss}.
Draw a sequence of $n$ samples $X_1, \dots , X_n$ 
from a standard $d$-dimensional normal
distribution. Let $P = Conv(X_1, \dots X_n)$
denote the convex hull of these
$n$ points; this is a random convex polytope.

Suppose that $d$ and $n$ are both large,
and let $\delta = d/n$.  Figure \ref{fig-9}
presents a black curve to be called
$\rho(\delta; T)$ and formally defined\footnote{The auxilary 
parameter $T$ in $\rho(\delta; T)$ 
associates this curve with $T^{N-1}$, the standard 
simplex   \eqref{eq:simplex};  see \S \ref{sec:GeoComb}.} in \S \ref{sec:GeoComb}.
It has the following interpretation.

Let $\rho = k/d$. Suppose that $\rho < \rho(d/n; T)$
and that $n$ and $d$ are both large, with $d < n$.
For the typical $k+1$ tuple $(X_{i_1}, \dots X_{i_{k+1}})$, 
\bitem
 \item every  $X_{i_j}$ is a vertex of $P$, $1 \leq j \leq k+1$;
 \item every line segment $[X_{i_j}, X_{i_{j'}}]$
 is an edge of $P$, $1 \leq j , j'\leq k+1$;
 \item ...
 \item the convex polytope $Conv(X_{i_1}, \dots X_{i_{k+1}})$
 is a $k$-dimensional face of $P$.
 \eitem
 
Figure \ref{fig-9} also presents a second, lower, black curve,
which is actually the one we have seen
in our three surprises.
This curve, denoted by $\rho(\delta; C)$ and defined
in \S \ref{sec:GeoComb},
has the following interpretation.  Draw the same 
$n$ samples $X_1, \dots , X_n$ 
from a standard $d$-dimensional normal
distribution. Now let  $Q = Conv(X_1, \dots ,X_n,-X_1, \dots, -X_n)$
denote the convex hull of the
$2n$ points including the original $n$ points
and their reflections through the origin.
This is a random {\it centrosymmetric} convex polytope.
 
 Let $\rho = k/d$ and $\eps > 0$. Suppose that $\rho < \rho(d/n; C)(1 -\eps)$
and that $n$ and $d$ are both large.
For the typical $k+1$ tuple $(X_{i_1}, \dots X_{i_{k+1}})$, 
\bitem
 \item every  $\pm X_{i_j}$ is a vertex of $Q$, $1 \leq j \leq k+1$;
 \item every line segment $[\pm X_{i_j},  \pm X_{i_{j'}}]$
 is an edge of $Q$, $1 \leq j , j'\leq k+1$;
 \item ...
 \item the typical convex polytope $Conv(\pm X_{i_1}, \dots \pm X_{i_{k+1}})$
 is a $k$-dimensional face of $Q$.
 \eitem

In short, the curve arising each time in Surprises 1-3
involves convex hulls of
symmetrized Gaussian point clouds.  
The curve involving $\rho(\delta ; T)$ would arise
if we had instead positivity constraints on the objects
to be recovered (Surprises 1 and 3) or on the outliers
(Surprise 2).

\subsection{This paper}

The curve $\rho( \delta;C)$
describes properties of high-dimensional polytopes
deriving from the Gaussian distribution. 
What  {\it now} seems surprising about
Surprises  1-3 is the lack of formal connection to those
polytopes in the applications.  For example
\bitem
 \item In Surprise 1, stepwise regression as practised by statisticians seems
           unrelated to convex polytopes.
 \eitem
In Surprises 2 and 3, the appearance of
 the $\ell_1$ norm establishes a connection with convex polytopes
 (as the unit ball of the $\ell_1$ norm is in fact a regular polytope).
 Yet
 \bitem
\item In Surprise 2, neither a Hadamard design nor outliers have any formal connection
           to any Gaussian distribution.
 \item In Surprise 3, observing random frequencies of the Fourier transform
           of a two-dimensional signal has no visible connection to any Gaussian
           distribution.
  \eitem   
  
The curves $\rho( \delta; C)$ and $\rho(\delta;T)$
accurately describe thresholds in many situations where the Gaussian distribution is
 not present; in fact we have witnessed it in cases where the points of a point cloud were 
 chosen deterministically.  We believe this signals a new 
 kind of  {\it limit theorem} in probability theory
 that, when formalized, will make precise a 
 new kind of universality phenomenon in high-dimensional geometry.
       
\section{Geometric combinatorics and phase transitions}\label{sec:GeoComb}

\subsection{Polytope terminology}

Let $P$ be a convex polytope in $\bR^N$, i.e.
the convex hull of points $p_1, \dots , p_m$. Let $A$ be
an $n \times N$ matrix.  The image $Q=AP$
lives in $\bR^n$; it is a convex set, in fact
a polytope, the convex hull of
points $Ap_1, \dots , Ap_m$.
$Q$ is the result of `projecting' $P$ from $\bR^N$ down to $\bR^n$
and will be called the {\it projected polytope}.

The polytopes $P$ and $Q$ are have vertices,
edges, $2$-dimensional faces, ... . Let $f_k(P)$ and
$f_k(Q)$ denote the the number of such $k$-dimensional
faces; thus $f_0(P)$ is the number of vertices of $P$
and $f_N(P)$ the number of facets, while $f_0(Q)$ is
the number of vertices and $f_n(Q)$ the number of facets.
Projection can only reduce the number of faces, so
\[
    f_k(Q) = f_k(AP) \leq f_k(P), \qquad k \geq 0.
\]

Three very special families of polytopes $P$
are available in every dimension $N>2$,
the so-called regular polytopes.  Here we consider 
two of the three:
\bitem
\item  the {\it simplex} (an $(N-1)$-dimensional analogue of the equilateral triangle)
\begin{equation}\label{eq:simplex}
T^{N-1}:=\left\{x\in\bR^N |\quad \sum_{i=1}^N x_i=1,\quad x_i\ge 0 \right\},
\end{equation}
and  the
\item  {\it cross-polytope}  (an $N$-dimensional analogue of the octahedron)
\begin{equation}\label{eq:crosspolytope}
C^{N}:=\left\{x\in\bR^N |\quad \sum_{i=1}^N |x_i| \leq1,\right\}
\end{equation}
\eitem
Statements concerning the hypercube similar to those made here for the
simplex and cross-polytope are available to an interested reader in
Donoho \& Tanner (2008$a$). 

\subsection{Connection to underdetermined systems of equations}
The regular polytopes are simple and beautiful objects,
but are not commonly thought to be {\it useful} objects.
However, their face counts reveal solution properties
of underdetermined systems of equations.
Such underdetermined systems arise frequently in
modern applications and the existence of {\it unique}
solutions to such systems is responsible for the three surprises
given in the introduction.  Consider the case of the simplex.

Consider the underdetermined system of equations $y_0 = Ax$, 
where $A$ is $n \times N$, $n < N$, and the optimization problem (LP):
\begin{equation}\tag{LP}
   \min 1'x \mbox{ subject to } y_0 = Ax, \quad x \geq 0 .
\end{equation}
Of course, ordinarily, the system has an infinite number
of solutions as does the problem \lp\ .

\begin{lemma}
Let $A$ be a fixed matrix with $n$ columns
in general position in $R^N$. Consider 
vectors $y_0$ with
a sparse solution $y_0 = A x_0$ where $x_0 \geq 0$ 
has $k$ nonzeros.
The fraction of systems $(y_0,A)$ where \lp\ has that underlying $x_0$ as its only solution is
\[
    fraction\{ \mbox{Exact Reconstruction using \lp\ } \} = \frac{f_k(AT^{N-1})}{f_k(T^{N-1})}.
\]
\end{lemma}

In short, the ratio of face counts between the projected simplex 
and the unprojected simplex tells us the probability that
the \lp\ correctly reconstructs a $k$-sparse $x_0$.

Consider the case of the cross-polytope.
Apply the optimization problem \p\  to
the problem instance $(y_0,A)$ generated by
\begin{equation}\tag{P1}
   \min  \|x\|_1 \mbox{ subject to } y_0 = Ax,
\end{equation}
an underdetermined system of equations $y_0 = Ax_0$, 
where $A$ is $n \times N$, $n < N$ . 
Of course, ordinarily, both the linear system and
the problem \p\ have an infinite number 
of solutions. 
\begin{lemma} \label{lem-C}
Let $A$ be a fixed matrix with $n$ columns
in general position in $R^N$. Consider 
vectors $y_0$ with
a sparse solution $y_0 = A x_0$ where $x_0$
has $k$ nonzeros.
The fraction of systems $(y_0,A)$  where \p\ has that 
underlying $x_0$ as its unique solution is
\[
    fraction\{ \mbox{Reconstruction using \p\ } \} = \frac{f_k(AC^N)}{f_k(C^N)}.
\]
\end{lemma}

In short, the ratio of face counts between the projected cross-polytope 
and the unprojected cross-polytope tells us
the probability that \p\ can successfully recover
a true underlying $k$-sparse object.




In short,  it is essential to know whether or not
\[
    \frac{f_k(AQ)}{f_k(Q)} \approx 1 
\]
for $Q$ the simplex or cross-polytope.

\subsection{Asymptotics of face counts with Gaussian matrices $A$}

We now consider the case where the $n$ by $N$ matrix $A$
has iid  Gaussian random entries.  Then the mapping $P \mapsto Q$
is a random projection. In this case, rather amazingly, tools
from polytope theory and probability theory
can be combined to study the  expected face counts
in high dimensions. The results demonstrate rigorously
the existence of sharp thresholds in face count ratios.
 
\begin{theorem}[Donoho (2005$a,b$), Donoho \& Tanner (2005$a,b$, 2009$a$)]\label{thm-1}
Let the $n \times N$ random matrix $A$ have iid $N(0,1)$ Gaussian elements.
Consider  sequences of triples  $(N,n,k)$ where $n = \delta N$ , $k = \rho n $,
and $N \goto \infty$. There  are functions $\rho(\delta;Q)$
for $Q  \in \{ T, C\}$ demarcating phase transitions
in face counts:
\[
      \lim_{N \goto \infty}  \frac{ f_k (AQ) }{f_k(Q)} = \left \{ 
         \begin{array}{ll}
             1 & \rho < \rho(\delta, Q) \\
             0 & \rho > \rho(\delta,Q) .
         \end{array} \right.
\]
\end{theorem} 

Figure \ref{fig-9} displays the two curves referred to in this theorem.
The simplex's  transition is higher than
the cross-polytope's:  $\rho(\delta,T)>\rho(\delta,C)$ for 
$\delta\in (0,1)$.


\section{Empirical results for non-Gaussian ensembles}

\subsection{Explorations}

Over the last few years we ran 
computer experiments generating
millions of underdetermined
systems of equations of various kinds, using standard optimization tools
to select specific solutions, and checking whether or not the solution 
was unique and/or sparse.  In overwhelmingly many cases, 
Gaussian polytope theory accurately matches the 
experimental results, even when the matrices
involved are not Gaussian. We here summarize results
about experiments with the non-Gaussian ensembles 
listed in table \ref{tab:ensembles}. Further detail is provided in 
the Electronic Materials Supplement (Donoho \& Tanner 2009$b$).
\begin{table}
\begin{tabular}{|l|l|l|r|}
\hline
Suite & Ensemble Name &  Coefficients & Matrix Ensemble \\
\hline
3 & Bernoulli & + & iid elements equally likely to be 0 or 1 \\
4 & Bernoulli & $\pm$ & iid elements equally likely to be 0 or 1 \\
\hline
\hline
5 & Fourier & +        & $n$ rows chosen at random from $N$ by $N$ DCT matrix \\
6 &Fourier & $\pm$ & $n$ rows chosen at random from $N$ by $N$ DCT matrix  \\
\hline
\hline
7 & Ternary (1/3)& +        & iid elements equally likely to be -1, 0 or 1 \\
8 & Ternary (1/3) & $\pm$ & iid elements equally likely to be -1, 0 or 1  \\
\hline
\hline
9 & Ternary (2/5)& +        & iid elements taking values  -1, 0 or 1  with $P(0)= 3/5$\\
10 & Ternary (2/5)& $\pm$ & iid elements taking values  -1, 0 or 1  with $P(0)= 3/5$  \\
\hline
\hline
11 & Ternary (1/10)& +        & iid elements taking values  -1, 0 or 1  with $P(0)= 9/10$\\
12 & Ternary (1/10)& $\pm$ & iid elements taking values  -1, 0 or 1  with $P(0)= 9/10$  \\
\hline
\hline
13 & Hadamard & +        & $n$ rows chosen at random from $N$ by $N$ Hadamard matrix \\
14 & Hadamard & $\pm$ & $n$ rows chosen at random from $N$ by $N$ Hadamard matrix  \\
\hline
\hline
15 & Expander & +        & special binary matrices, here with $P(0)=14/15$\\
16 & Expander & $\pm$ &special binary matrices, here with $P(0)=14/15$\\
\hline
\hline
19 & Rademacher & + & iid elements equally likely to be 0 or 1 \\
20 & Rademacher & $\pm$ & iid elements equally likely to be 0 or 1 \\
\hline
\end{tabular}
\caption{{\it Suites of problem instances.} Problem suite specifies random matrix ensemble and coefficient  sign pattern ($+$/$\pm$).
$Ternary (p)$ ensemble has $P(0)=1-p$ and $P(\pm1) = (1-p)/2$.  Expanders with $P(0)=p$ have
  distinct columns, each with fixed fraction $p$  of entries equal to one and
other entries zero.}\label{tab:ensembles}
\end{table}

\begin{figure}
\begin{center}
\begin{tabular}{c}
\includegraphics[bb=76 215 543 691,width=3in,height=2.1in]{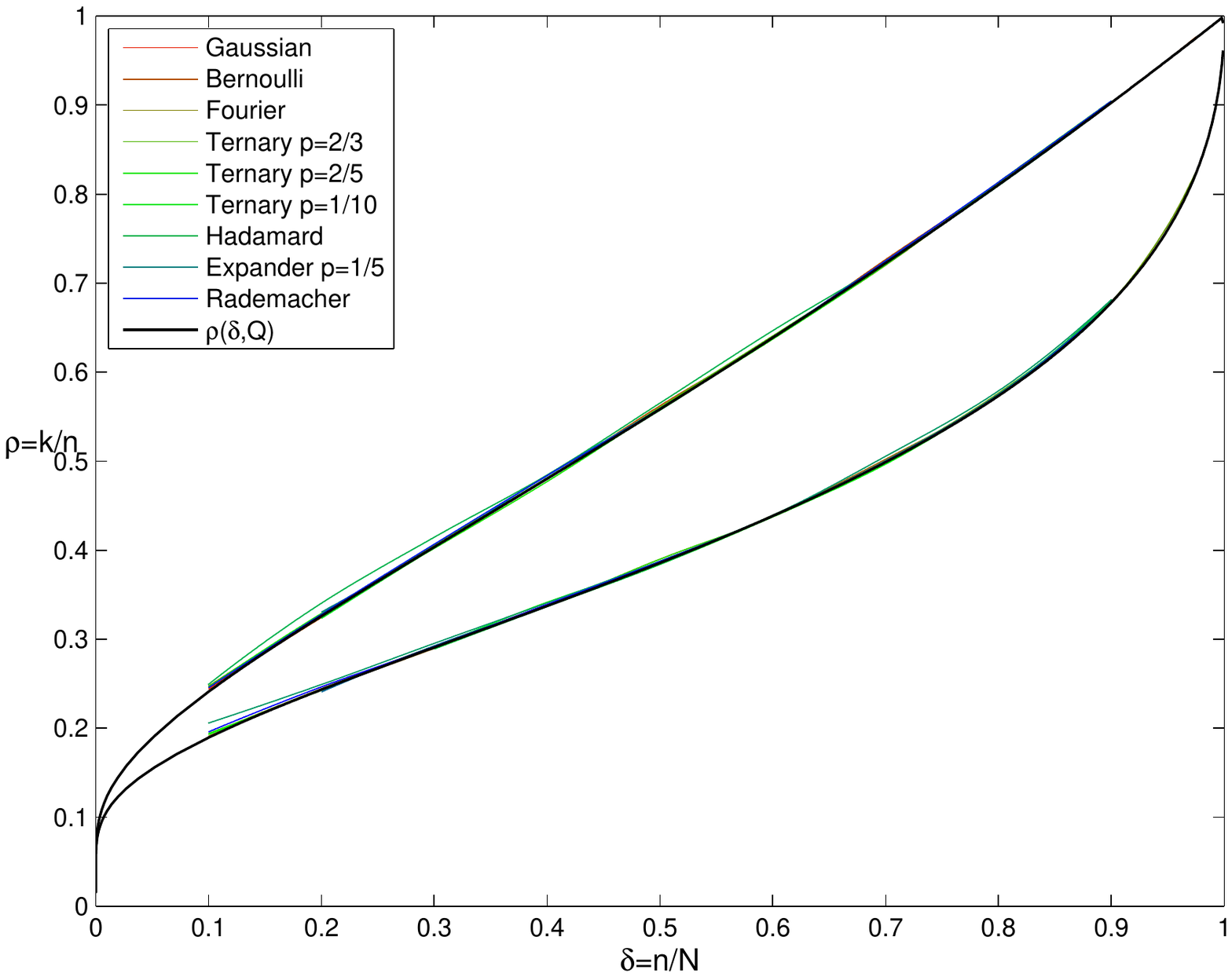} \\
\end{tabular}
\end{center}
\caption{Upper curves: level curves of 50\% success rate  for each 
  non-Gaussian suite in Table 1 with nonnegative coefficients, as well as for 
the Gaussian suite.  Asymptotic phase transition
$\rho(\delta;T)$ overlaid in black. Lower curves: level curves for 50\% success 
rate for each suite in Table \ref{tab:ensembles} with coefficients
of either sign, and for 
Gaussian suite. Asymptotic phase transition
$\rho(\delta;C)$ overlaid in black.}\label{fig-9}
\end{figure}

We varied
the matrix shape $\delta = n/N$, and the solution sparsity levels $\rho = k/n$.
 At problem size $N = 1600$ we varied $n$ systematically through a grid
 ranging from $n = 160$ up to $n = 1440$ in 9 equal steps.
 At each combination $N,n,k$ 
 we considered $M=200$ different problem instances
 $x_0$ and $A$, each one drawn randomly as above.
We both generated nonnegative sparse vectors
and solved \lp ,  and generated signed sparse vectors and solved \p .
The `signal processing language' event `exact reconstruction'
corresponds to the `polytope language'  event `specific $k$-face
of $Q$ is also a $k$-face of $AQ$'. In both cases we speak of {\it success},
and we call the frequency of success in  $M$ empirical trials at a given
$(k,n,N)$ the {\it success rate}.
At each combination $N,n$, we varied $k$ systematically to sample the success 
rate transition region from 5\% to 95\%.
Figure \ref{fig-9} presents summary results, showing the level curves
for 50\% success rate, for each of the 9 ensembles above.
The  appropriate theoretical curves  $\rho(\delta; Q)$ are overlaid.
The uppermost nine curves  give the case of nonnegative solutions, 
 $Q = T$, where we solve  \lp ; and the nine lower curves 
 present the data for $Q = C$, where we solved \p .
 (Note: the Hadamard case is exceptional and uses $N=512$.)
 
At  first glance, figure \ref{fig-9} shows excellent agreement between
the actual empirical results in each matrix ensemble\footnote{
Visual evidence, similar to figure \ref{fig-9}, of qualitative agreement was presented 
at conferences in 2006-2009 by Donoho and Tanner for all but the Expander 
ensemble.  Inclusion of the Expander ensemble in the results presented 
here was motivated by evidence in Bernide {\em et al.} (2008) 
for an Expander ensemble (with a different choice of $P(0)$) which also showing
qualitative agreement with the asymptotic phase transition $\rho(\delta;C)$.} 
and the asymptotic theory for the Gaussian.  This is not very surprising for 
$A$ from the Gaussian ensemble;
it merely proves that the large-$N$ polytope theory works
accurately already at moderate $N=1600$.
For the other ensembles there is not, to our knowledge, any existing theory 
suggesting what we see so clearly here:  {\it phase transition behaviour 
in non-Gaussian ensembles that 
accurately matches the Gaussian case} (compare \S \ref{sec:Conclusion}).

\subsection{Universality hypothesis}

Figure \ref{fig-9}  suggests to us the following

\begin{quotation}
{\bf Hypothesis.}  {\it Universality of Phase Transitions.}
Suppose that the $n$ by $N$ matrix $A$ is sampled
randomly from a ``well-behaved'' probability distribution. Suppose that
the $N$ by $1$ vector $x_0$ is sampled randomly from the set of $k$-sparse vectors,
either with or without positivity constraints on the nonzeros of $x_0$.
The observed behaviour of  solutions to $(LP)$ and $(P_1)$
will exhibit, as a function of $(N,n,k)$,
success probabilities matching those
which are proven to hold when sampling from the Gaussian distribution
with large $N$.
\end{quotation}

This hypothesis really contains two assertions: (a) that many matrix ensembles
behave like the Gaussian; and (b) that moderate-sized $N$ exhibit behaviour 
in line with the $N \goto \infty $ asymptotic.

The hypothesis also contains an element of  vagueness,
since we do not know at the time of writing how to delineate the
ensembles of random matrices over which Gaussian-like behaviour
will hold.  Of course universality results are well known
in probability theory; the Central Limit Theorem is the most
well-known universality result for the distribution
of sums of independent random variables.
The precise universality class of the 
Gaussian distribution  for such sums was only discovered two
centuries after the phenomenon itself was identified.
Apparently we are here at the stage of just identifying a comparable
phenomenon.  We hope it does not take two centuries to identify
the corresponding universality class!

{\it Clarification 1.}  In fact, our hypothesis could also be called a {\it rigidity} of the phase transition --
it is invariantly located at the same place in the phase diagram
across a range of matrix ensembles.   In statistical physics, universality of a phase transition
means something different, and much weaker -- not a rigidity, but instead a flexibility
of the location of a phase transition while preserving an underlying 
structural similarity.  Our hypothesis is far stronger.

{\it Clarification 2.}  In fact,  there are trivial counterexamples
to the hypothesis; for example the matrix $A$ of all ones does
not generate any useful phase transition behaviour.

\subsection{Experimental procedure}

We conducted a Monte Carlo experiment
to test the Universality Hypothesis.

The general procedure was like our earlier exploratory studies.
We call a {\it suite} a distribution of problem instances  $(y,A)$
fully specified by two factors: (1) the ensemble of matrices $A$
and (2) the ensemble of coefficients $x_0$ generating $y = Ax_0$.
Matrix ensembles include Bernoulli, Ternary, ...  Coefficient ensembles
studied here have vectors of $N$ coefficients
with only $k$ nonzeros, in sites chosen at random.  The positive sign
coefficient ensemble
indicated by $+$ has all nonzeros drawn uniformly from $[0,1]$. The signed
coefficient ensemble indicated by $\pm$ has nonzeros drawn uniformly from 
$[-1,1]$.  For each suite we 
visited a collection of triples $(N,n,k)$.
At each triple we drew a sequence of random problem
instances of the given size and shape 
from the given problem suite. We then
ran optimization software to compute the 
solution of the random problem.  We computed observables
from the obtained solution, in particular
the binary observable {\tt ExactRecon}, which takes the value 1 when
the obtained solution is equal to the
true solution within 6 digits accuracy,
and zero otherwise.

We aimed to be {\it confirmatory}
rather than {\it exploratory}:  to use formal inferential tools, 
and carefully explain apparent 
departures from our hypothesis.
Our experiments differed from earlier
efforts in scope and attention to detail.
\bitem 
\item {\it Scale.} We performed  2,948,000
separate optimizations spanning 16984 different
situations.  Our computations
required the use of as many as 200 CPUs
in an available cluster and overall required 6.8 CPU {\it years}.
We considered 16 problem suites based on 8 different matrix ensembles;
see table \ref{tab:ensembles}.
The scope of previous 
exploratory studies, which can be measured in CPU-days,
is tiny by comparison.
\item {\it Calibration.}  
The vast majority of our experimental
computations relied on Mosek, a commercial package. 
We made runs comparing the results with CVX,
a popular open-source optimization package.
We believe our results are consistent across optimizers.

\eitem

\subsection{Inferential formulation}
Rather than go on a fishing expedition, from the outset we
chose to frame our evaluation of the evidence using standard inferential 
procedures.
\bitem
\item  {\it Two-sample comparisons.} The strict form of
the Universality Hypothesis  says that the probability of unique solution under
the Gaussian Ensemble is the same as the probability
of unique solution at each other ensemble in the
universality class.  It follows that we may compare two sets 
of results at the same problem
size, one with the Gaussian ensemble and one where everything 
else in the problem is the same
except that a specific non-Gaussian ensemble is used.
If at each ensemble we generate $M$ problem instances 
and obtain $M$ realizations
of the observable {\tt ExactRecon}, strict universality requires
that the number of successes in each ensemble have a
binomial probability distribution with the same success probability in both ensembles.
Hence, the hypothesis really amounts to the assertion that two
binomial distributions are the same.  We proceed
with traditional tests for equality of two
binomial distributions.  We chose to work with the $Z$-score:
\begin{equation} \label{z-score}
     Z( \hat{p}_0, \hat{p}_1; M) =
     \frac{\hat{p}_0 - \hat{p}_1}{\hat{SD}(\hat{p}_0 - \hat{p}_1; M_0,M_1)}.
\end{equation}
Here $\hat{p}_i$ denotes ``the fraction of cases where {\tt Exact Recon} = 1
in ensemble $i$'', and 
$\hat{SD}(\hat{p}_0 - \hat{p}_1; M_0,M_1)$ is the appropriate 
standard error for comparing proportions with possibly unequal sample sizes $M_i$.
In this comparison $\hat{p}_0$ describes the Gaussian baseline experiment, and $\hat{p}_1$ 
describes the non-Gaussian alternative experiment.
Under the Universality
Hypothesis, $Z$ has an approximate standard normal distribution.

Reducing our problem merely to consideration of $Z$-scores we can
formalize  our hypothesis:
\begin{quotation}
{\bf Strong Null Hypothesis:} {\sl  The $Z$ scores have an approximate $N(0,1)$ 
distribution at each value of $(k,n,N)$.}
\end{quotation}
 
\item {\it Study of asymptotics with problem size.}
The Strong Null Hypothesis seems implausible {\it a priori}
on the strength of experience  from other settings.

Consider another setting where the Gaussian
distribution is universal: the  central limit theorem.
There, although  the Gaussian distribution provides the 
correct limiting behaviour, there are well-understood departures
from Gaussian behaviour at small problem sizes.
Such departures of course decay  with
increasing problem size. The
 theory of Edgeworth 
expansions  shows that such deviations
from Gaussianity decay with problem size
according to a specific  power of  size.
Hence, for a symmetric distribution, we will
see deviations of order $1/(\mbox{problem size})$
and, for an asymmetric one, deviations
of order $1/\sqrt{\mbox{problem size}}$ occur.

Analogously, in this setting
we may see systematic behaviour of the $Z$-scores
varying with problem size $N$ and perhaps
also with $k$.   We used three problem sizes \-- $N=200, 400$ and $1600$ \--
so we might identify trends 
in the $Z$-scores with problem size.

\begin{quotation}
{\bf Weak Null Hypothesis:} {\sl  The $Z$-scores  exhibit
discrepancies from the standard $N(0,1)$ distribution (e.g. in means, variances, tail probabilities)
which decay to zero  with increasing $N$.}
\end{quotation}

\eitem

\subsection{Results}

Results of our experiment were already summarized
in  figure \ref{fig-9}. For each suite in table \ref{tab:ensembles},
for each value of $\delta = n/N$, we measured the
value of $\rho = k/n$ at which the empirical probability of success
crossed 50\%.  Each of the 18 different curves in figure \ref{fig-9} presents 
results for one suite; each one depicts the 50\%
success rate curves as a function of $\delta$.
These ``empirical $\rho$ curves'' exhibit very strong visual agreement with
the corresponding theoretical curves $\rho(\delta;T)$ and $\rho(\delta;C)$.


\subsubsection{Raw $Z$-scores}
\begin{figure}
\begin{center}
\begin{tabular}{cc}
\includegraphics[bb=76 215 543 591,height=1.5in,width=2.3in]{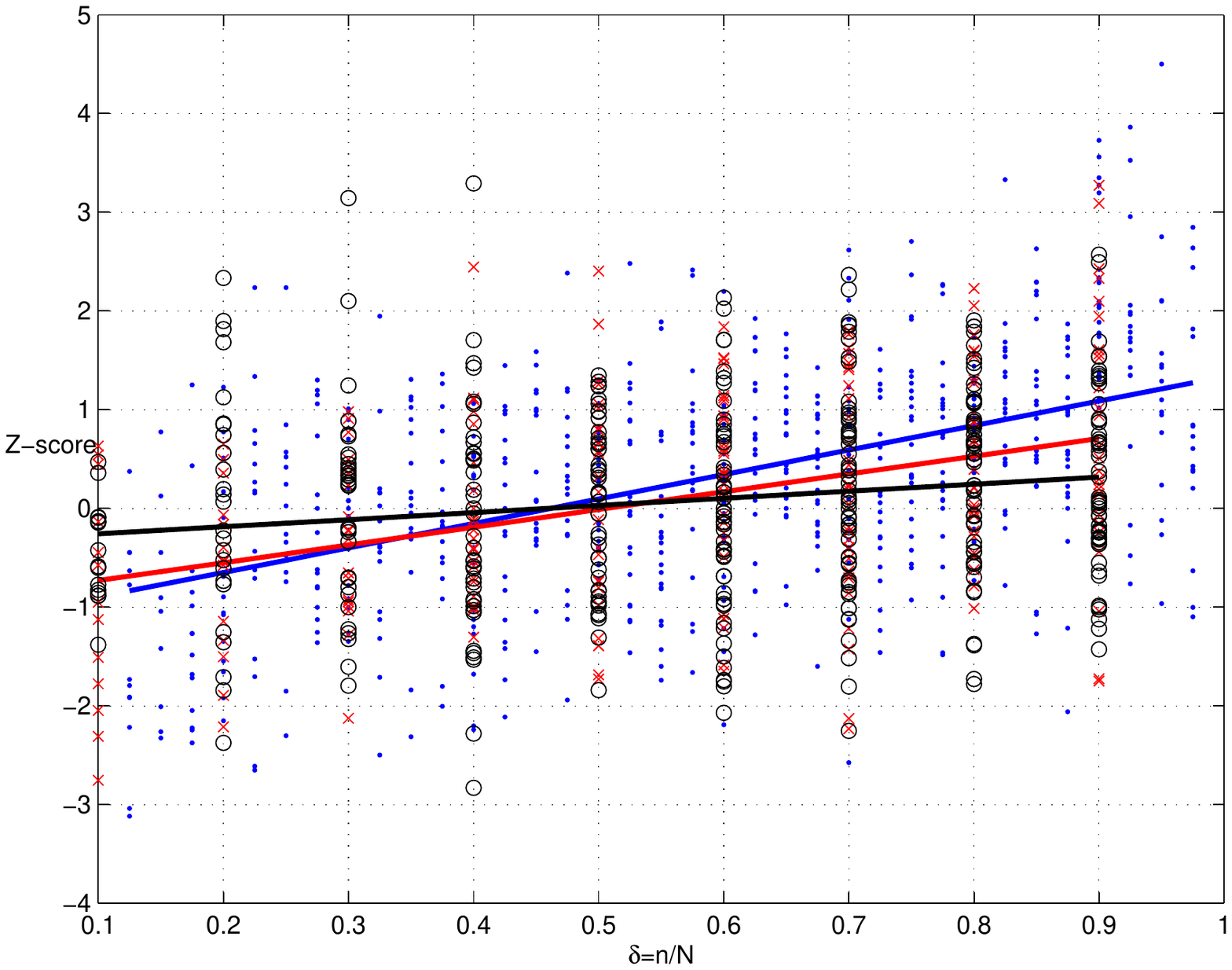} &
\includegraphics[bb=76 215 543 591,height=1.5in,width=2.3in]{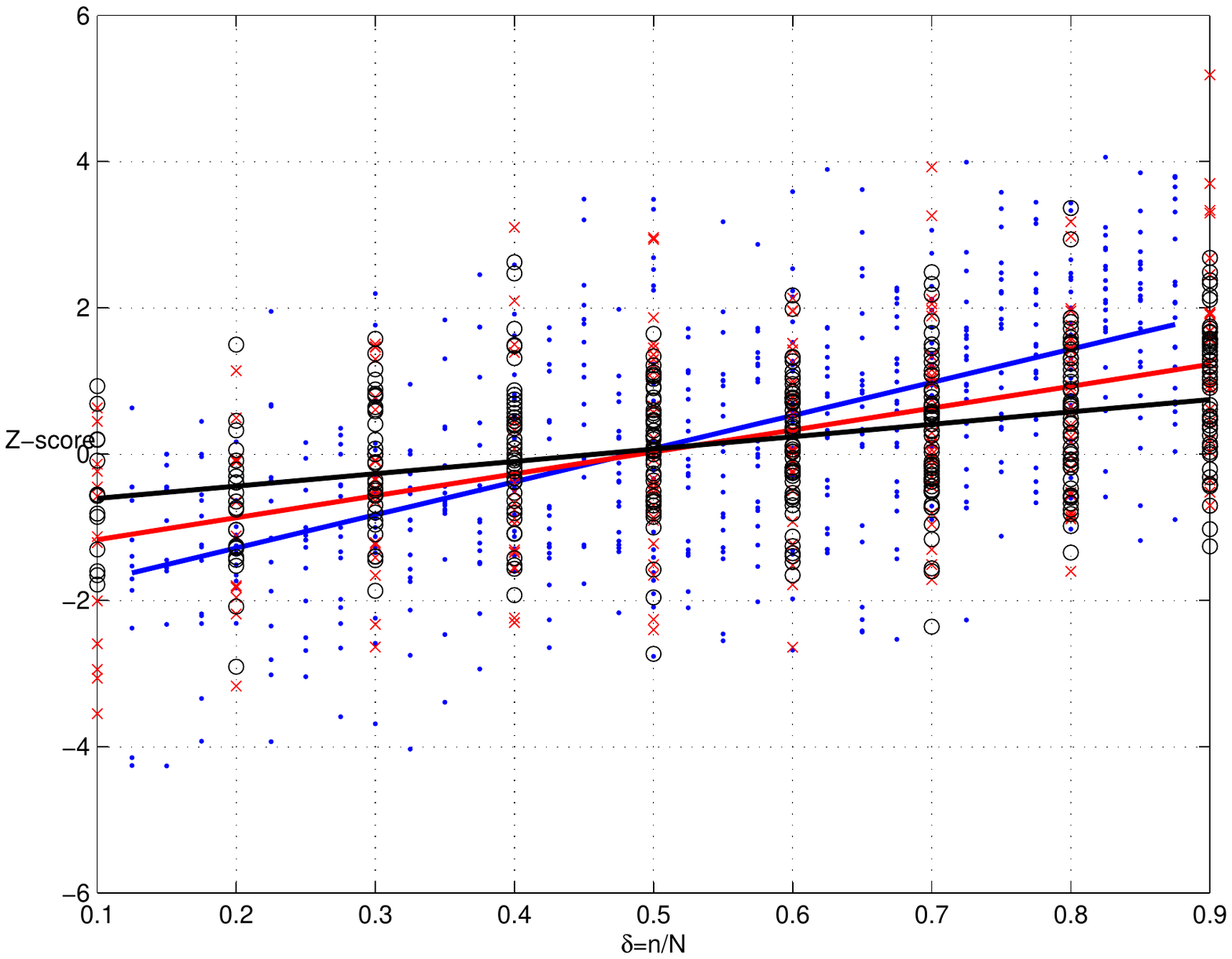} \\
(a) Bernoulli & (b) Fourier \\
\includegraphics[bb=76 215 543 591,height=1.5in,width=2.3in]{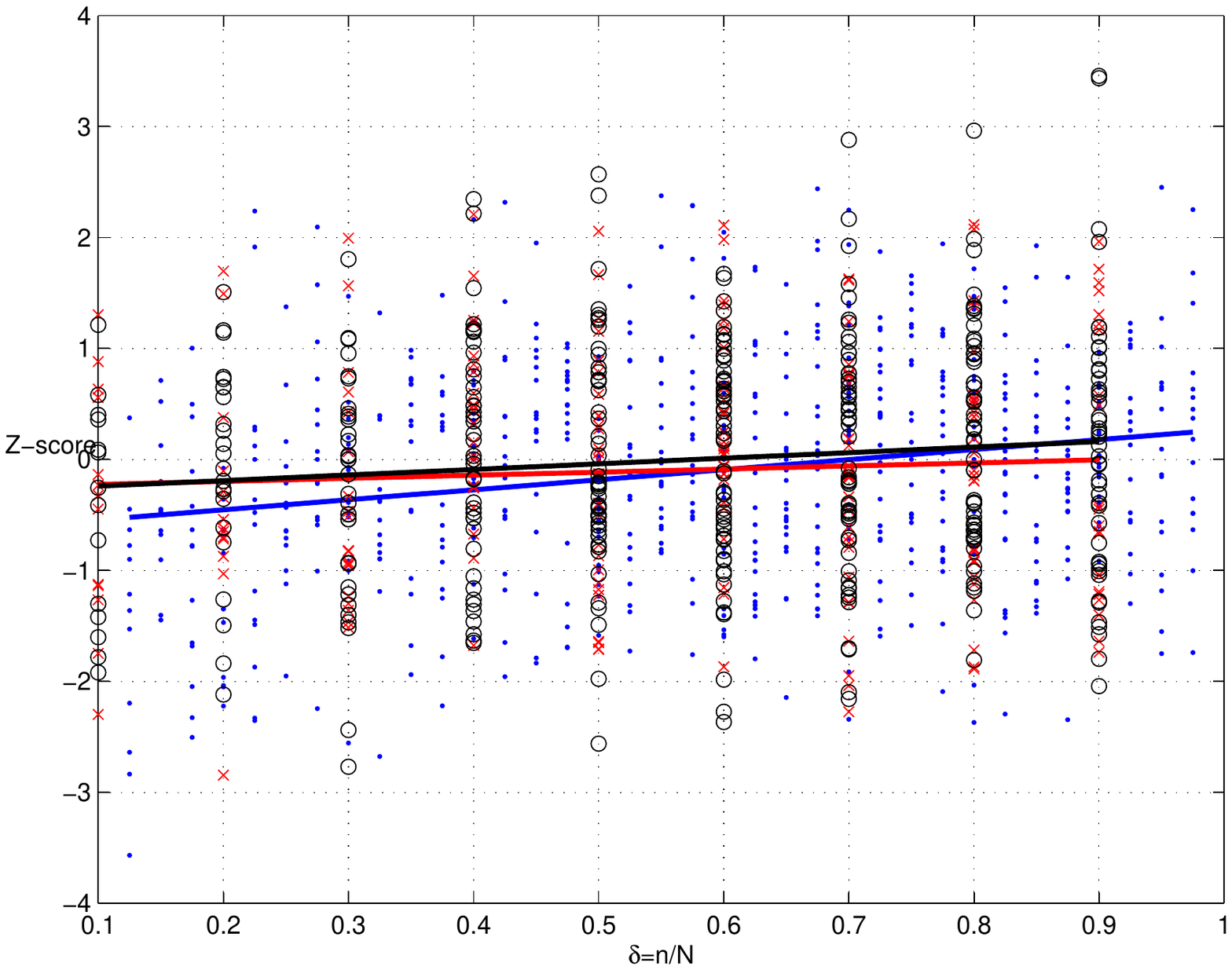} &
\includegraphics[bb=76 215 543 591,height=1.5in,width=2.3in]{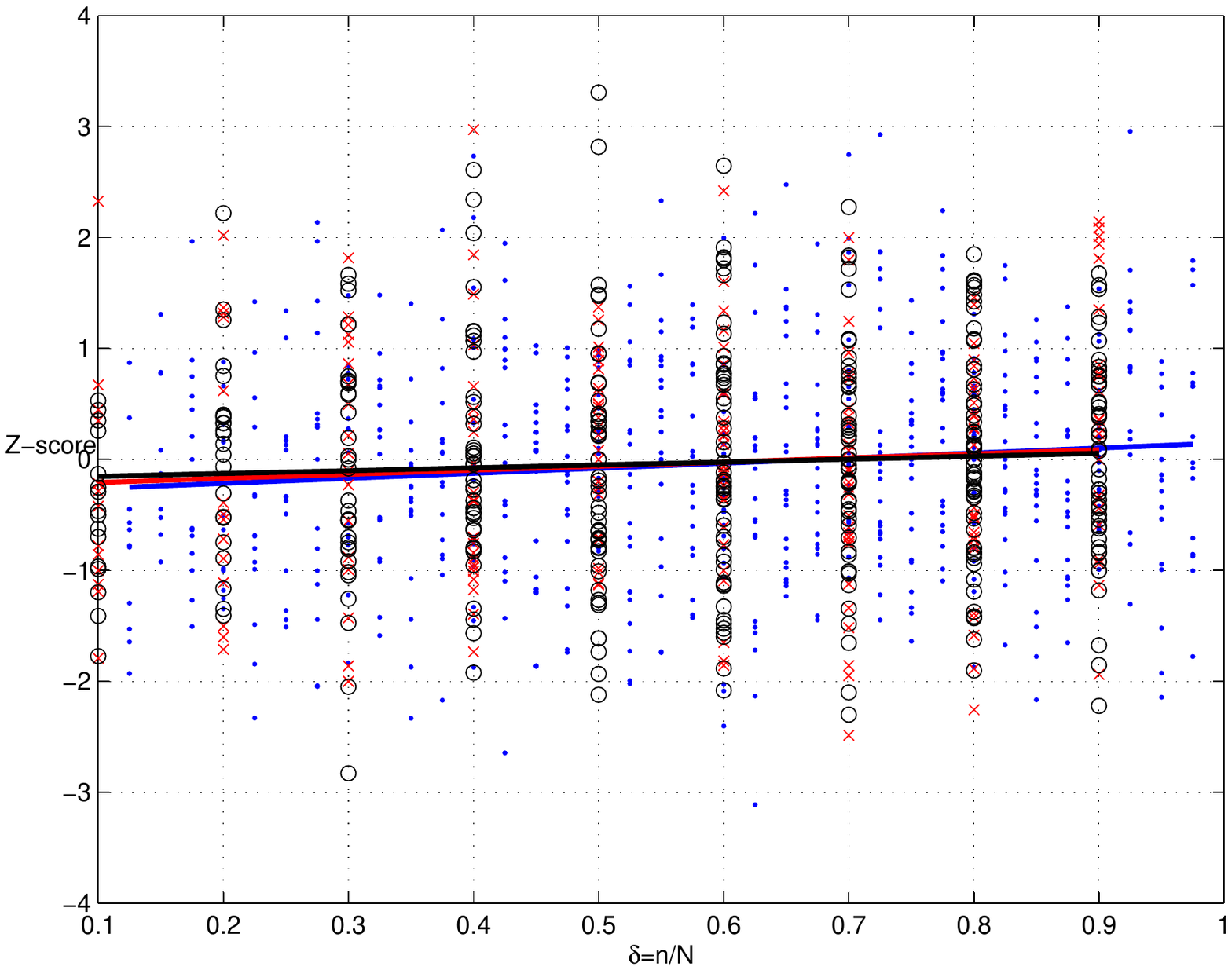} \\
(c) Ternary (1/3) & (d) Ternary (2/5) \\
\includegraphics[bb=76 215 543 591,height=1.5in,width=2.3in]{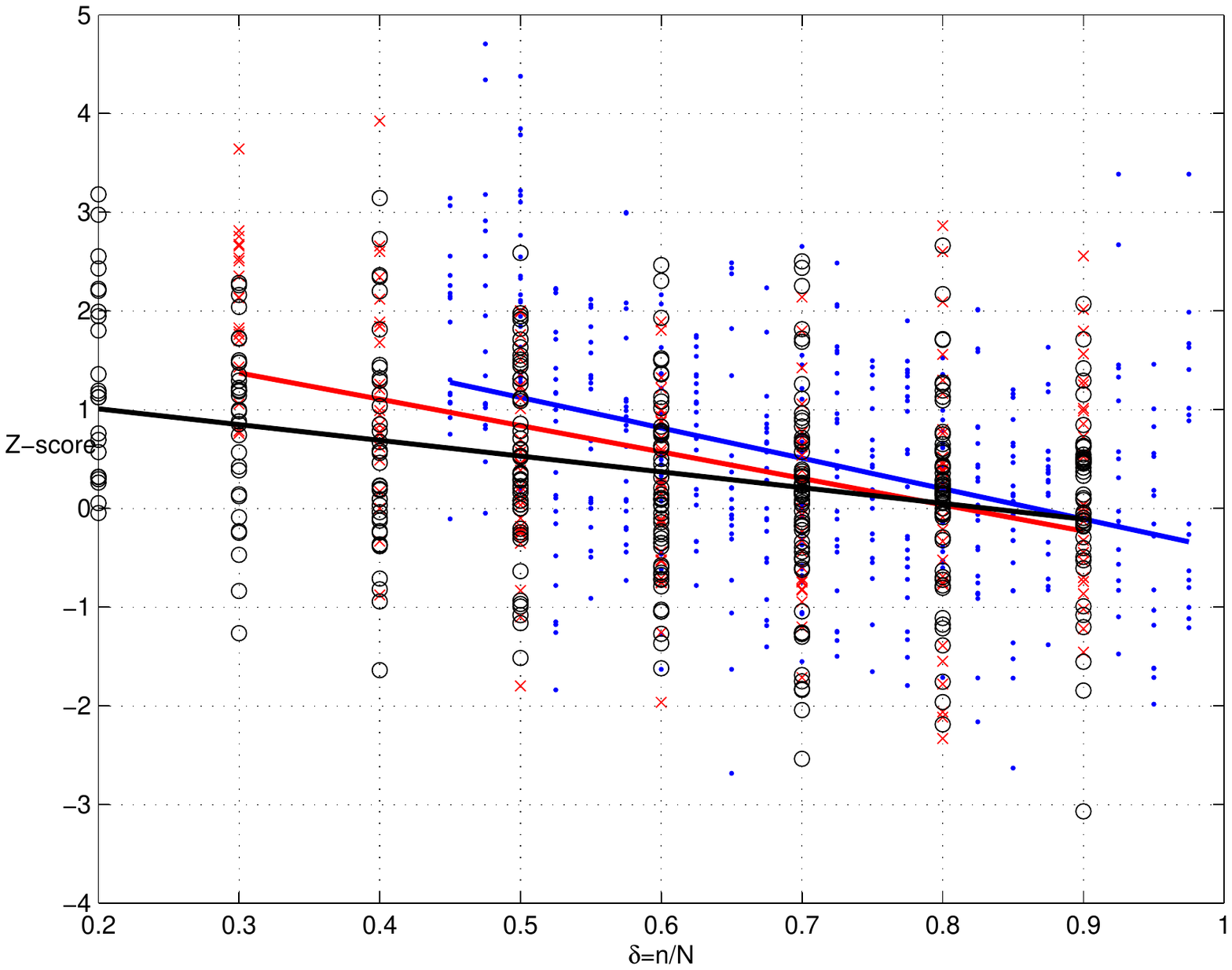} &
\includegraphics[bb=76 215 543 591,height=1.5in,width=2.3in]{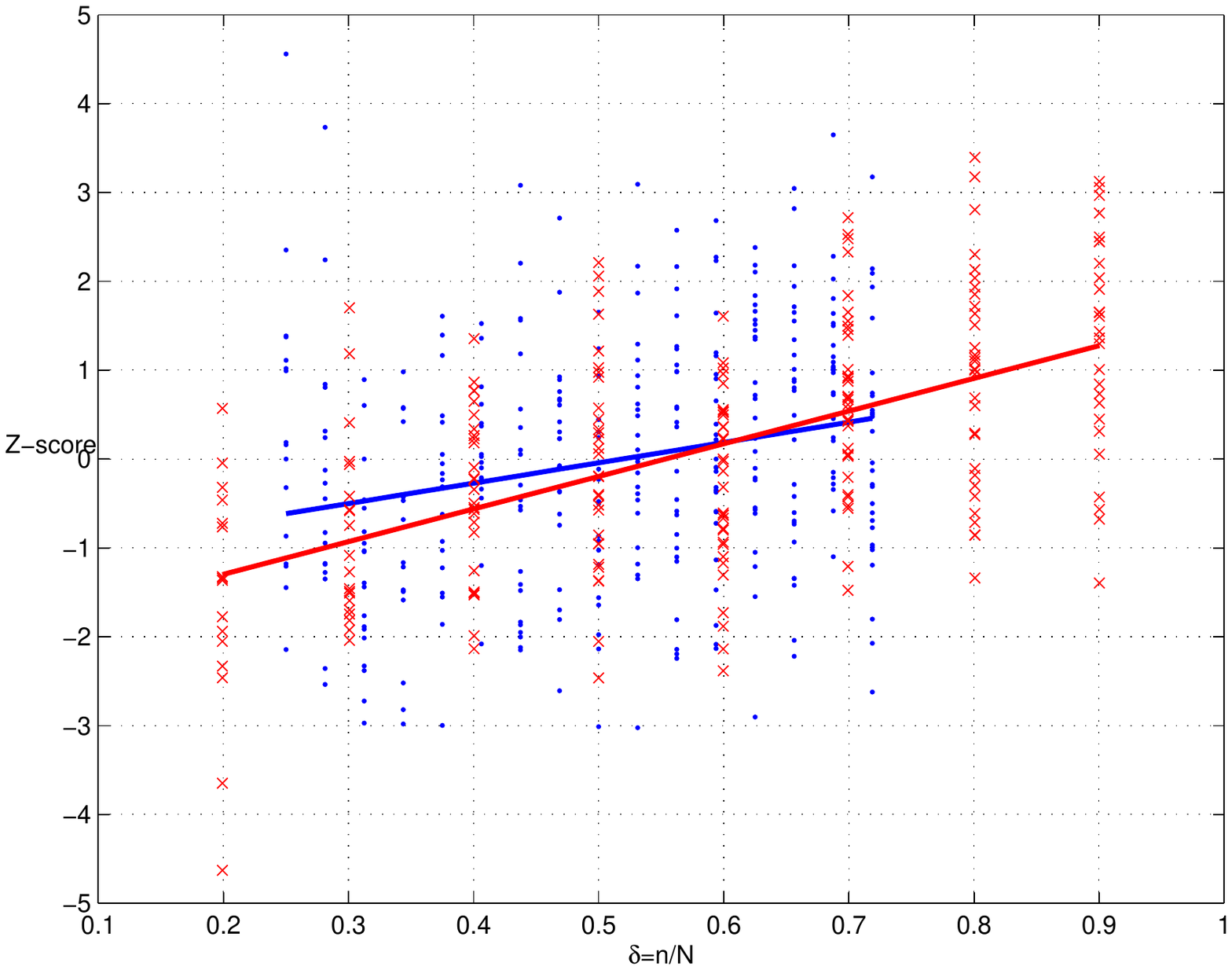} \\
(e) Ternary (1/10) & (f) Hadamard \\
\includegraphics[bb=76 215 543 591,height=1.5in,width=2.3in]{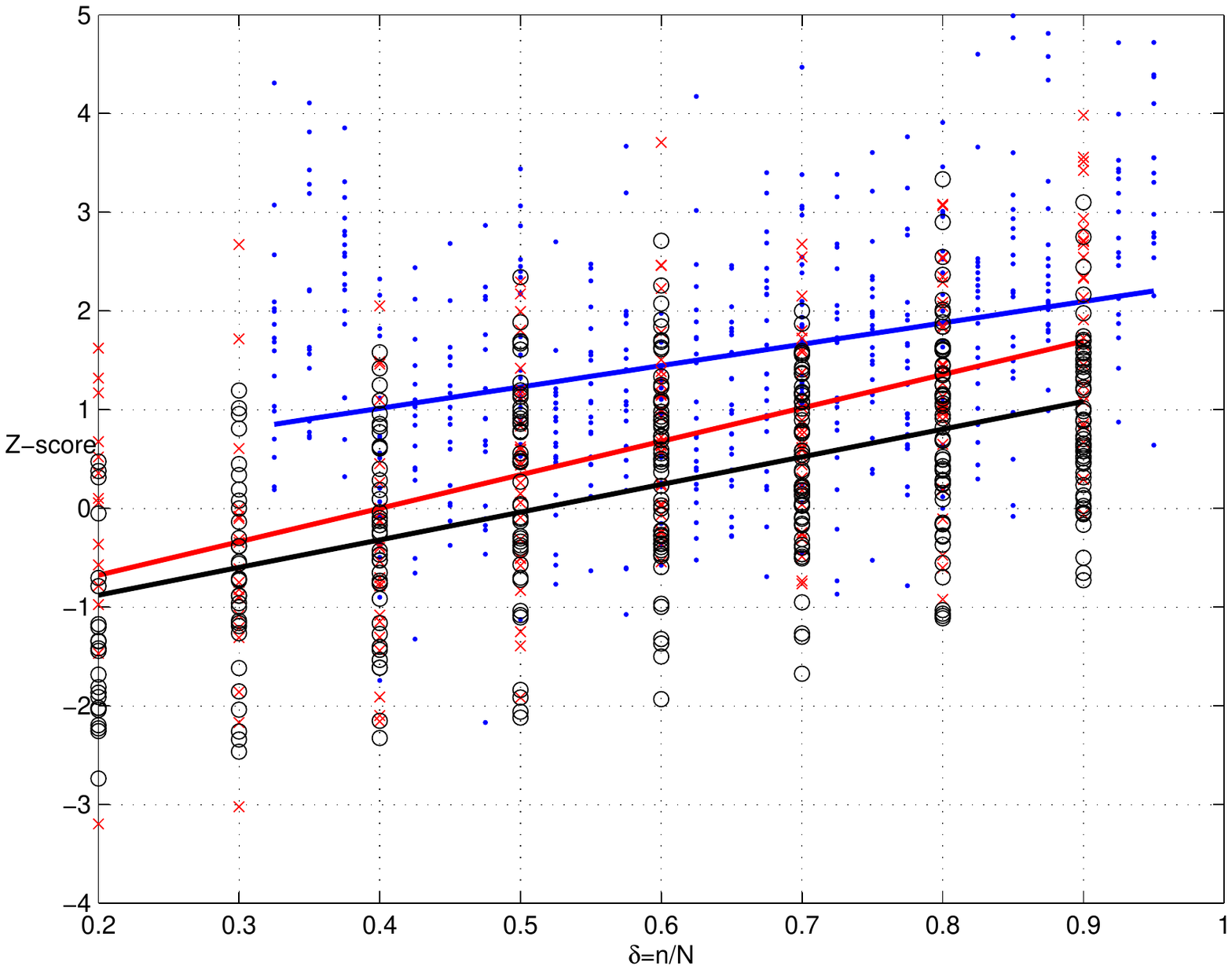} &
\includegraphics[bb=76 215 543 591,height=1.5in,width=2.3in]{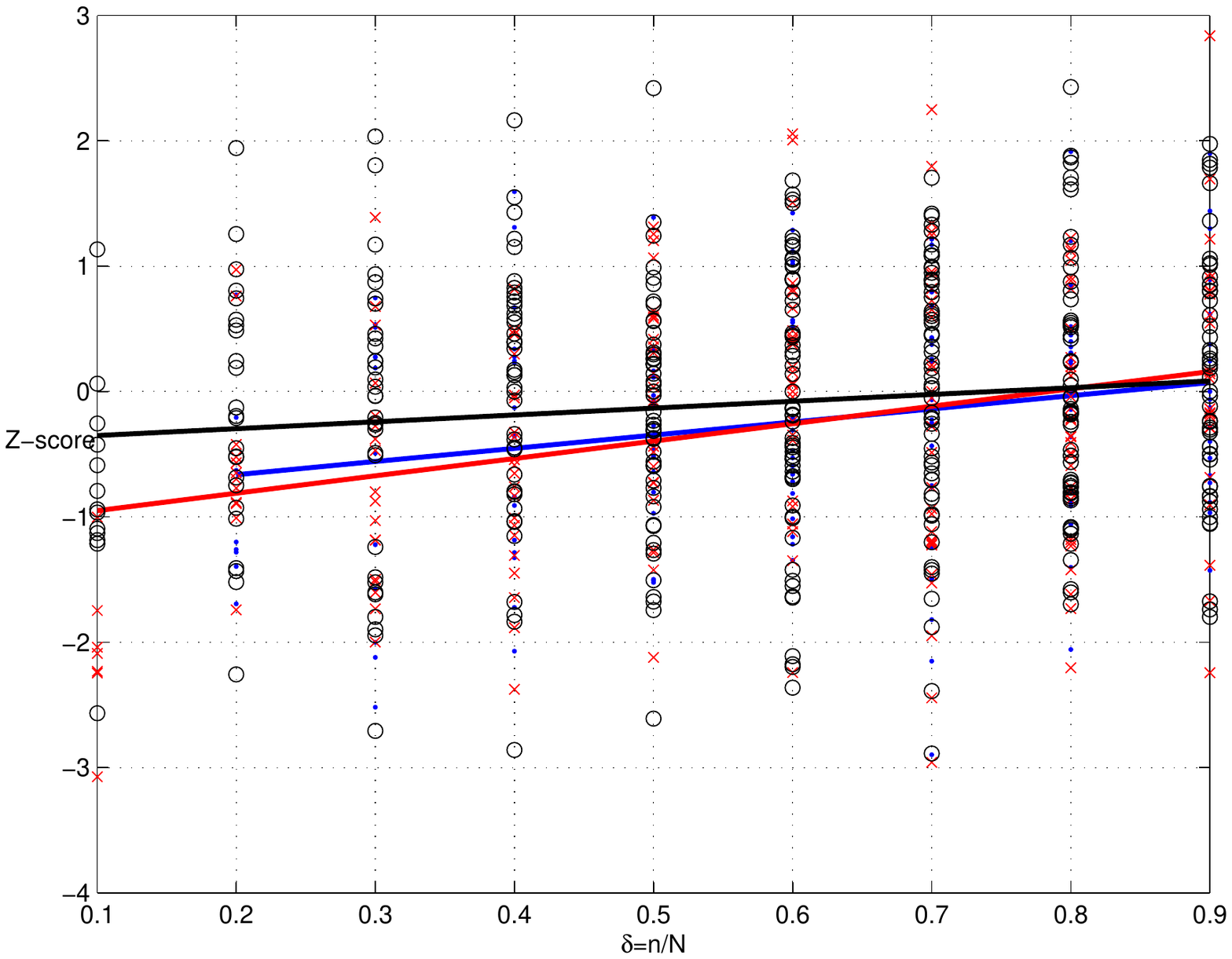} \\
(g) Expander & (h) Rademacher
\end{tabular}
\end{center}
\caption{{\it Raw $Z$-scores comparing success rate of
    reconstruction by \lp\ for $A$
    from the Gaussian ensemble vs. 
    success rate for $A$ from suites 3,5,7,9,11,13,15, and 19
    in Table \ref{tab:ensembles}.} Panels (a-e,g-h): suites
  3,5,7,9,11,15, and 19, respectively. Blue dots: $N=200$; red crosses:
   $N=400$; black circles: $N=1600$.   Linear fits to the raw 
$z$-scores in matching colors. Panel (f): Hadamard case,
dyadic $N$ only. Blue dots:  $N=256$; red crosses: $N=512$. 
} \label{fig:z1}
\end{figure}

\begin{figure}
\begin{center}
\begin{tabular}{cc}
\includegraphics[bb=76 215 543 591,height=1.5in,width=2.3in]{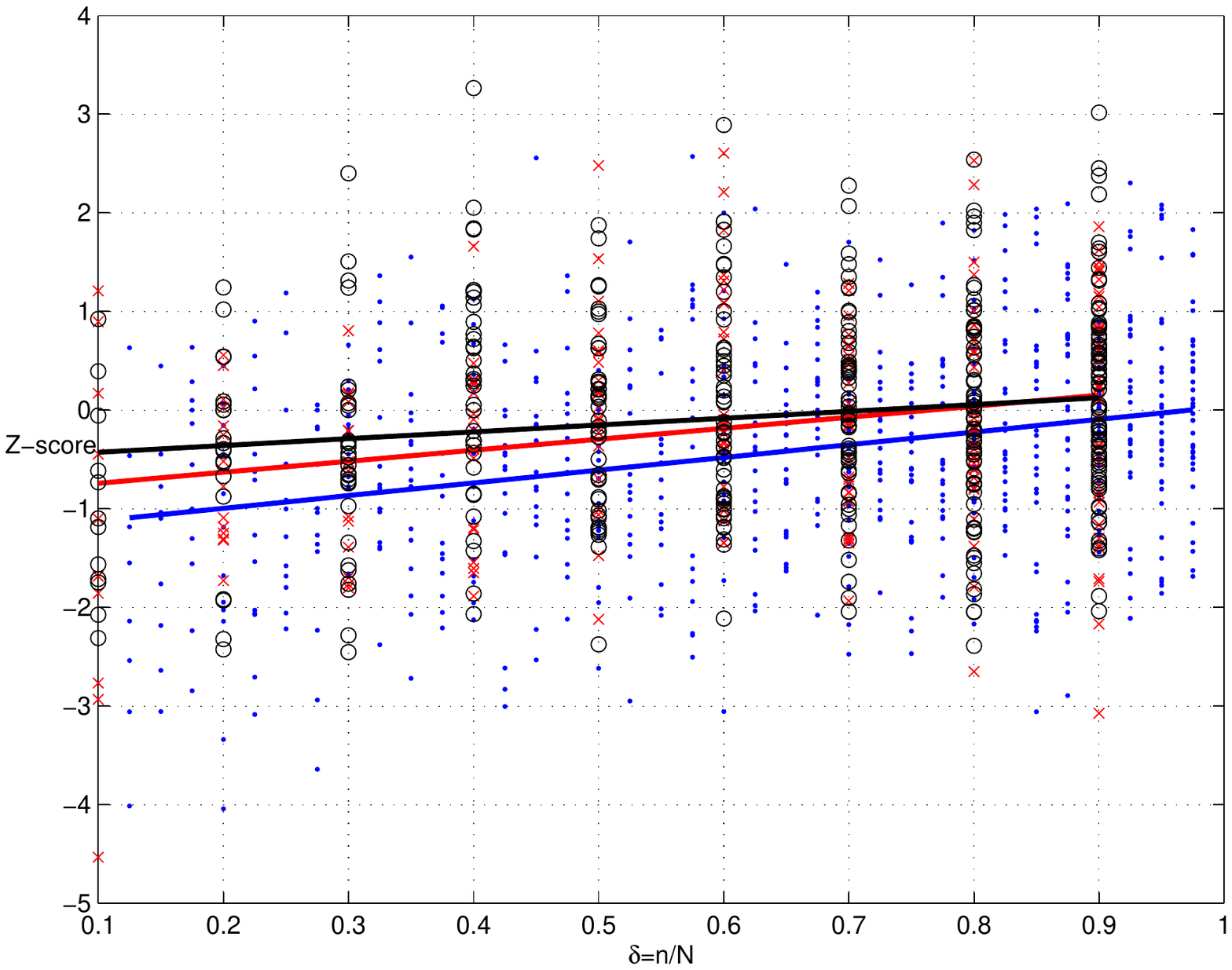} &
\includegraphics[bb=76 215 543 591,height=1.5in,width=2.3in]{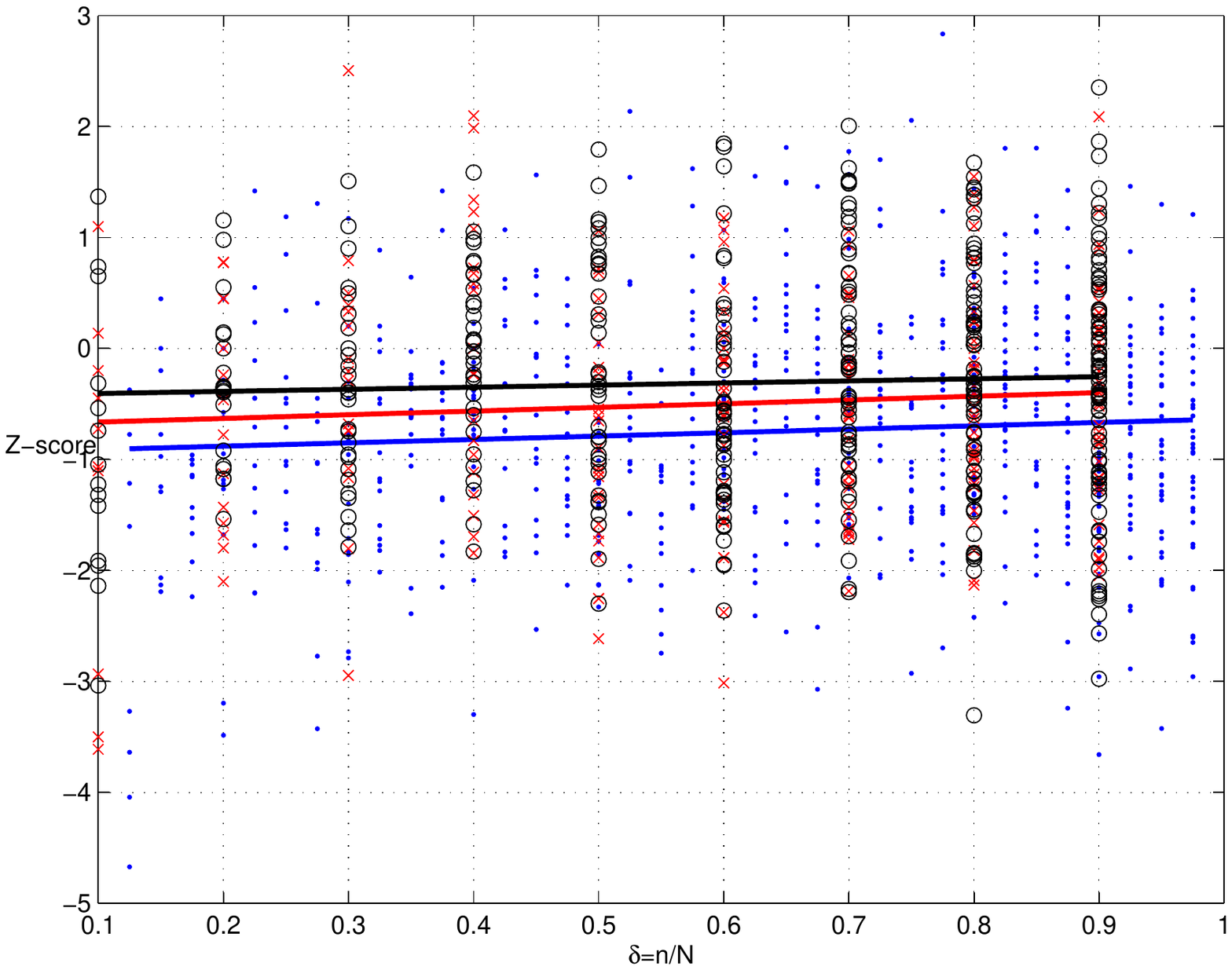} \\
(a) Bernoulli & (b) Fourier \\
\includegraphics[bb=76 215 543 591,height=1.5in,width=2.3in]{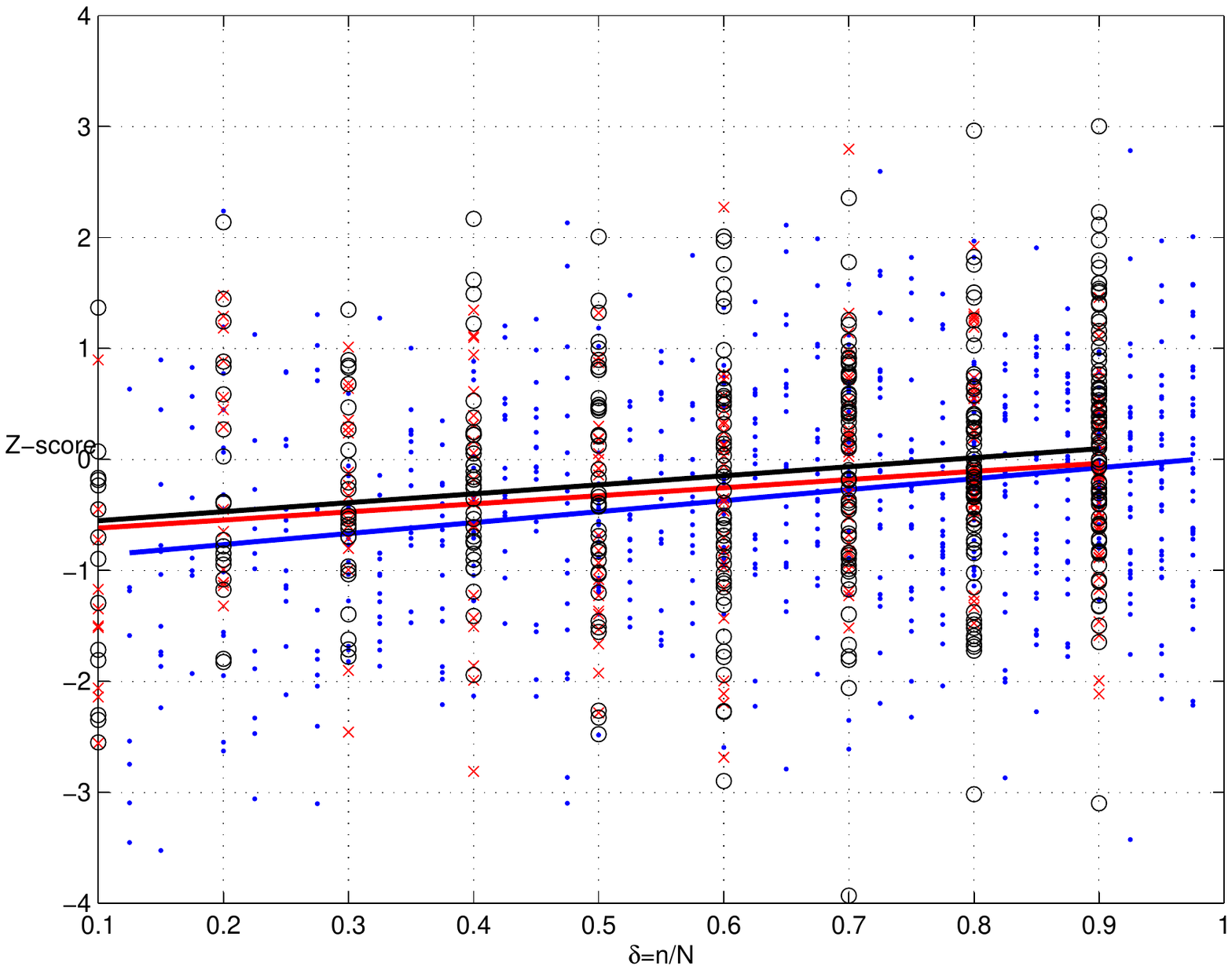} &
\includegraphics[bb=76 215 543 591,height=1.5in,width=2.3in]{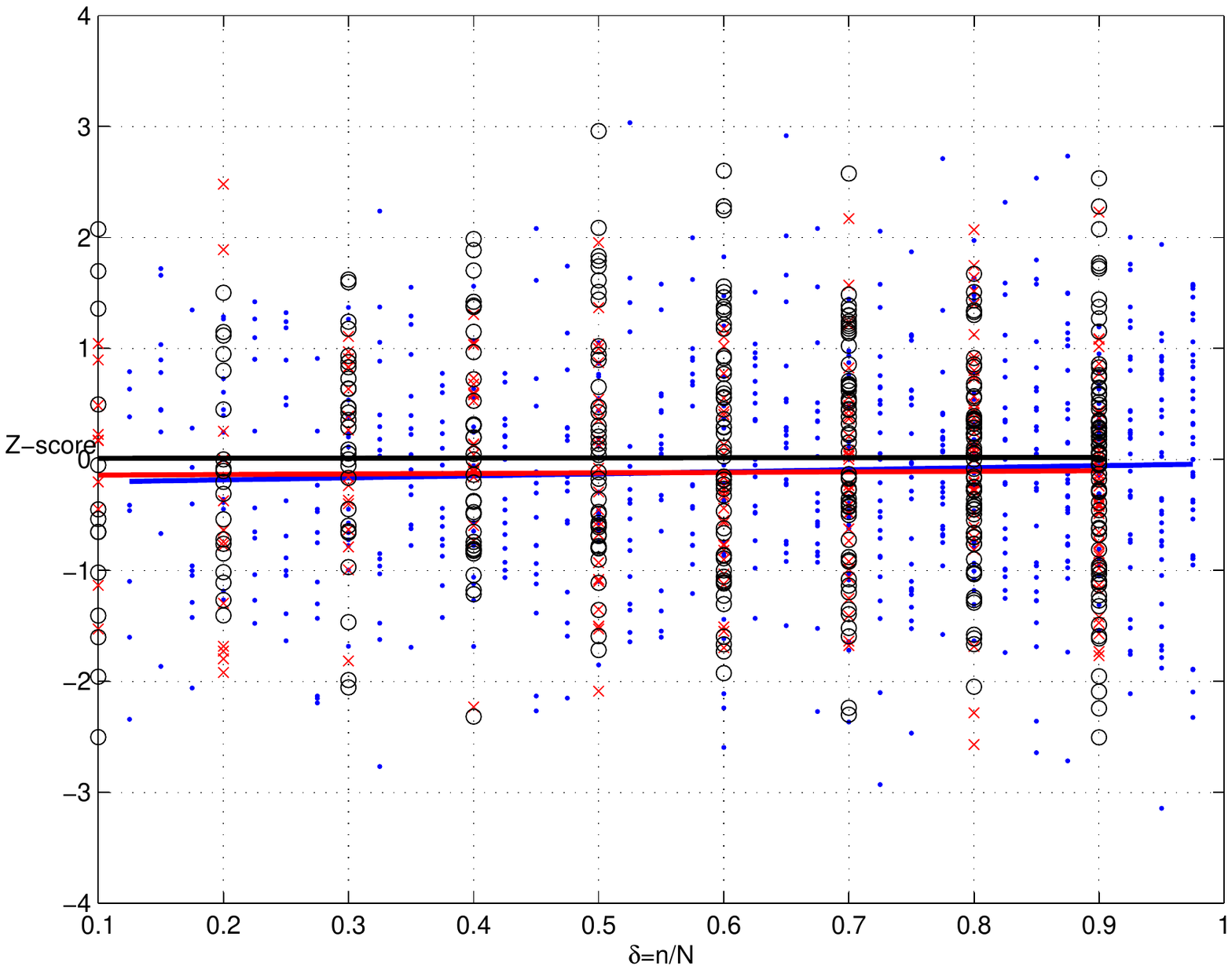} \\
(c) Ternary (1/3) & (d) Ternary (2/5) \\
\includegraphics[bb=76 215 543 591,height=1.5in,width=2.3in]{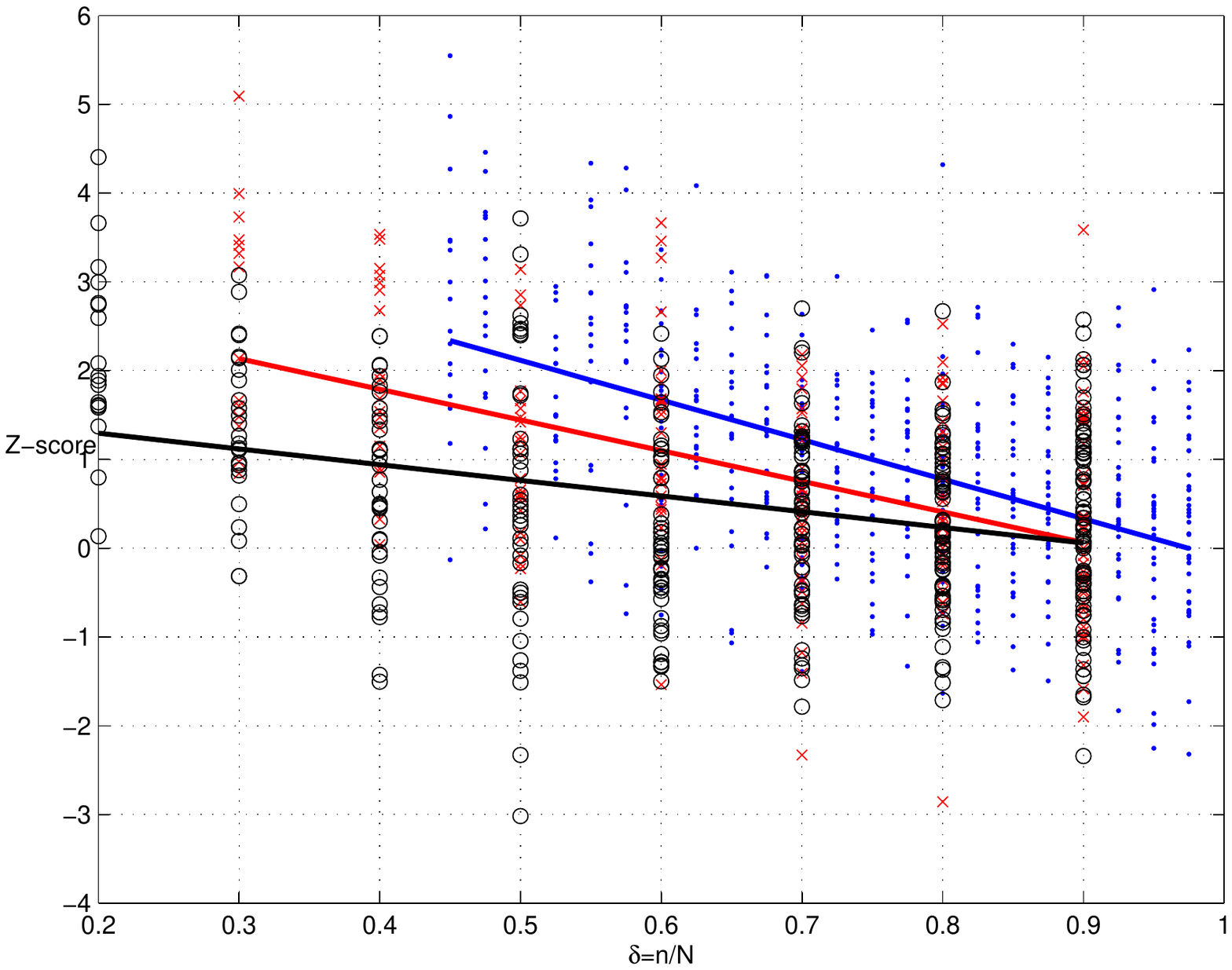} &
\includegraphics[bb=76 215 543 591,height=1.5in,width=2.3in]{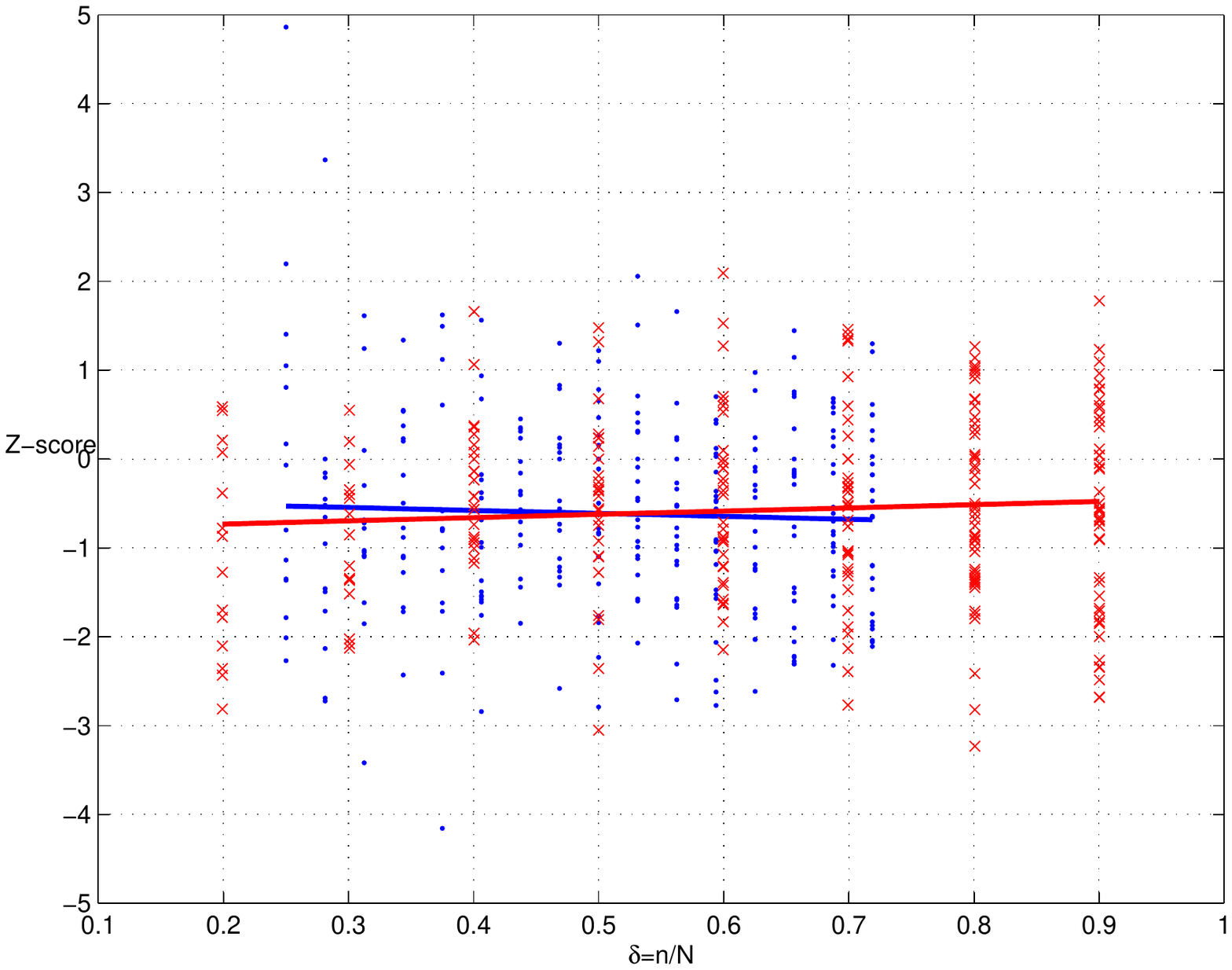} \\
(e) Ternary (1/10) & (f) Hadamard \\
\includegraphics[bb=76 215 543 591,height=1.5in,width=2.3in]{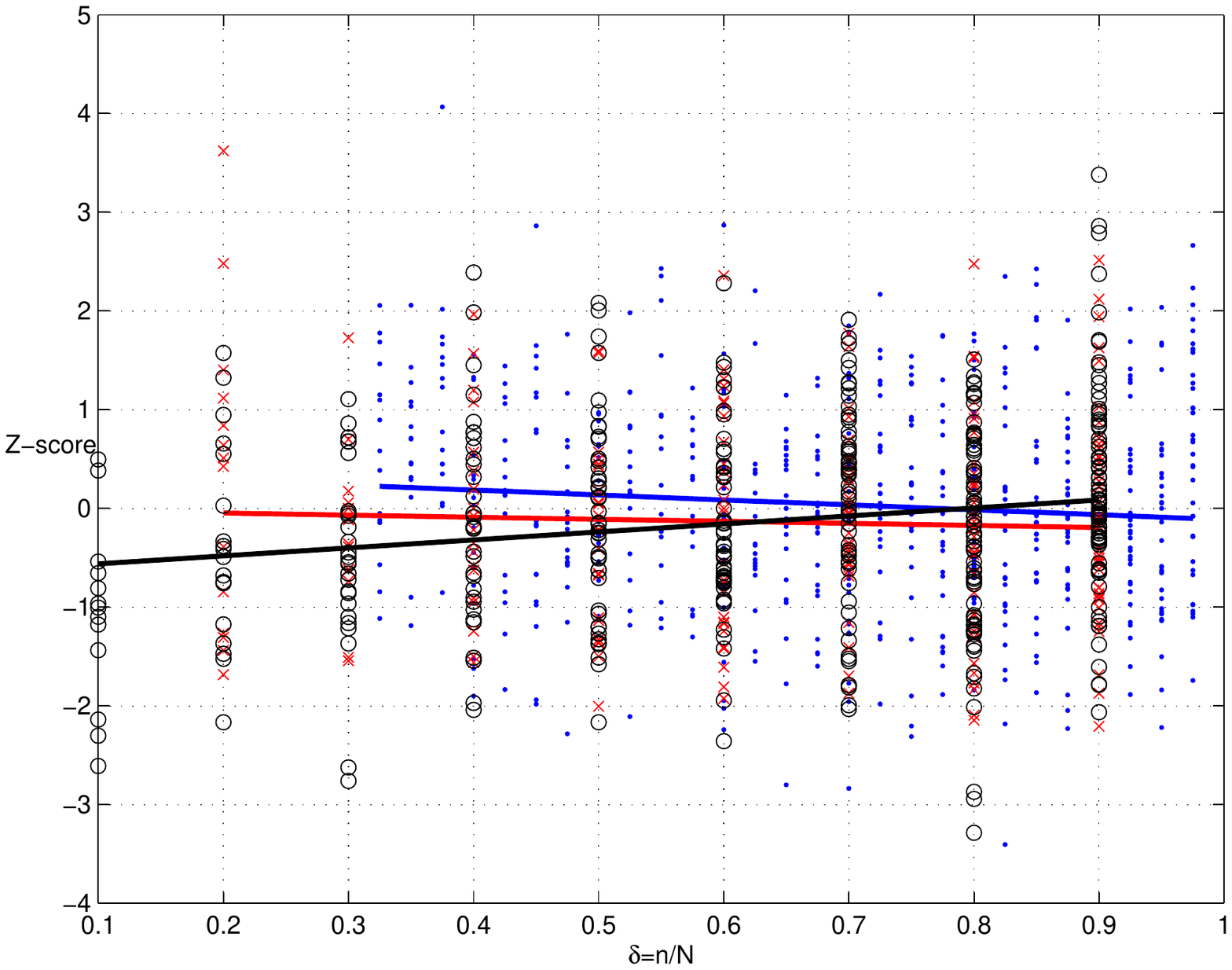} &
\includegraphics[bb=76 215 543 591,height=1.5in,width=2.3in]{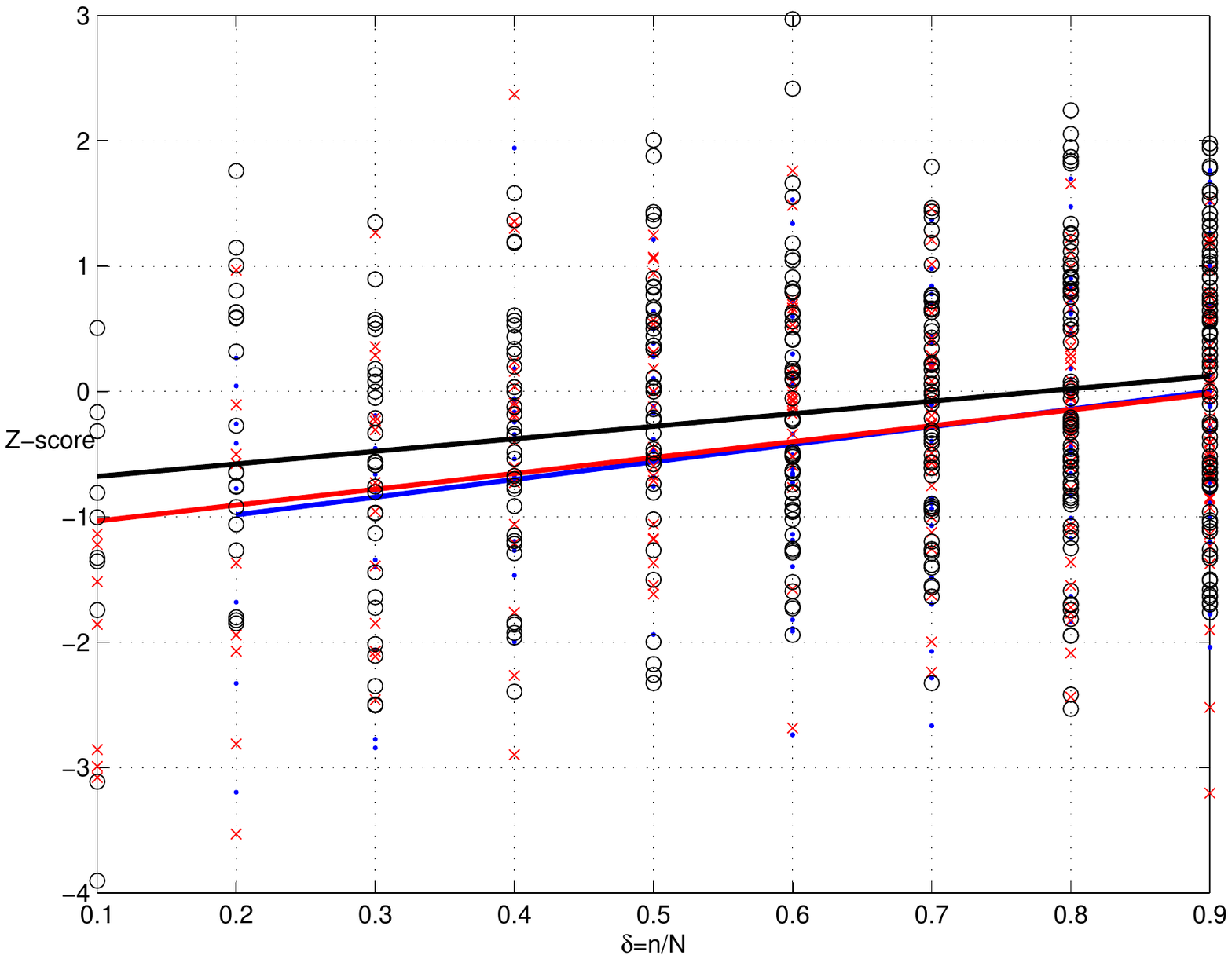} \\
(g) Expander & (h) Rademacher
\end{tabular}
\end{center} 
\caption{{\it Raw $Z$-scores comparing success rate in
    reconstruction by \p\ for $A$
    from the Gaussian ensemble vs. 
    success rate for $A$ from suites 4,6,8,10,12,14,16, and 20
    in Table \ref{tab:ensembles}.}  Panels (a-e,g-h): suites
  4,6,8,10,12,16, and 20 respectively. Blue dots: $N=200$; red crosses:
   $N=400$; black circles: $N=1600$.   Linear fits to the raw 
$z$-scores in matching colors. Panel (f): Hadamard case,
dyadic $N$ only.
Blue dots:  $N=256$; red crosses: $N=512$. 
} \label{fig:z2}
\end{figure}

Formal statistical tests are much more
sensitive and objective 
than visual impressions. Figure \ref{fig:z1}
displays the $Z$-scores (\ref{z-score}) for two-sample comparisons
between the Gaussian ensemble with nonnegative coefficients and 
the odd numbered suites; similarly, figure \ref{fig:z2} presents
two-sample comparisons between the Gaussian ensemble and the odd 
numbered suites.
The eight panels in each figure depict differences 
between each of the non-Gaussian ensembles
and the Gaussian ensemble.
The problem shape $\delta = n/N$ runs along the 
horizontal axis; these plots display results for $N=200,400,1600$
all combined.

The vast majority of the $Z$-scores in these displays fall in
the range $-2,2$.
\begin{quotation}
 {\bf Finding 1:} {\sl The $Z$-scores {\bf in bulk} are consistent with our hypothesis
of {\bf no difference} between the distributions.} 
\end{quotation}
In effect, our experiment conducted  16,984 hypothesis tests
and found relatively few `significant differences' at the
individual test level.

\subsubsection{Rejection of strict universality}
\newcommand{\cZ}{{\cal Z}}
There are marked  `tilts'
in the display of $Z$-scores in figures
\ref{fig:z1} and \ref{fig:z2};
linear trends with $\delta$ are visually
evident. Consider the general mean-shift model
\[
    Z(\delta,\rho; N, E) =  \mu(\delta; N, E) + {\cal Z} ,
\]
where $\cZ\sim N(0,1)$ is standard normal.  This expresses the 
idea that the observed $Z$-scores exhibit `drift'
as a function of $\delta$ and $N$,
but otherwise have the expected statistical properties of such scores.

If $\mu$ is not truly zero in this model, then of course the
null hypothesis of no difference fails.
In our setting this means that the Gaussian ensemble 
does not give truly the same success probabilities as the ensemble 
being compared to it.  Our analysis below rejects the hypothesis that $\mu=0$:
\begin{quotation}
 {\bf Finding 2:} {\sl The $Z$-scores are {\bf not consistent} with the hypothesis
of {\bf Strict Universality}.  Exact finite-problem-size agreement of 
success probabilities $p_1$ and $p_0$ between each alternative
ensemble and the Gaussian ensemble is not supported by our experiments.}
\end{quotation}

\subsubsection{Non-rejection of weak universality}

We reformulate the weak universality hypothesis
in terms of moments:

\begin{quotation}
{\bf Refined Null Hypothesis:} {\sl for each ensemble $E$, $\mu(\delta; N,E)$
tends to zero with increasing $N$, and the standard deviation of $Z$
scores approaches $1$.}
\end{quotation}

This hypothesis has a clear motivation.
Earlier displays aided the eye with lines
fitted to the means.
Evidently, the mean $Z$-scores within
a given suite are generally closer to zero for 
$N$ large than for $N$ small.  Hence,
in an informal appraisal,  the refined
null hypothesis seems quite plausible.  
Inspired by the `Higher Criticism'  (see Donoho \& Jin (2009)),
we compare the observed bulk distribution of $Z$ scores
with the theoretical distribution  $N(0,1)$. Table \ref{tab:HC}  shows 
roughly as many large $Z$ scores as one would expect under the null hypothesis.

\begin{table}
\begin{center}
\begin{tabular}{|l|l|l|l|l|}
\hline
$N$ &  \#$\{i: |Z_i| < 1\}$ & \#$\{i:|Z_i|<2\}$ & \#$\{i:|Z_i|<3\}$ & $M$ \\
\hline
 200 &101 / 99.67 & 142 / 139.4 & 148 / 145.6 & 148 \\ 
 400 & 146 / 140.6 & 199 / 196.6 & 208 / 205.4 & 208 \\ 
1600 &258 / 267.6 & 380 / 374.2 & 394 / 390.9 & 394 \\ 
\hline
\end{tabular}
\caption{Higher Criticism style analysis of residual $Z$-scores in suite 19.
Cell contents: observed counts/expected counts.}
\label{tab:HC}
\end{center}
\end{table}

\begin{quotation}
 {\bf Finding 3:} {\sl The $Z$-scores {\bf do not reject} the hypothesis
of {\bf Weak Universality}.  The difference between
success probabilities $p_1$ and $p_0$ for each alternative
ensemble and the Gaussian ensemble 
can be adequately modelled as a matrix-dependent random variable
with stochastic order $p_1(A) - p_0(A') = O_p(N^{-1/2})$,
where $A$ and $A'$ are realizations in the two matrix ensembles.}
\end{quotation}

\subsection{Results not presented in the main text}

In the appendix, 
we present a fuller record of our analyses.  Key points include the following.
\bitem
\item {\it Transition zone scaling with $N$}. We verified that
the width  $w(\delta,N; Q)$ of the zone where success probability
drops from 1 to 0 scales as $w\propto N^{-1/2}$.
\item {\it Adequacy of probit model}. We verified that
the  success rate varies with $\rho$ as a Probit
function $\overline{\Phi}( (\rho - \rho(\delta; Q))/w(\delta,N; Q) )$. Here 
$\overline{\Phi}$ is the Gaussian survival function and 
$w$ is the transition width.
\item {\it Exceptional Ensembles.} It is evident from 
figures \ref{fig:z1}-\ref{fig:z2}
that, at small $N= 200$, certain ensembles offer a
relatively poor match to the Gaussian case.
 Most of these
 discrepancies can be accounted for by saying that in these
 exceptional ensembles, at small problem sizes,
 the level curve for 50\% success
 rate is shifted noticeably below the 50\% curve for the Gaussian ensemble.
 However, at the larger problem size $N=1600$,
 both the shift and the exceptional character of the ensembles 
are no longer evident. For details, see the appendix.
\eitem
\subsection{Limitations of our conclusions}

We considered a limited set of matrix ensembles in this
study. The ensembles are not all based on iid elements --
there are dependencies among rows in the Fourier and Hadamard
ensembles  and among columns in the Expander ensemble.
Even so, there is a certain air of `orthogonality' or `weak independence'
in these examples. 

There are exceptions to the pattern presented
here. The classical example of cyclic polytopes shows that
\lp\ can have a notably higher success rate for very special matrices
than it does for random matrices (Donoho \& Tanner, 2005$a$).

In addition to forward stepwise regression, \lp\ and \p , there are 
several competing algorithms which we do not study here.
Maleki \& Donoho (2009) conducted extensive 
empirical testing of many such algorithms and observed  clear
phase-transition-like behaviour, which varies from algorithm to algorithm
and from the results presented here.
Unlike the phase transitions presented here, which match theoretical results in
combinatorial geometry, the phase transitions observed for competing algorithms
are not yet supported by theoretical derivations.

\section{Conclusion, and a glimpse beyond}
\label{sec:Conclusion}

Certain phase transitions in high-dimensional
combinatorial geometry have been
derived assuming a Gaussian distribution.
We had informally observed
that the Gaussian theory seemed
approximately right even in some non-Gaussian cases.
In this study, we made extensive computational
experiments with more than a dozen
matrix ensembles considering millions
of instances at a range of problem sizes.  
Empirical results for both Gaussian
and non-Gaussian ensembles show finite-$N$ transition
bands centred around the asymptotic phase transition derived
from a Gaussian assumption. The bands have a width
of size $O(N^{-1/2}$), consistent with the proven behaviour for 
the Gaussian ensemble (Donoho \& Tanner 2008$b$).
Such behaviour at non-Gaussian ensembles
goes far beyond current theory.
Adamczak {\em et al.} (2009) proved that, for a range
of random matrix ensembles with independent columns, there 
is a region in the phase diagram
where the expected success fraction tends to one, but
there is no suggestion that this region matches
the region for the Gaussian.

We used standard two-sample statistical
inference tools to compare results from  non-Gaussian
ensembles with their Gaussian counterparts at
the same problem size and sparsity level. We
observed fairly good agreement of the two-sample $Z$-scores
with the null hypothesis of no difference; however,
fitting a linear model to an array of such $Z$-scores
we were able to identify statistically significant trends of the $Z$-scores
with problem size and with undersampling fraction $\delta=n/N$.
The fitted trends vary from ensemble
to ensemble, decay with problem
size, and are consistent with weak, `asymptotic', universality
but not with strong,  finite-$N$, universality.

Our evidence points to
a new form of `high-dimensional limit theorem'.
There is some as-yet-unknown class
of matrix ensembles that yield phase transitions
at the same location
as the Gaussian polytope transitions.
Delineating this universality
class seems an important new task for future work in 
stochastic geometry.

\section{Appendix: suplementary statistical analysis}

We present details of the data analysis.

\vspace*{.1in}
\noindent
{\sl Gaussian ensemble.} We study the basic properties of success
probabilites at the Gaussian ensemble, as a function of $\delta$ and
$\rho$.
\bitem
\item {\it Transition zone scaling with $N$}. We quantify
the width  $w(\delta,N; Q)$ of the zone where success probability
drops from 1 to 0. We verify that our measurement scales as $N^{-1/2}$.
\item {\it Adequacy of Probit/Logit models}. We verify that
the  success probability varies with $\rho$ approximately as a Probit
function $\overline{\Phi}( (\rho - \rho(\delta; Q))/w(\delta,N; Q) )$. Here 
$\overline{\Phi}$ is the Gaussian survival function and 
$w$ is the transition width.  A logit function fits just about as well.
\eitem

\noindent
{\sl Analysis of $Z$-scores.} We compare success probabilities
at the non-Gaussian ensembles
to those at the Gaussian using $Z$-scores arising from two-sample tests for binomial
proportions.
\bitem
\item {\it Methodology of $Z$-score comparison.}
We verify that in the null case of no difference, our methodology indeed
finds no difference; we also verify that in the case of known difference,
it indeed finds a difference.
\item {\it Scaling of moments with $N$.}  We identify nonzero
means in the $Z$-scores, and show that the
scaling law $\mu(\delta, N ; E) = O(1/N^{1/2})$
best describes the data.
\item {\it Exceptional ensembles.} Two matrix ensembles
exhibit substantial lack-of agreement with the others
(e.g. some individual $Z$-scores as large as 20) for $n$ small.
In effect the location for 50\% success in those ensembles
obeys a slight shift away from $\rho(\delta; Q)$, of order $w$.
While this is an asymptotically negligible shift, failing to model 
it causes a noticeable lack of fit at $N=200$ and $n$ small.  This lack of agreement 
is observed to dissipate as $n$ increases.
\item {\it Validation ensembles.} Two matrix ensembles
were studied only after all other analysis had been completed.
Using the models arrived at in the prior analysis without
changing the model form,  we found that the same models
describe the validation ensembles adequately,
reinforcing the validity of our analysis.
\eitem

All noticeable elements of lack of fit are best accounted for as
evidence of effects consistent with the weak universality hypothesis.

\subsection{Experiments conducted}

\subsubsection{Framework}

Terminology:
\bitem
 \item We study $n\times N$ random matrices $A$, $N > n$. 
\item A matrix ensemble is a generating device for $n \times N$ matrices.
We report here results on the 9 different random matrix ensembles , 
listed in table \ref{tab:ensembles}.
  \item We generate vectors $y=Ax_0$ where $x_0$ has $k$ nonzeros.
  \item The nonzeros in $x_0$ are either drawn uniformly 
from $(0,1)$ or $(-1,1)$, designated as {\em coefficients} $+$ and $\pm$ 
respectively.
  \item The instances $(y,A)$ where the underlying $x_0$ is nonnegative
  by intent are then processed using an optimizer to approximately solve $(LP)$. 
  The instances where the underlying
  $x_0$ can be of both signs are processed by using an optimizer to
  approximately solve $(P_1)$.
 \footnote{ In principle, the precise values of the nonzeros do not matter for
 properties of $(P_1)$ and $(LP)$. (n.b. For other sparsity seeking algorithms 
 this would not be the case.)  }
 \item The optimizer is presented with the problem instance $(y,A)$, but not $x_0$.
 \item After running the optimizer, we measure
 {\tt ExactRecon}, which takes the value 1 when
the obtained solution $x_1$, say, is equal to the
desired solution $x_0$, within 6 digits accuracy.
It is zero otherwise.
\item We conduct $M$ independent replications at each fixed combination of $N,n,k$,
matrix ensemble, and coefficient type. 
\item The variable $S$ totals the number of times  {\tt ExactRecon}
was $1$ in the $M$ replications.  Results are tabulated in a data file with
column headings
\begin{verbatim}
E N n k M S
\end{verbatim}
Here $E$ is an integer code specifying the {\em suite} of problem instances.
Such a suite specifies both the matrix ensemble (eg Gaussian, Bernoulli, Rademacher, ...) and the coefficient type ($+$ or $\pm$).  For each matrix ensemble we consider both coefficient types.
\item For analysis and presentation, we use coordinates 
 $\delta = n/N$, the matrix `shape', and the solution sparsity level $\rho = k/n$.
 We generally consider the success fraction $\hat{p} = S/M$.
 \item The most important structure of the dataset concerns the {\it constant-$\delta$ slices},
 where $n$, $N$, $E$ are held constant and $k$ is varying.
  The success fraction $\hat{p}(k,n,N;E)$ is generally monotone
 decreasing in such a slice: monotone decreasing in $k$ for fixed $n$, $N$, and $E$.  
 \item We focus on what is called the $LD50$ in bioassays, the 50\% quantal response more generally.
 It is the value of $k/n$ where $S/M$ is expected
 to be $1/2$, for fixed $n$,$N$, $E$. 
 \eitem
 
\subsubsection{Range of experiments}

For each suite in table \ref{tab:ensembles} and each combination of $N,n,k$, 
we consider $M=200$ different 
problem instances $(y,A)$ each one drawn randomly as above.
Each suite is compared with a ``baseline'' of either suite 1 or 2 
for the same combination of $N,n,k$ but a larger \underline{independent} 
draw of $M=1000$ problem instances.
We varied the matrix shape $\delta = n/N$, and the solution sparsity 
levels $\rho = k/n$.
 At problem size $N = 1600$, we varied $n$ systematically through a grid
 ranging from $n = 160$ up to $n = 1440$ in 9 equal steps.
The `signal processing language' event `exact reconstruction'
corresponds to the `polytope language'  event `specific $k$-face
of $Q$ is also a $k$-face of $AQ$'. In both cases we speak of {\it success},
and we call the frequency of success in $M$ empirical trials at a given
$(k,n,N)$ the {\it success rate}.
At each combination $N,n$, we varied $k$ systematically to sample the success 
rate transition region from $5\%$ to $95\%$.

\subsubsection{Suites studied}

In sections 3-7 we report details
about experiments with the non-Gaussian suites $3-12$ and $15-16$
listed in table \ref{tab:ensembles}.  As it happens, after the analysis of these
suites was conducted, data became available for four other suites,
based on the Hadamard and Rademacher ensembles. The analysis
of those suites will be reported only in section 8. 
 The Rademacher ensemble, suites 19 and 
20, generated data with the same problem sizes and 
other parameters as suites $3-12$ and $15-16$.
Our study of the Hadamard ensemble, suites 13 and 14, 
is restricted since only two problem sizes $N=256$ and $N=512$
were run.

\subsection{Behaviour of the Gaussian ensemble}

In this section, we restrict attention to the Gaussian ensemble, suites 1 and 
2, and investigate these questions:
\bitem
 \item {\it Width of transition zone.}  How does the width $w(\delta; N, Q)$
 of the transition zone at phase transition vary,
 as a function of $N$? 
 \item {\it Quantal Response Profile.}   How does the probability of success
 vary as a function of the reduced $(\rho - \rho(\delta; Q))/w(\delta; N,Q)$?  
 Where $\rho=k/n$ and $\delta=n/N$.
 \item {\it Behaviour of $LD50$.}  In what manner does the empirical 50\% point 
 of the quantal response function approach the underlying asymptotic limit
 $\rho(\delta; Q)$?
 \eitem
 
 The Gaussian ensemble is an appropriate place to focus attention, because:
 \bitem
  \item A complete, rigorous understanding of the asymptotic behaviour exists
  in the Gaussian case (Donoho \& Tanner 2009$a$); we know that as $N \goto \infty$
  with $\delta$ fixed and $\rho$ fixed away from $\rho(\delta,Q)$, 
  the success probability tends to either zero or 1. So we know that 
  there {\it is} a
  transition zone, and that its width tends to $0$ as $N \goto \infty$.
  \item A rigorous set of finite-$N$ bounds has been rigorously proven 
  (Donoho \& Tanner 2008$b$); we know that the width scales like $O(1/\sqrt{n})$.
  \eitem
  Hence, there are rigorous theoretical constraints: we know the phase 
transition exists asymptotically and we can constrain its width.
  
 \subsubsection{Modelling the quantal response function}
 In the field of bioassays, the Quantal response function gives the
 probabiliy of organism failure (eg death) as a function of dose.  
 In our setting, the analogous concept is the probability of
 algorithm failure at a fixed problem size $n,N$ as a function of
 $\rho=k/n$, the ``complexity dose''.
 
\begin{figure}
\begin{center}
\includegraphics[bb=76 215 543 591,width=4in,height=2.1in]{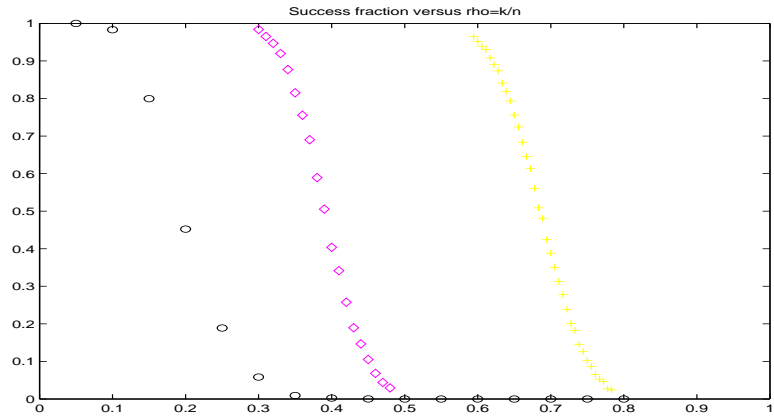} 
\caption{Empirical success fraction at Gaussian ensemble, suite 2,
for $\delta = .1$ (black circles), $\delta = .5$ (magenta triangles),
$\delta = .9$ (yellow  diamonds). Horizontal axis: $\rho = k/n$;
vertical axis: success fraction $S/M$. Data result from $M=1000$
trials at $N=1600$.}
\label{DoseResponse}
\end{center}
\end{figure}

 Considering a constant-$\delta$ slice at the Gaussian ensemble, suite 2,
 we see a roughly monotone increasing probability of failure
 as a function of $k/n$. Figure \ref{DoseResponse} presents
 the fraction of success $S/M$ as a function of $k/n$, at three
incompleteness ratios $\delta= n/N$.
 
\begin{figure}
\begin{center}
\includegraphics[width=4in]{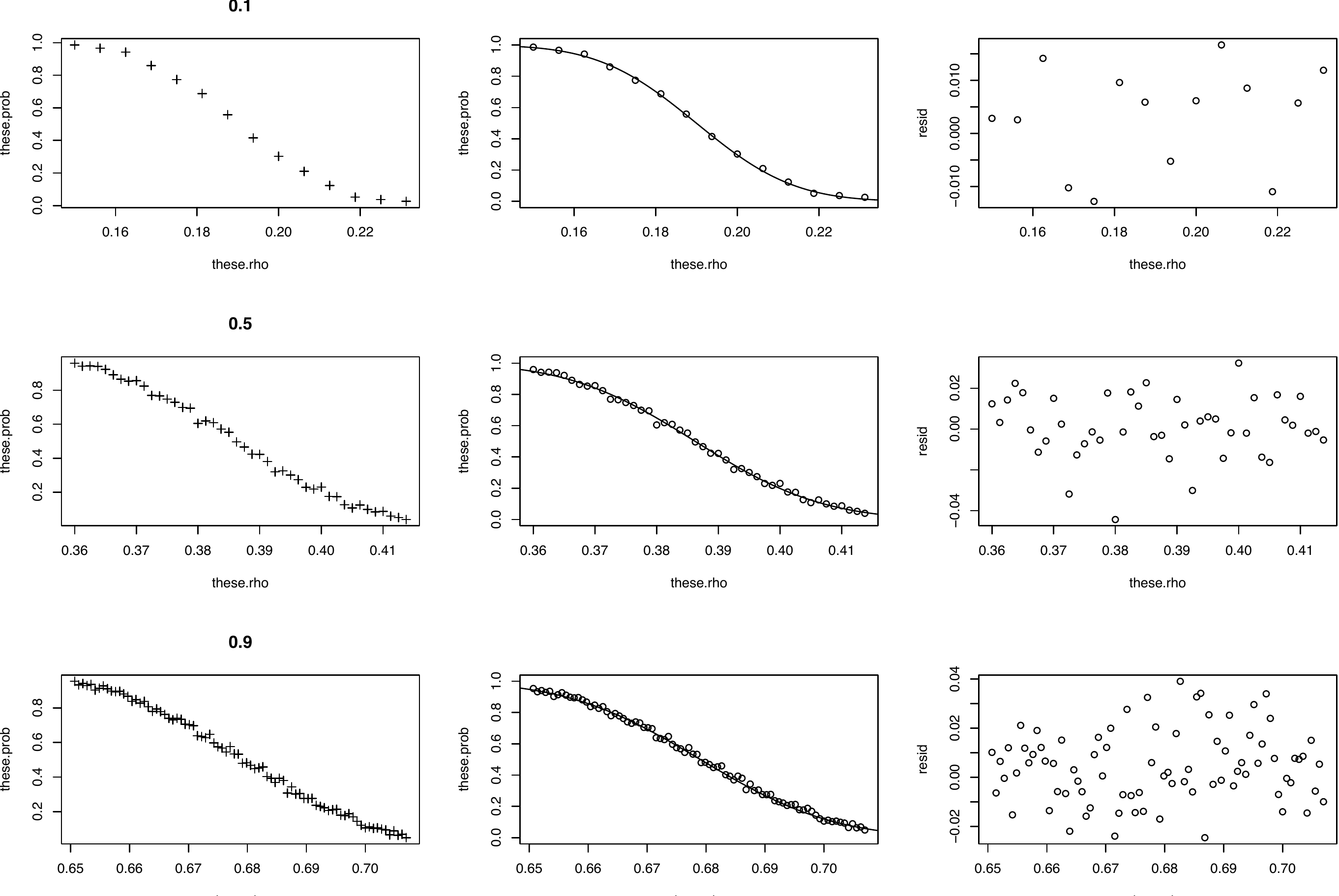} 
\caption{Crude probit modeling of dose response.
Rows for $\delta = .1$, $\delta = .5$,
$\delta = .9$. First column: same
as in figure \ref{DoseResponse} for suite 2, $N=1600$ and $M=1000$. Second Column: data and probit model. Third Column; Residuals. }
\label{CrudeFitProbit}
\end{center}
\end{figure}

A Probit model for the dose response states that, for parameters $a$, and $b$,
the expected fractional success rate is given by
\begin{equation}\label{crudeProbit}
     {\cal E}(S/M) = \bar{\Phi}( a(\delta) + b(\delta)\rho) 
\end{equation}
where $\bar{\Phi}$ is the complementary normal distribution, and ${\cal E}$ denotes expectation.
Figure \ref{CrudeFitProbit} presents a 
first pass at checking the suitability
of such a model. It identifies empirical estimates
of the points where $ {\cal E}(S/M) = \alpha$ for $\alpha \in \{  1/4, 1/2, 3/4 \}$.
and then chooses $a$ and $b$ in relation  (\ref{crudeProbit})
to match those. It displays the raw data, the model curves, and
residuals from the model.

 It is standard in biostatistics to fit generalized linear models
 to such data; the binomial response model is appropriate here.
 Such models take the form
 \[
       S \sim Bin(p(\delta,\rho), M)
 \]
 where the success probability $p(\delta)$, after
 a fixed transformation $\eta()$ , obeys a linear model:
 \[
    \eta(p(\delta,\rho))  = a(\delta) + b(\delta)\rho;
 \]
 the function $\eta$ is called the {\it link} function.

\begin{figure}
\begin{center}
\includegraphics[width=4in]{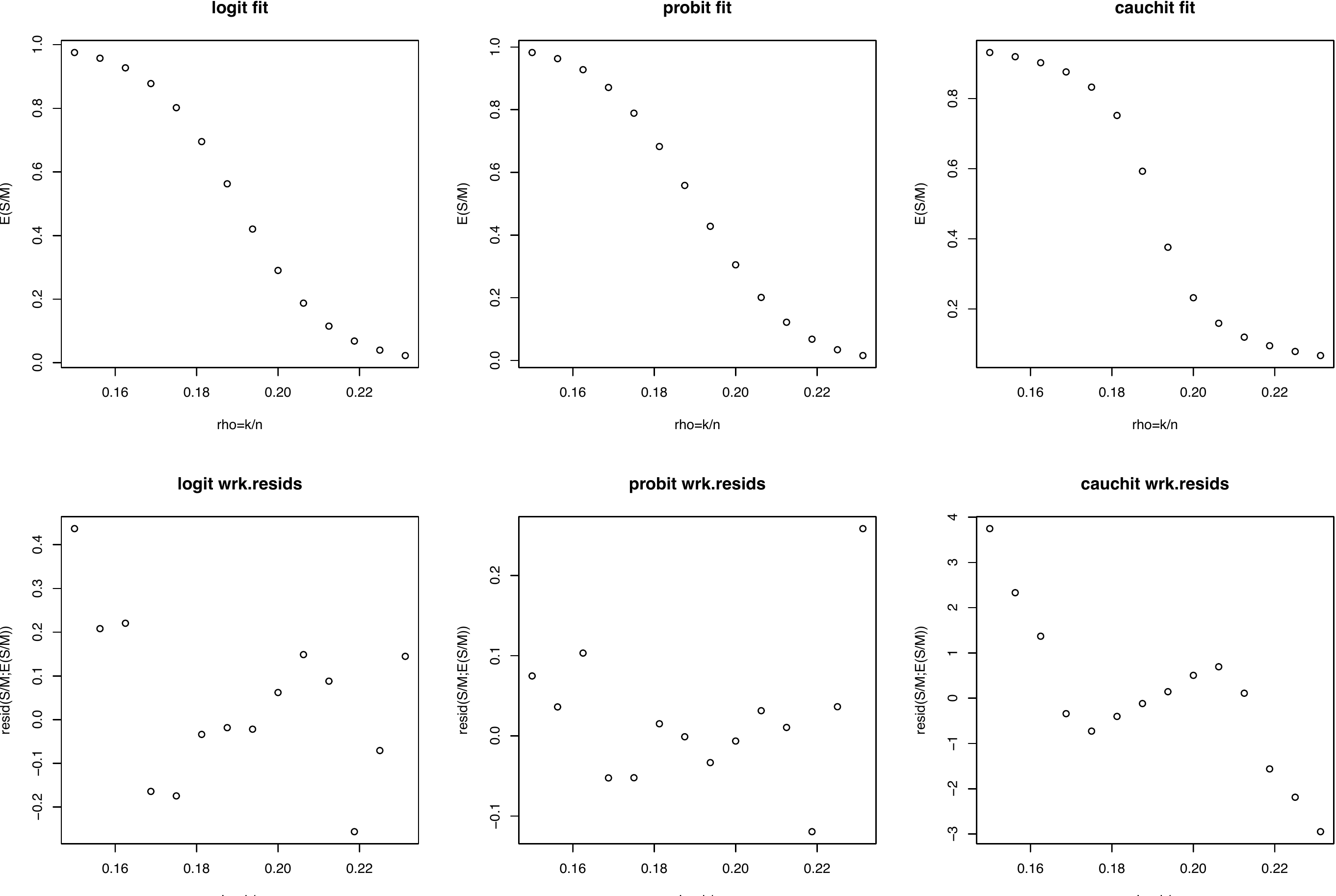} 
\caption{GLM modeling of dose response for suite 2 ($N=1600$ and $M=1000$) with three links;
case $\delta = .1$. First Column: fit with Logit link. Second Column: fit with 
Probit link. Third Column: fit with Cauchyit Link. 
First Row: Fits. Second Row: Working Residuals.}
\label{GLMFits}
\end{center}
\end{figure}

We considered three standard link functions:
the logit, probit  and cauchyit links.
Figure \ref{GLMFits}  presents fitted models
and what the statistics analysis package R calls 
the {\it working} residuals for these three links.
In fact the best loglikelihood is achieved 
among the three at the probit link,
but there is a large residual at the most extreme
response; it seems the probit link goes to zero
too fast (this is not unexpected, owing to the 'thin tails'
of the normal distribution, and also owing to
the finite-N large deviations analysis in Donoho and Tanner (2008$b$)).  
The logistic link is nearly as good in
deviance or likelihood senses,  makes sense
on theoretical grounds and gives more balanced residuals.
(Note however, that as figure \ref{CrudeFitProbit}
showed, the Probit fit is adequate as long as we look at
ordinary rather than the more statistically sensitive working residuals.)

We used R to
 fit these models; this has the advantage of automatically
 providing standard inferential tools -- confidence bounds
 for $a$ and $b$ and   goodness-of-link tests.

\subsubsection{Behaviour of $LD50$}

In bioassays, the $LD50$ is the dose that
corresponds to 50/50 chance of failure. 
This can be estimated from binomial response data
in two ways.

The first, `nonparametric' method
finds the largest ratio $k/n$ where
 $S > M/2$ for a given $n$, $N$ and $E$.
 We found a slight refinement useful:
we fit a linear spline to the success ratios
and solved for the (smallest) value of $\delta$
where the spline crosses  $50$\%.
 
 Using this method, we obtained table \ref{tab-2},
 which presents, for suite 2 and
 $M=1000$, the difference between
 the estimated $LD50$ and the theoretical large $N$ limit.

\begin{table}
\begin{center}
\begin{tabular}{|l|l|l|l|}
\hline
 $\delta$ &     200     &      400      &    1600 \\
\hline
0.2 & 0.0114 & 0.00410  &  0.001799 \\
0.3 & 0.0079 & 0.00476  &  0.001179 \\
0.4 & 0.0043 & 0.00444  &  0.001037 \\
0.6 & 0.0053 & 0.00198  &  0.000760 \\
0.7 & 0.0048 & 0.00380  & 0.001649 \\
0.8 & 0.0082 & 0.00280  & 0.001171 \\
\hline
\end{tabular}
\caption{Difference between empirically estimated $LD50$ of suite 2
with $M=1000$ samples and the corresponding asymptotic limit
$\rho(\delta; C)$}\label{tab-2}
\end{center}
\end{table}
Evidently, the $LD50$ is approaching the
expected phase transition with increasing $N$.
To quantify this effect, we have table \ref{tab-3},
which shows that the $LD50$ typically approaches
its limit at roughly the rate $1/N$.
\begin{table}
\begin{center}
\begin{tabular}{|l|l|l|l|}
\hline
Quantity &        200   & 400 & 1600 \\
\hline
median $N^{1/2} (LD50 - \rho(\delta)) $ & 0.094&  0.079 &0.047 \\ 
median $N           (LD50 - \rho(\delta)) $ & 1.324 &  1.581   & 1.881  \\
\hline 
\end{tabular}
\caption{Median across $\delta$ of the scaled difference  between empirically estimated $LD50$ of suite 2 with $M=1000$ 
and the corresponding asymptotic limit $\rho(\delta; C)$}\label{tab-3}
\end{center}
\end{table}

%
 
\subsubsection{Transition zone width}

We can define the  $\alpha$-width of the transition
 zone as  the horizontal distance
 between $p = \alpha$ and $p = 1-\alpha$
 on the dose-response. 
  
 We again can measure this nonparametrically
 and parametrically.   We present here a nonparametric
 analog based on fitting splines to the empirical success fractions 
and measuring the $\alpha=0.1$ and $1-\alpha=0.9$
quantile locations.  We then normalize by the corresponding 
distance on the standard Probit curve
\[
     w=w_{0.1} = \frac{ q_{0.9} - q_{0.1}  }{\Phi^{-1}(0.9) - \Phi^{-1}(0.1)}.
\]
 Table \ref{tab-4} presents values of $\sqrt{N} \cdot w(\delta,N)$ for 
suite 2 with $M=1000$.
 \begin{table}
 \begin{center}
\begin{tabular}{|l|l|l|l|}
\hline
$\delta$ & $N=200$ & $N=400$ & $N=1600$ \\
\hline 
0.2 & 0.7277 & 0.6691 & 0.6367 \\
0.3 & 0.6567 & 0.6535& 0.6780 \\
0.4 & 0.6530 & 0.6464 & 0.6240 \\
0.6 & 0.6328 & 0.6663 & 0.6420 \\
0.7 & 0.6413 & 0.6690 & 0.6384 \\
0.8 & 0.6708 & 0.6841 & 0.6809 \\
\hline
\end{tabular}
\caption{Scaled values  $\sqrt{N} w(\delta; N)$ for suite 2 with $M=1000$.}\label{tab-4}
\end{center}
\end{table}

%

\subsection{Methodology of $Z$-score comparison}

How does the methodology of $Z$-score comparison work
on cases where we know the ground truth -- both
where we know there is no difference and we know
there is an asymptotic difference?  For each of the suites 
in table \ref{tab:ensembles} the suite with $M=200$ is compared 
against suite 1 or 2 with $M=1000$ independent problem instances 
for the same values of $N,n,k$.
Suites 1 and 2 with $M=1000$ form the baseline against which all 
$Z$-scores are calculated.
Unless specified otherwise, suites with the same coefficient sign 
are compared; for example suite 9 is compared with suite 1 and suite 
16 is compared with suite 2.

\subsubsection{Under a true null hypothesis}

When the two problem settings being compared are
simply replications of the same underlying conditions,
we can be sure that $H_0$: {\it no difference} is true.
We compare here suites 1 and 2 with $M=200$ against the 
baseline experiments of suites 1 and 2 with $M=1000$.
This follows the same procedure as will be conducted later 
when comparing non-Gaussian ensembles against the baseline.



\begin{figure}
\begin{center}
\includegraphics[width=3in]{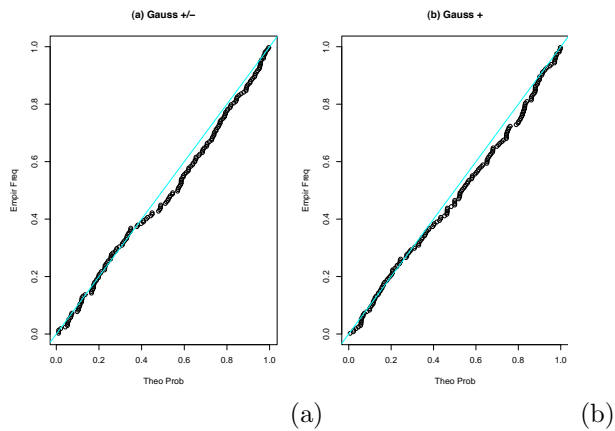} 
\end{center}

\vspace{-0.2in}
\hspace{1.75in}
(a) \hspace{1.25in} (b)
\caption{PP Plots of $Z$-scores comparing success frequencies 
under one set of realizations from the Gaussian matrix ensemble ($M=200$)
with an independent realization of success frequencies
from the Gaussian ensemble ($M=1000$).
Panel (a): suite 2.
Panel (b): suite 1.}
\label{fig-bulk-ZScores}
\end{figure}

Figure \ref{fig-bulk-ZScores} Panel(a)
presents the bulk distribution of $Z$-scores for 
suite 2;
Panel (b) presents the bulk distribution of $Z$-scores for 
suite 1.
The figure presents PP-plots: the fraction of $Z$-scores
exceeding a threshold versus the fraction to be expected at the standard Normal.
If the $Z$-scores were exactly standard normal, these  
plots would be close to the identity line,
which is, in fact what we see.

In quantitative terms, we have table \ref{tab:Znull}.

\begin{table}
\begin{center}
\begin{tabular}{|l|l|l|l|l|}
\multicolumn{5}{c}{Suite 1: Gaussian ensemble, positive coefficients}\\
\hline
$N$ &  \#$\{i: |Z_i| < 1\}$ & \#$\{i:|Z_i|<2\}$ & \#$\{i:|Z_i|<3\}$ & $\#$ \\
\hline
 200 &29/54 &  49/54  & 54/54 & 54 \\ 
 400 & 40/63 &  60/63 & 63/63 & 63 \\ 
1600 &45/63.6 & 61/63 & 63/63 &63 \\ 
\hline
\multicolumn{5}{c}{Suite 2: Gaussian ensemble,  coefficients of either sign}\\
\hline
 200 & 41/55 & 53/55 & 55/55 & 55 \\ 
 400 & 38/63 & 60/63 & 63/63  & 63 \\ 
1600 &42/63 & 62/63 & 63/63  & 63 \\ 
\hline
\end{tabular}
\caption{Occurrences of $Z$-scores under the true null hypothesis. Suites 1 and 2 self comparison of $S/M$ with $M=200$ against the basline suites 1 and 2 with $M=1000$.
Cell contents: Observed Counts/Total Counts.}
\label{tab:Znull}
\end{center}
\end{table}

Table \ref{tab:Znull} shows that 
of 180 $Z$-scores associated with 
comparisons of suite $1$,
170 were less than $2$ in absolute value;
for $94.4$ \%  -- very much in line 
with an assumed standard $N(0,1)$ null distribution.
It also shows that of
181 $Z$-scores associated with 
comparisons of suite $2$,
175 were less then $2$ in absolute value;
for $96.7$ \%  -- very much in line 
with an assumed standard $N(0,1)$ null distribution.
We thus see that under a true null hypothesis, our
$Z$-scores behave largely as if they were $N(0,1)$.  This is an observation
that needed to be checked, since $Z$-scores,
when constructed in the way we have done so here, only are
known to have an asymptotically normal distribution.

Figures \ref{fig:z1} and \ref{fig:z2}  presented $Z$-scores
in scatterplots of $Z(k,n,N; E)$ versus $\delta$ comparing non-Gaussian 
ensembles with Gaussian ensembles; what happens
when we have a true null hypothesis?

\begin{figure}
\begin{center}
\begin{tabular}{cc}
\includegraphics[bb=76 215 543 591,height=1.5in,width=2.3in]{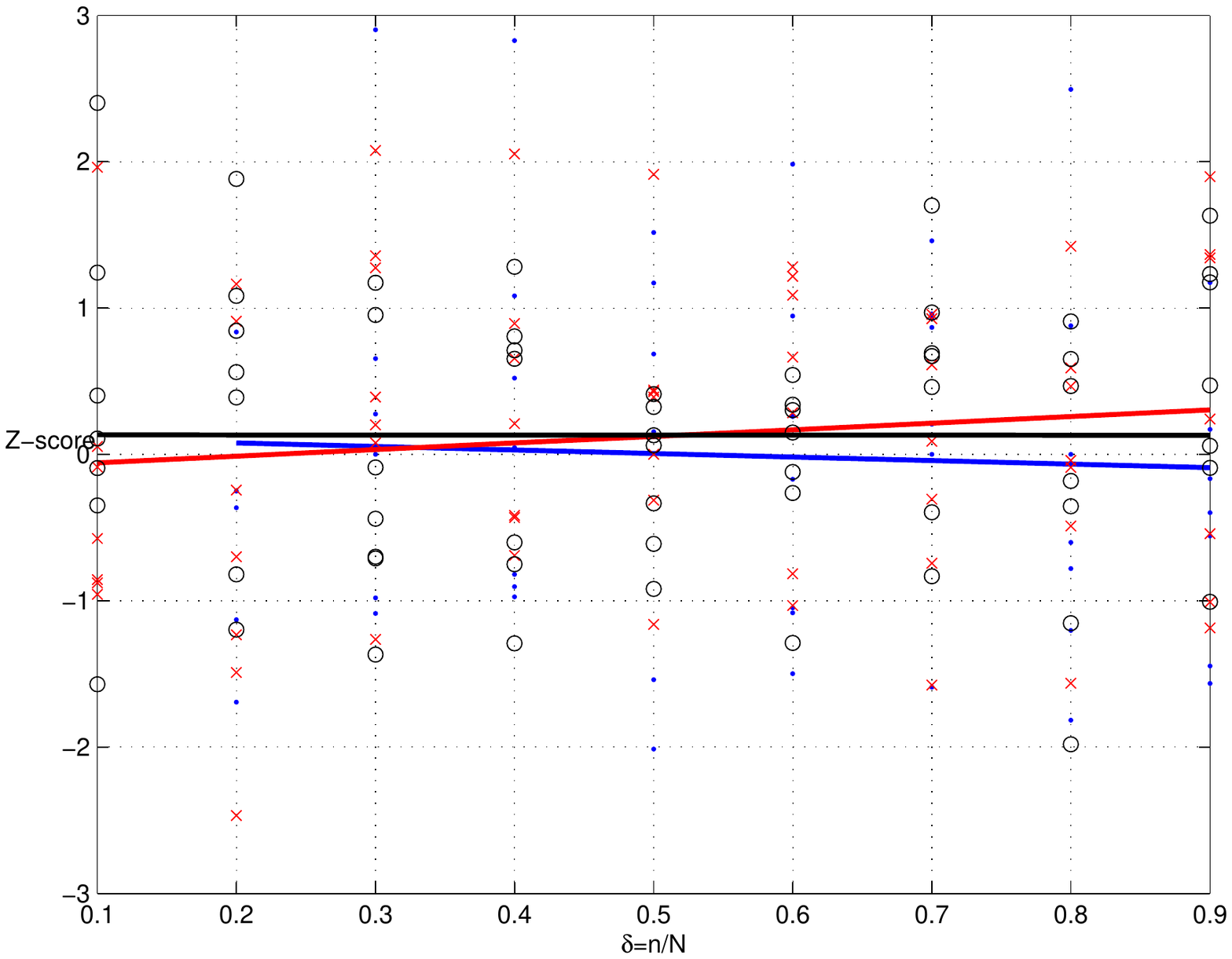} &
\includegraphics[bb=76 215 543 591,height=1.5in,width=2.3in]{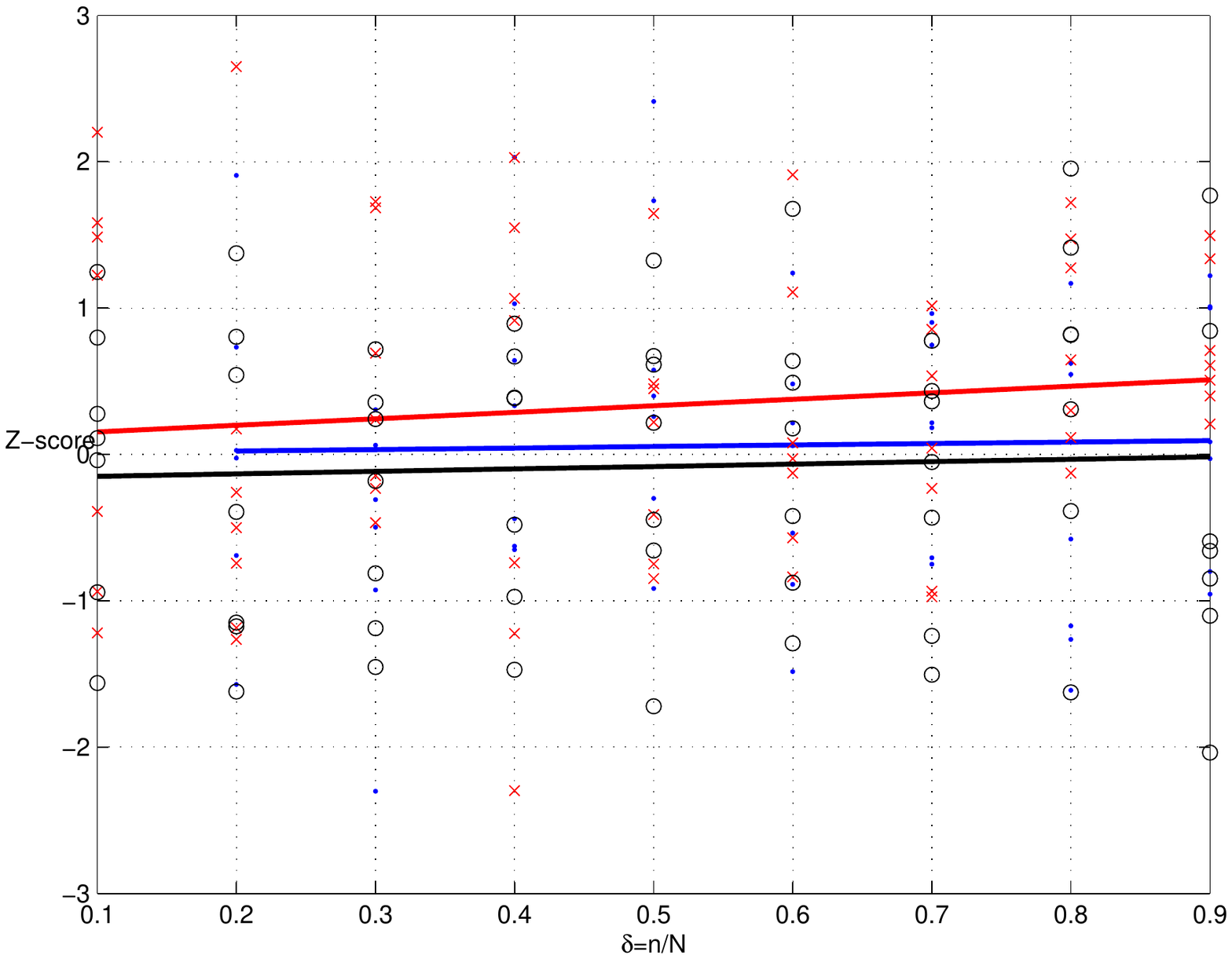} \\
(a) & (b)
\end{tabular}
\caption{$Z$-scores under true null hypothesis, and linear fit.  
Panel (a): suite 1.
Panel (b): suite 2.}
\label{fig-linearfit-Z-Gaussian-unsigned}
\end{center}
\end{figure}

Figure \ref{fig-linearfit-Z-Gaussian-unsigned} (a) and (b) presents 
the $Z$-scores for suites 1 and 2 respectively with fitted lines modelling 
the dependence on $\delta$.
In principle, the $Z$-scores all have mean zero
and there is no expected trend. However, owing to sampling fluctuation,
we obtain nonzero intercepts and slopes.
Table \ref{tab:Zlinfit} shows the results that obtained
in this truly null case. We learn from this that fitted intercepts and
slopes of about the size indicated in the table can be viewed
as consistent with a true null hypothesis.

\begin{table}
\begin{center}
Suite 1
\hspace{1.45in}
Suite 2

\vspace{0.03in}
\begin{tabular}{|l|l|l|l|l|}
\hline
$N$ & $a$ &$b$\\
\hline
200&  0.127 &  -0.242\\
400& -0.103 &  0.451\\
1600& 0.133&  -0.002\\
\hline
\end{tabular}
\hspace{0.5in}
\begin{tabular}{|l|l|l|l|l|}
\hline
$N$ & $a$ &$b$\\
\hline
200&  0.002 &  0.102\\
400& 0.109  &  0.445\\
1600&-0.166 & 0.166\\
\hline
\end{tabular}
\caption{Linear fits to $Z$-scores $Z \sim a + b \delta $ 
under the true null hypothesis. Suites 1 and 2.
$a$ intercept; $b$ slope}
\label{tab:Zlinfit}
\end{center}
\end{table}

\subsubsection{Under a true alternative hypothesis}

Is our methodology powerful? Can it detect any differences 
from null?

To study this question, we considered a simple and blatant
mismatch: compare suite 1 with $M=200$ against the suite 2 
with the baseline $M=1000$.

\begin{figure}
\begin{center}
\includegraphics[width=3in]{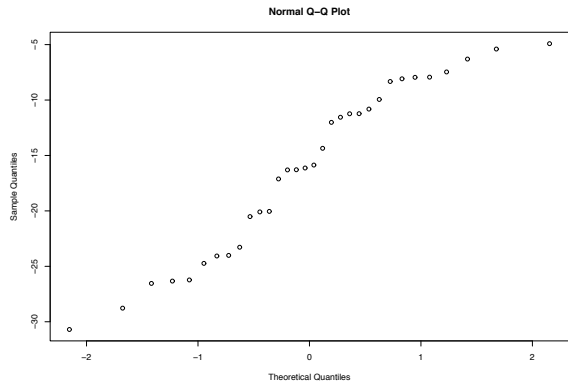} 
\caption{QQ plot of $Z$-scores comparing success frequencies 
of suite 1 ($M=200$)
to an independent realization of success frequencies
from suite 2 ($M=1000$).
In this plot, the $M=200$-based frequencies
describe success in solving $(LP)$ when the 
nonzeros in $x_0$ are positive, 
 the $M=1000$-based frequencies
describe success in solving $(P_1)$ when the 
nonzeros in $x_0$ are of either sign. 
The two ensembles are proven to have
different phase transitions. 
The values on the $Y$ axis are far away from $0$.}
\label{fig-bulk-ZScores-nonnull}
\end{center}
\end{figure}

The baseline, suite 2, should reflect a transition near $\rho(\delta; C)$
while the comparison group, suite 1, should reflect a transition 
near $\rho(\delta ; T)$.
As these two curves are very different we should see this 
reflected in the $Z$-scores.
And we do.  Figure \ref{fig-bulk-ZScores-nonnull} shows 
the QQ-plot, which is noticeably far from the identity line.

\subsection{Bulk behaviour of $Z$-scores}

\begin{figure}
\begin{center}
\begin{tabular}{c}
\includegraphics[width=3in]{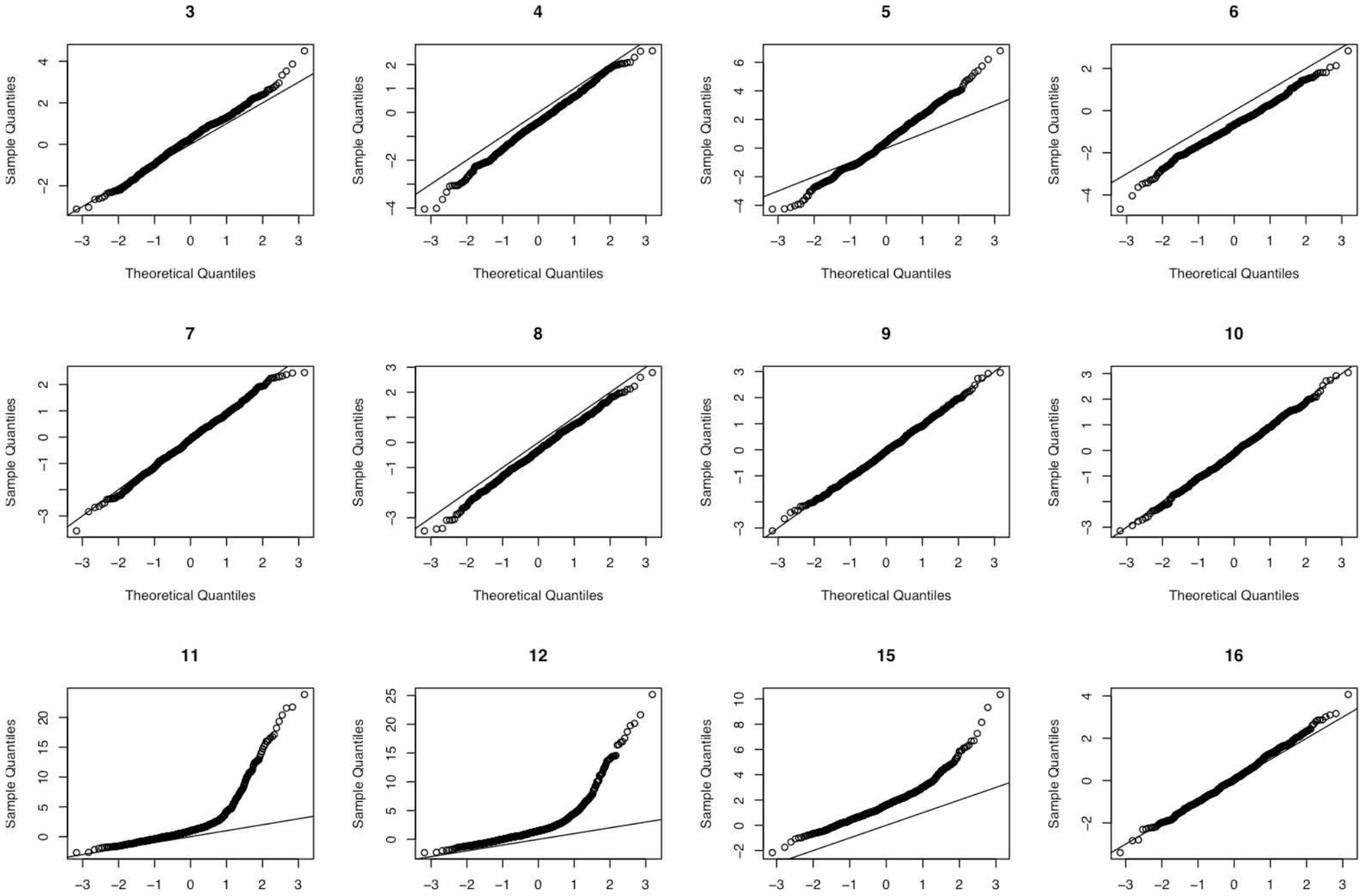} \\
\end{tabular}
\caption{$Z$-scores at $N=200$ for suites 3-12,15-16.}
\label{fig-Bulk-Zscores-N-200}
\end{center}
\end{figure}

\begin{figure}
\begin{center}
\begin{tabular}{c}
\includegraphics[width=3in]{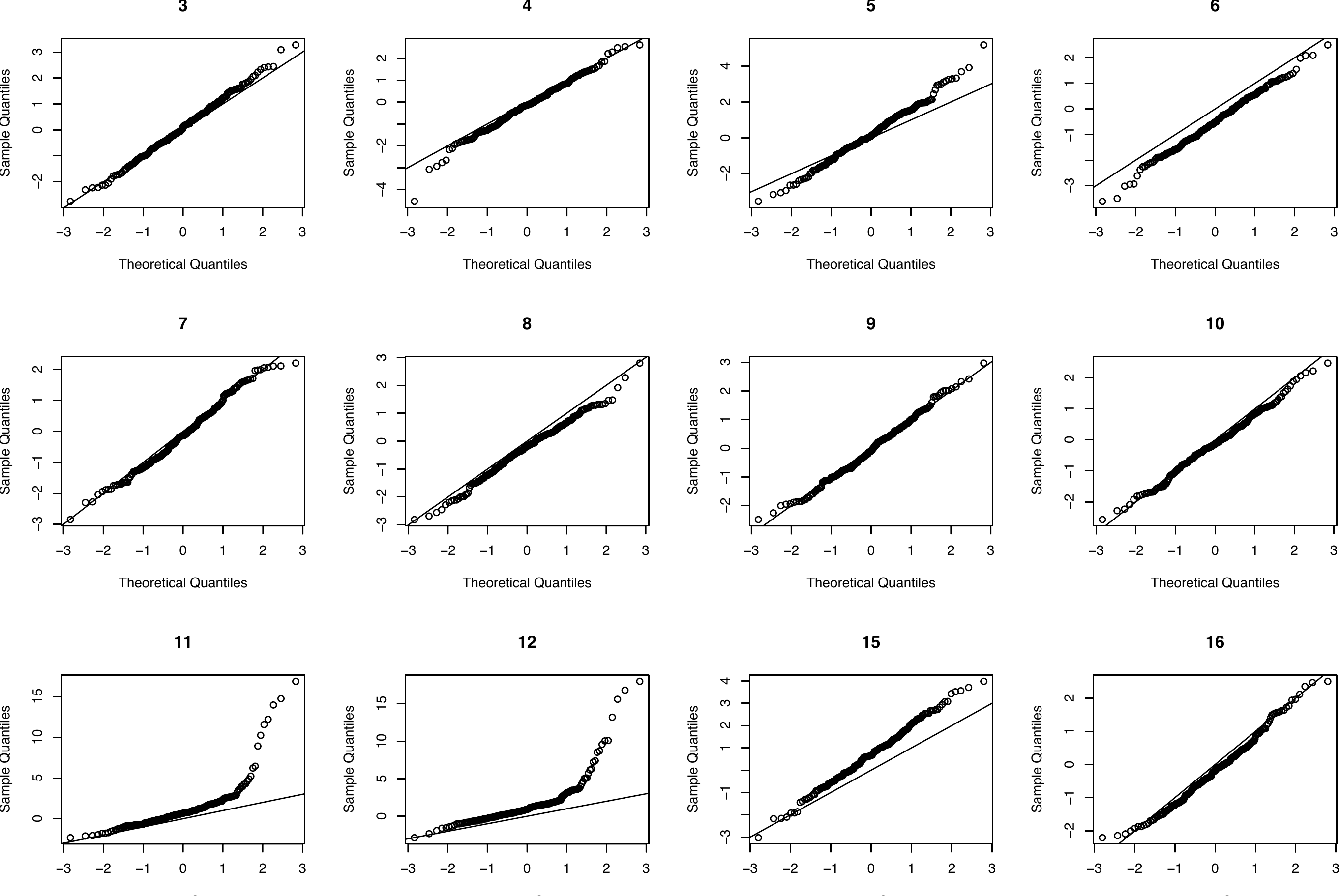} \\
\end{tabular}
\caption{$Z$-scores at $N=400$ for suites 3-12,15-16.}
\label{fig-Bulk-Zscores-N-400}
\end{center}
\end{figure}

\begin{figure}
\begin{center}
\begin{tabular}{c}
\includegraphics[width=3in]{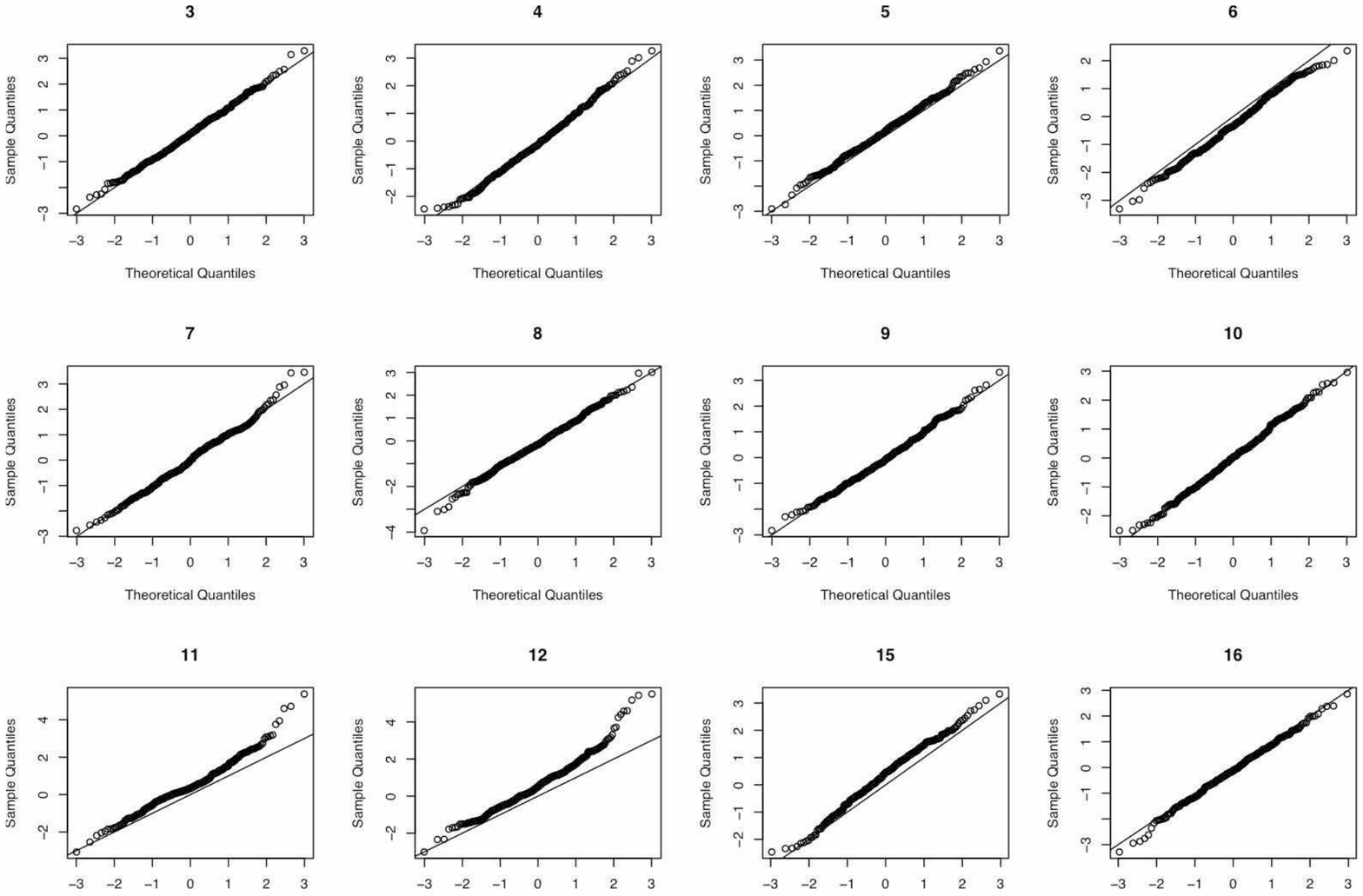} \\
\end{tabular}
\caption{$Z$-scores at $N=1600$ for suites 3-12,15-16.}
\label{fig-Bulk-Zscores-N-1600}
\end{center}
\end{figure}

In figures \ref{fig-Bulk-Zscores-N-200}-\ref{fig-Bulk-Zscores-N-1600}
we present a sequence of  QQ plots showing the bulk distribution
of $Z$-scores for suites $3-12,15-16$ at $N=200$, $400$ and $1600$.
In each plot the identity line is also displayed; if the $Z$-scores
had truly a standard normal distribution, they would oscillate around this line.

While in most suites we see a good match
between the $Z$-scores and the standard normal
already at $N=200$, it is not until $N=1600$
when every ensemble seems to yield 
approximately $N(0,1)$ scores.  Even then,
suites 11 and 12 (Ternary (1/10)) exhibit some noticeable
deviations.  These effects are apparent for small $n$.  In fact, at small $n$, 
trivial linear dependencies are found 
among the columns of typical random $A$ at these ensembles.
Such linear dependencies force affine dependencies
which force lost polytope faces. Consequently,
the observed $S/M$ will be decreased, this effect is discussed 
further in \S \ref{sec:exceptional}.

This is consistent with our conclusions in the main text that
strict finite-$N$ universality does not hold, but a weaker
notion of asymptotic universality does hold.

\subsection{Linear modeling of the $Z$-scores}\label{LinearFits}

The QQ Plots of $Z$-scores in figures 
\ref{fig-Bulk-Zscores-N-200}-\ref{fig-Bulk-Zscores-N-1600}
show that suites 11, 12, 15, and 16 (the highly sparse matrix suites)
behave somewhat differently than other suites.

We now report results \underline{excluding} those suites,
and focus instead on suites 3-10.  The excluded suites $11-12,15-16$
will be discussed in \S \ref{sec:exceptional}.

In the main text we reported fits explaining the
observed $Z$-scores using two models of the form,
\begin{equation}
\label{satmodel}
     \mu(\delta; N, E)  = \alpha(N,E) + \beta(N,E) (\delta - 1/2),
\end{equation}
where
\begin{equation} \label{linmodel}
      \alpha(N,E) = \alpha_0(E) +\alpha_1(E)/N^{1/2}, \qquad \beta(N,E) = \beta_0(E) + \beta_1(E)/N^{1/2}.
\end{equation} 
We considered two such models, one with 
the restrictions $\alpha_0 = 0$
and $\beta_0 = 0$, and one without.

We first report results using the submodel
$\alpha_0 = 0$ and $\beta_0 = 0$, and later the full model,
without the restriction. The display below
shows the output of the $R$ linear model
command for the restricted model.  
The symbol {\tt de} denotes  $\delta - 1/2$
and the symbol {\tt probSize} denotes 
$1/\sqrt{N}$.   The output lists interaction
terms such {\tt probSize:En}.
It turns out that  $\alpha_1(N,E)$ is the sum of   {\tt probSize} and  {\tt probSize:En}.
Similarly $\beta_1(N,E)$
is the sum of   {\tt probSize:de} and  {\tt probSize:En:de}.

\begin{verbatim}
lm(formula = Z ~ probSize * E * de - E * de - E - de - 1, subset = subE)

Residuals:
      Min        1Q    Median        3Q       Max 
-3.966108 -0.703283 -0.003336  0.691872  4.309933 

Coefficients: (1 not defined because of singularities)
                        Estimate    Std. Error t value Pr(>|t|)    
probSize         -1.7714     0.5752  -3.080  0.00208 ** 
probSize:E3       2.5960     0.7955   3.263  0.00111 ** 
probSize:E4      -6.1540     0.8141  -7.560 4.40e-14 ***
probSize:E5       3.3350     0.7972   4.183 2.90e-05 ***
probSize:E6      -9.0306     0.8128 -11.110  < 2e-16 ***
probSize:E7      -0.6122     0.7957  -0.769  0.44168    
probSize:E8      -4.8878     0.8117  -6.022 1.79e-09 ***
probSize:E9       0.4795     0.7964   0.602  0.54715    
probSize:E10          NA         NA      NA       NA    
probSize:de      32.4608     2.2249  14.590  < 2e-16 ***
probSize:E4:de  -13.3427     3.0552  -4.367 1.27e-05 ***
probSize:E5:de   38.3394     3.1584  12.139  < 2e-16 ***
probSize:E6:de  -25.7357     3.0616  -8.406  < 2e-16 ***
probSize:E7:de  -19.7373     3.1520  -6.262 3.96e-10 ***
probSize:E8:de  -17.3618     3.0514  -5.690 1.31e-08 ***
probSize:E9:de  -25.3606     3.1504  -8.050 9.22e-16 ***
probSize:E10:de -30.1162     3.0533  -9.863  < 2e-16 ***
---
Signif. codes:  0 Ô***Õ 0.001 Ô**Õ 0.01 Ô*Õ 0.05 Ô.Õ 0.1 Ô Õ 1 

Residual standard error: 1.042 on 9980 degrees of freedom
Multiple R-Squared: 0.1801,	Adjusted R-squared: 0.1788 
F-statistic:   137 on 16 and 9980 DF,  p-value: < 2.2e-16 
\end{verbatim}

One sees that the coefficients in this model are significant,
and that the residual variance is roughly $1$, as is to be 
expected of proper $Z$-scores. (The reduction in sum of squares
in the null case due to fitting 16 degrees of freedom would
negligible and has no significant bearing on 
our evaluation of the adequacy of residuals).

In contrast, here is the report on the fit of the full
model, where there is no restriction
$\alpha_0 \neq 0$ and $\beta_0 \neq 0$.

\begin{verbatim}
Call:
lm(formula = Z ~ probSize * E * de, subset = subE)

Residuals:
      Min        1Q    Median        3Q       Max 
-3.986320 -0.704510 -0.003003  0.688163  4.288311 

Coefficients:
                  Estimate Std. Error t value Pr(>|t|)    
(Intercept)      -0.009221   0.088655  -0.104 0.917164    
probSize          0.977059   1.560776   0.626 0.531324    
E4                0.122821   0.129093   0.951 0.341418    
E5                0.041659   0.126111   0.330 0.741152    
E6               -0.081222   0.129162  -0.629 0.529471    
E7                0.046731   0.125177   0.373 0.708921    
E8               -0.083599   0.129377  -0.646 0.518183    
E9               -0.016288   0.125517  -0.130 0.896752    
E10               0.086785   0.128709   0.674 0.500151    
de               -0.064338   0.362145  -0.178 0.858995    
probSize:E4     -10.772603   2.266890  -4.752 2.04e-06 ***
probSize:E5       0.054531   2.217975   0.025 0.980386    
probSize:E6     -10.292009   2.269018  -4.536 5.80e-06 ***
probSize:E7      -3.980349   2.204676  -1.805 0.071040 .  
probSize:E8      -6.119916   2.270653  -2.695 0.007046 ** 
probSize:E9      -1.849487   2.209092  -0.837 0.402491    
probSize:E10     -4.025803   2.264442  -1.778 0.075461 .  
probSize:de      33.518507   6.351755   5.277 1.34e-07 ***
E4:de             0.454335   0.506908   0.896 0.370121    
E5:de            -0.292063   0.518197  -0.564 0.573030    
E6:de             0.093784   0.505083   0.186 0.852699    
E7:de             0.204607   0.516213   0.396 0.691847    
E8:de             0.721077   0.506570   1.423 0.154637    
E9:de             0.196951   0.515310   0.382 0.702322    
E10:de           -0.038394   0.505668  -0.076 0.939478    
probSize:E4:de  -20.696043   8.826233  -2.345 0.019055 *  
probSize:E5:de   43.119875   9.068423   4.755 2.01e-06 ***
probSize:E6:de  -27.292071   8.811615  -3.097 0.001958 ** 
probSize:E7:de  -23.088609   9.031837  -2.556 0.010592 *  
probSize:E8:de  -29.097174   8.824071  -3.297 0.000979 ***
probSize:E9:de  -28.591599   9.022533  -3.169 0.001535 ** 
probSize:E10:de -29.477628   8.814444  -3.344 0.000828 ***
---
Signif. codes:  0 Ô***Õ 0.001 Ô**Õ 0.01 Ô*Õ 0.05 Ô.Õ 0.1 Ô Õ 1 

Residual standard error: 1.043 on 9964 degrees of freedom
Multiple R-Squared: 0.1751,	Adjusted R-squared: 0.1726 
F-statistic: 68.25 on 31 and 9964 DF,  p-value: < 2.2e-16 

\end{verbatim}

Every fitted term associated to $\alpha_0$ or $\beta_0$
is not statistically significant even at the 0.10 level.
In contrast, coefficients associated to $\alpha_1$ or $\beta_1$
terms are mostly significant.

The adjusted $R^2$ of the unrestricted model is actually
worse than the adjusted $R^2$ of the restricted model.
The standard analysis of variance table produced by the 
R software {\tt anova} command gives:
\begin{verbatim}
Model 1: Z ~ probSize * E * de - E * de - E - de - 1
Model 2: Z ~ probSize * E * de
  Res.Df   RSS   Df Sum of Sq      F Pr(>F)
1   9980 10843                             
2   9964 10830   16        13 0.7467 0.7475
\end{verbatim}

The improvement in variance explained using the 
unrestricted model is definitely not significant, as the
$P$ value exceeds 1/2.

This lack of significance justifies our finding
in the main text that weak universality holds,
at least for suites 3-10.

\subsection{Justification of scaling model with exponent $1/2$}

We fit models 
\begin{equation}
\label{satmodelb}
     \mu(\delta; N, E)  = \alpha(N,E) + \beta(N,E) (\delta - 1/2), \nonumber
\end{equation}
where
\begin{equation} \label{linmodelb}
      \alpha(N,E) = \alpha_0(E) +\alpha_1(E)/N^{1/2}, \qquad \beta(N,E) = \beta_0(E) + \beta_1(E)/N^{1/2}. \nonumber
\end{equation} 
in words we are saying that, within one suite, 
the means vary as a function of $\delta$ at a given sample
size, and that across sample sizes there is a root-n scaling 
of the means as a function of $N$.

The root-n scaling can be motivated both theoretically and empirically.

We first sketch a theoretical motivation.
Donoho \& Tanner (2008$b$) proved that, for $\delta$ fixed, the  
success probability for the Gaussian ensemble $p_1(A)\in(0,1)$ 
has a transition zone near the asymptotic phase transition $\rho(\delta;Q)$ 
of root-N width.  Namely, we showed that for $\rho$ below the asymptotic transition,
$p_1(A) \leq C \cdot \exp( - (\rho - \rho(\delta;Q))/w(\delta,N, Q) ) $,
where the width $ w \asymp N^{-1/2}$. 

From the viewpoint of probability theory, a width $w = O(N^{-1/2})$
would be typical when describing the frequency of success of an event among
$O(N)$ weakly dependent indicator variables.  Perhaps there is some such interpretation here;
although the rigorous proof in Donoho \& Tanner (2008$b$) does not provide one.
If this were the case, it would not be surprising for
the hypothesized underlying indicator variables 
to have success probabilities differing by order $1/\sqrt{N}$ at the non-Gaussian ensembles
from the corresponding ones at the Gaussian ensembles; such discrepancies
in limit theorems are common.
This would generate $N^{-1/2}$ scaling in the 
mean $Z$-scores between the Gaussian and other ensembles.
%

We now turn to  empirical motivation. We fit model (\ref{satmodelb})
to each suite and sample size separately, 
obtaining the full collection of coefficients
$\alpha(N,E)$ and $\beta(N,E)$ for $N= 200, 400,$ and $1600$.
We then fit an intercept-free  linear model to the collection of
fitted intercepts:
\begin{equation}\label{alphafit}
   \alpha(N,E) = \alpha_1(E,\gamma) N^{-\gamma} + error,
\end{equation}
and another to the slopes:
\begin{equation}\label{betafit}
   \beta(N,E) = \beta_1(E,\gamma) N^{-\gamma} + error.
\end{equation}
Table \ref{tab:gamma} presents the $R^2$ of the different fitted models.
Evidently exponent $\gamma = 1/2$ gives the best fit
both to intercepts and to slopes.  
\begin{table}
\begin{center}
\begin{tabular}{|l|l|l|l|l|l|l|l|}
\hline
         & \multicolumn{7}{c|}{power $\gamma$} \\
 \hline
 coef &    1.50     &  1.25  &     1.00   &   0.75  &     0.50   &    0.33     & 0.25 \\
 \hline   
intercept & 0.8477 & 0.8814& 0.9118 & 0.9482 & {\bf 0.9760} & 0.9533 & 0.9627 \\
slope      & 0.8413 & 0.8676 & 0.9122& 0.9343 & {\bf 0.9822} & 0.9626 & 0.9491 \\
\hline
 \end{tabular}
 \caption{$R^2$ of fits in models (\ref{alphafit}) and (\ref{betafit}). 
 The exponent $\gamma=1/2$ gives the best fit in both cases}\label{tab:gamma}
\end{center}
\end{table}
 
 We don't see an easy way to attach statistical
 significance to this finding, using standard software.
 
This fitting exercise also shows that 
the exponent $\gamma=1/2$ fits well
enough to make the case for
$\alpha_0 \neq 0$ and $\beta_0 \neq 0$
very weak.
We  fit two joint models to two datasets,
one containing all the fitted intercepts $\hat{\alpha}(N,E)$
and the other containing all the fitted slopes  $\hat{\beta}(N,E)$
from the same collection of ensembles and 
the standard problem sizes $N = 200,400,1600$.

The first joint model for  absolute intercepts included terms
that do not vanish as $N$ grows large
 \begin{equation} \label{linmodel-int}
      |\hat{\alpha}(N,E)| = \gamma_0(E) +\gamma_1(E)/N^{1/2}  + error.
\end{equation} 
The second joint model did not include such terms:
\begin{equation} \label{linmodel-noint}
     |\hat{\alpha}(N,E)| = \gamma_1(E)/N^{1/2}  + error.
\end{equation} 
The traditional t-statistics associated
with $\gamma_0$ terms were nonsignificant
with one exception; this must be treated as 
unimpressive owing to multiple comparison effects.
In contrast, the majority of $\alpha_1$ terms were
significant.
The analysis of variance comparing the two fits
gave an F statistic of $1.53$ on $12$ and $23$ degrees
of freedom with a nominal $P$-value of $0.2394$.
In short the larger model which includes a constant term
independent of $N$ explains little.
Although this cannot be used with the usual interpretation
it does show that any tendency to not vanish must be very weak.

The first joint model for slopes included intercept terms:
  \begin{equation} \label{linmodel-slp}
       \hat{\beta}(N,E) = \beta_0(E) + \beta_1(E)/N^{1/2} + error
\end{equation} 
The second joint model did not include intercept terms:
\begin{equation} \label{linmodel-noslp}
    \hat{\beta}(N,E) =  \beta_1(E)/N^{1/2} + error
\end{equation} 
The traditional t-statistics associated
with $\beta_0$ terms were nonsignificant
with one exception; this again must treated as 
unimpressive owing to multiple comparison effects.
In contrast, the majority of $\beta_1$ terms were
significant.
The analysis of variance comparing the two fits
gave an F statistic of $1.35$ on $12$ and $23$ degrees
of freedom with a nominal $P$-value of $0.30$.
Again the larger model which includes a constant term
independent of $N$ explains little.
Although this cannot be used with the usual interpretation
it does show that any tendency to not vanish must be very weak.

\subsection{Exceptional suites}\label{sec:exceptional}

\begin{figure}
\begin{center}
\includegraphics[width=4in]{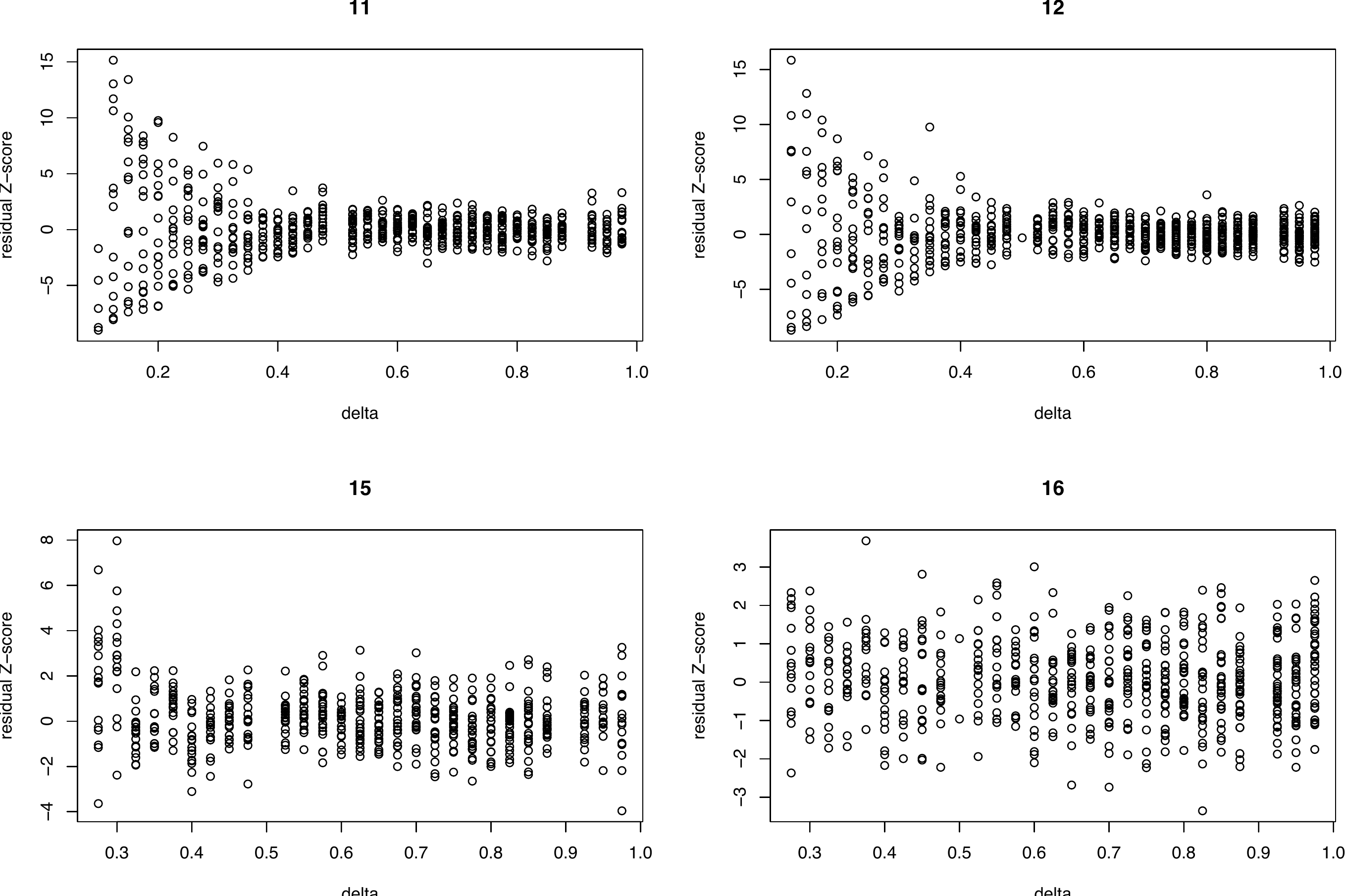}
\caption{Residuals from scaling hinged fit, $N=200$.
Panels: Suites 11,12,15,16.
Note the apparent increase in variance of residuals at small $\delta$}
\label{fig-HingedResidualsVsDelta-200}
\end{center}
\end{figure}

\begin{figure}
\begin{center}
\includegraphics[width=4in]{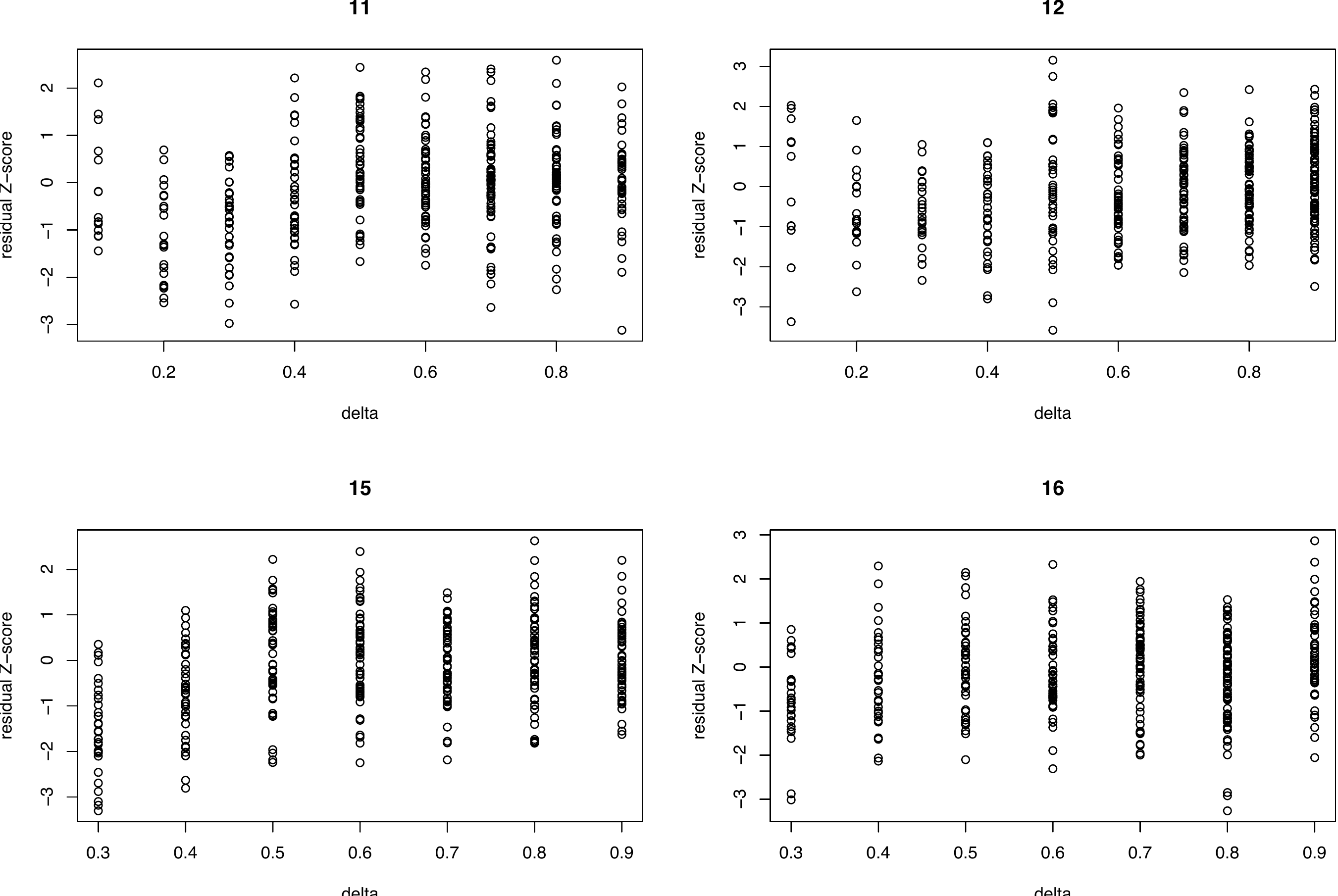}
\caption{Residuals from scaling hinged fit, $N=1600$.
Panels: Suites 11,12,15,16.
Note the apparent curvilinearity in mean residuals at small $\delta$.
}
\label{fig-HingedResidualsVsDelta-1600}
\end{center}
\end{figure}

In \S \ref{LinearFits} we excluded  suites 11,12,15,16 from analysis.
Our decision was based on
the fact that three of these  -- 11,12,15 -- were to the naked eye
quite different than the others, by two criteria:
\bitem
 \item Behaviour of the bulk distribution of $Z$-scores 
 (see figures \ref{fig-Bulk-Zscores-N-200}-\ref{fig-Bulk-Zscores-N-1600}; 
panels in the lower row)
 \item Behaviour of the plots of $Z$-scores  versus $\delta$.
\eitem 
We grouped suite 16 with these based on the fact that the
underlying matrix ensemble was the same as suite 15;
they differ only in properties of the solution coefficients.

\subsubsection{Bulk distribution of $Z$-scores}

We remark that although the naked eye perceives
differences in the $Z$-scores for $\delta$ small, these suites
are overall consistent with weak universality.
As $N$ increases, at each fixed $\delta$, in each suite,
the distribution of $Z$-scores gets closer to
$N(0,1)$. For example, the evident discrepancies
that exist for small $\delta$ in the QQ Plots of $Z$-scores at $N=200$
are not apparent for large $\delta$; moreover, the discrepancies 
are dramatically attenuated by $N=1600$.

\subsubsection{Modelling of $Z$-scores}

A major part of the exceptional behaviour 
of these suites is the apparent nonlinear dependence
of mean $Z$-score on $\delta$.
We have seen that for suites 3-10,  an adequate description
is provided by linear dependence on $\delta$.
However, for the exceptional suites
this is no longer the case.
Plots of residual $Z$-scores versus $\delta$
show that the exceptional suites exhibit nonlinear dependence
on $\delta$; typically downward sloping at $\delta < .5$ and
no sloping or upward sloping at $\delta > .5$.

Mild nonlinearity in $\mu$ can be modelled through hinged models:
 \[
   \mu(N,E) = \alpha(N,E) + \beta(N,E)(\delta - 1/2) + \kappa(N,E) (1/2-\delta)_+ .
 \]
 These model $\delta$ dependence as a linear spline 
 with knot at $\delta = 1/2$.  As before we can model $\kappa$'s dependence
 on $N$ in the same way as we have done with $\alpha$ and $\beta$:
 \[
     \kappa(N,E) = \kappa_0(E) + \kappa_1(E)/\sqrt{N}.
 \]
 
 Such hinged models are to be preferred over the linear models
 in the 4 exceptional suites.
 
 We first present the R session transcript of the linear 
fit with $\sqrt{N}$ scaling terms but no constant terms:
  \[
   \mu(N,E) = \alpha_1(E)/\sqrt{N} + \beta_1(E)/\sqrt{N} \cdot (\delta - 1/2) 
  \]

 \begin{verbatim}
 Call:
lm(formula = Z ~ probSize * E * de - E * de - E - de - 1, subset = subX)

Residuals:
    Min      1Q  Median      3Q     Max 
-6.3990 -1.0153 -0.1640  0.6940 18.1441 

Coefficients: (1 not defined because of singularities)
                Estimate Std. Error t value Pr(>|t|)    
probSize           2.127      1.262   1.685    0.092 .  
probSize:E11      34.828      1.584  21.984   <2e-16 ***
probSize:E12      46.890      1.607  29.182   <2e-16 ***
probSize:E15      17.804      1.743  10.214   <2e-16 ***
probSize:E16          NA         NA      NA       NA    
probSize:de     -138.228      3.848 -35.923   <2e-16 ***
probSize:E12:de    3.747      5.274   0.710    0.477    
probSize:E15:de  159.013      6.410  24.806   <2e-16 ***
probSize:E16:de  128.493      6.065  21.186   <2e-16 ***
---
Signif. codes:  0 Ô***Õ 0.001 Ô**Õ 0.01 Ô*Õ 0.05 Ô.Õ 0.1 Ô Õ 1 

Residual standard error: 1.854 on 4736 degrees of freedom
Multiple R-Squared: 0.5366,	Adjusted R-squared: 0.5358 
F-statistic: 685.5 on 8 and 4736 DF,  p-value: < 2.2e-16 
\end{verbatim}

In this fit, there is substantial lack of fit: the standard deviation of
residuals, 1.854, is dramatically larger than 1.0 (the fit for
suites 3-10 gave instead a residual standard deviation close to 1).

We fit the hinged model with  $\sqrt{N}$ scaling:
 \[
   \mu(N,E) = \alpha_1(E)/\sqrt{N} + \beta_1(E)/\sqrt{N} \cdot (\delta - 1/2) 
                                                           + \kappa_1(E)/\sqrt{N}  \cdot (1/2-\delta)_+ 
 \]
 The R transcript follows:

\begin{verbatim}
Call:
lm(formula = Z ~ probSize * E * de + probSize * E * dea2 - E * 
    de - E * dea2 - 1, subset = subX)

Residuals:
     Min       1Q   Median       3Q      Max 
-9.00131 -0.86664 -0.07468  0.71868 15.83998 

Coefficients: (1 not defined because of singularities)
                  Estimate Std. Error t value Pr(>|t|)    
probSize            -2.497      1.904  -1.311 0.189834    
probSize:E11         8.493      2.557   3.322 0.000901 ***
probSize:E12        24.982      2.597   9.619  < 2e-16 ***
probSize:E15         7.455      2.677   2.785 0.005369 ** 
probSize:E16            NA         NA      NA       NA    
probSize:de        -10.276      7.000  -1.468 0.142203    
probSize:dea2      301.622     14.237  21.186  < 2e-16 ***
probSize:E12:de    -31.004      9.391  -3.302 0.000969 ***
probSize:E15:de     87.995     10.234   8.599  < 2e-16 ***
probSize:E16:de     15.900      9.630   1.651 0.098777 .  
probSize:E12:dea2  -50.774     20.168  -2.518 0.011850 *  
probSize:E15:dea2  -80.204     26.626  -3.012 0.002606 ** 
probSize:E16:dea2 -232.594     26.649  -8.728  < 2e-16 ***
---
Signif. codes:  0 Ô***Õ 0.001 Ô**Õ 0.01 Ô*Õ 0.05 Ô.Õ 0.1 Ô Õ 1 

Residual standard error: 1.706 on 4732 degrees of freedom
Multiple R-Squared: 0.6081,	Adjusted R-squared: 0.6071 
F-statistic: 611.9 on 12 and 4732 DF,  p-value: < 2.2e-16 
\end{verbatim}

Here {\tt dea2} denotes $(1/2 - \delta)_{+}$,
so $\kappa$ terms are associated with {\tt dea2},  in an impromptu notation:
$\kappa(N,E15) = $ {\tt probSize:E15:dea2} +  {\tt probSize:dea2 } /$\sqrt{N}$.

The key points to observe here are: (1) the individual coefficients 
associated with  hinge terms are significant; and (2) the adjusted $R^2$ (0.6071)
is substantially higher than it was for a linear fit (0.5358).  The analysis
of variance comparing the hinged model with the linear model
gives an $F$ statistic of 113 on 4720 and 8 DF, which is wildly significant;
see the transcript:
\begin{verbatim}
Model 1: Z ~ probSize * E * de + probSize * E * dea2
Model 2: Z ~ probSize * E * de
  Res.Df     RSS   Df Sum of Sq      F    Pr(>F)    
1   4720 13466.6                                    
2   4728 16055.4   -8   -2588.8 113.42 < 2.2e-16 ***
\end{verbatim}

In the above fits we considered only models imposing $N^{-1/2}$ scaling
on $\mu$.  Allowing terms which do not decay in $N$ does not
improve the fit. The following transcript shows  fits of a
a model for  mean $Z$-score in
suite $E$ containing non-scaling terms: $\mu(N,E) = \mu_0(N,E) + \mu_1(N,E)/\sqrt{N}$;
the fits shown in the paragraphs immediately above correspond 
instead to restrictions $\mu_0 = 0$.
\begin{verbatim}
Call:
lm(formula = Z ~ probSize * E * de + probSize * E * dea2, subset = subX)

Residuals:
     Min       1Q   Median       3Q      Max 
-9.28404 -0.79999 -0.02881  0.74525 15.75879 

Coefficients:
                    Estimate Std. Error t value Pr(>|t|)    
(Intercept)          0.03677    0.26074   0.141 0.887848    
probSize             5.25575    4.68185   1.123 0.261674    
E12                 -0.68992    0.37841  -1.823 0.068333 .  
E15                  0.13874    0.38490   0.360 0.718529    
E16                 -0.17611    0.38956  -0.452 0.651242    
de                   0.05261    1.12637   0.047 0.962751    
dea2                -4.81687    2.26618  -2.126 0.033593 *  
probSize:E12        27.98270    6.76604   4.136  3.6e-05 ***
probSize:E15        -3.43537    6.93196  -0.496 0.620211    
probSize:E16        -5.42887    6.98430  -0.777 0.437023    
probSize:de        -10.66394   19.77908  -0.539 0.589807    
E12:de               2.13317    1.52864   1.395 0.162940    
E15:de              -0.49034    1.61080  -0.304 0.760830    
E16:de               0.14309    1.59771   0.090 0.928643    
probSize:dea2      380.64608   40.15413   9.480  < 2e-16 ***
E12:dea2             6.96907    3.23697   2.153 0.031372 *  
E15:dea2           -13.34798    4.12398  -3.237 0.001218 ** 
E16:dea2            -0.58285    4.25568  -0.137 0.891070    
probSize:E12:de    -66.26910   26.69975  -2.482 0.013099 *  
probSize:E15:de     96.38232   28.49783   3.382 0.000725 ***
probSize:E16:de     12.93573   27.73609   0.466 0.640961    
probSize:E12:dea2 -165.95925   57.18373  -2.902 0.003723 ** 
probSize:E15:dea2  141.75943   73.66213   1.924 0.054358 .  
probSize:E16:dea2 -225.04810   74.88419  -3.005 0.002667 ** 
---
Signif. codes:  0 Ô***Õ 0.001 Ô**Õ 0.01 Ô*Õ 0.05 Ô.Õ 0.1 Ô Õ 1 

Residual standard error: 1.689 on 4720 degrees of freedom
Multiple R-Squared: 0.5418,	Adjusted R-squared: 0.5395 
F-statistic: 242.6 on 23 and 4720 DF,  p-value: < 2.2e-16 
\end{verbatim}

The key points to observe are: (1) The standard error of residuals
is not meaningfully improved by allowing the extra explanatory terms:
it drops from 1.706 for the $N^{-1/2}$ scaling model to 
1.689 for the full model; and (2) The adjusted $R^2$ is worse for the 
full model than it is for the scaling model.
At the level of individual effects, the bulk of the 
non-scaling terms are not significant, and 
the scaling hinge terms remain significant.
We view the few significant non-scaling terms as
possibly caused by the significant modeling error that still remains:
i.e. since the residual standard error, at about 1.7,
is about 70\% higher than it would be if 
everything were explained in a satisfactory way.

The next table summarizes the residuals from the
hinged scaling model, grouped by suite
and $N$.  The summaries include group
means, standard deviations, medians and median 
absolute values.

\begin{table}
\begin{center}
\begin{tabular}{|l|l|l|l|l|ll|}
\hline
E & N & cases & mean & sd & med & mav \\ 
\hline
11 & 200 & 654  & 0.0713 & 2.58 & -0.0432 & 0.942 \\
11 & 400 & 214  & -0.210 & 1.87 & -0.230 & 0.863 \\
11 & 1600 & 364  & -0.115 & 1.06 & -0.0651 & 0.625 \\
12 & 200 & 694  & 0.0528 & 2.39 & 0.055 & 0.92 \\
12 & 400 & 220  & -0.114 & 1.96 & -0.293 & 0.871 \\
12 & 1600 & 383  & -0.139 & 1.10 & -0.181 & 0.798 \\
15 & 200 & 546  & 0.134 & 1.32 & 0.0457 & 0.723 \\
15 & 400 & 191  & -0.38 & 1.30 & -0.349 & 0.868 \\
15 & 1600 & 337  & -0.185 & 1.08 & -0.213 & 0.734 \\
16 & 200 & 611  & 0.0615 & 1.07 & 0.0220 & 0.718 \\
16 & 400 & 196  & -0.176 & 1.01 & -0.257 & 0.681 \\
16 & 1600 & 334  & -0.112 & 1.02 & -0.132 & 0.66 \\
\hline
\end{tabular}
\caption{Group statistics of residuals from the hinged scaling model,
grouped by sample size}
\end{center}
\end{table}

\begin{figure}
\begin{center}
\includegraphics[width=4ins]{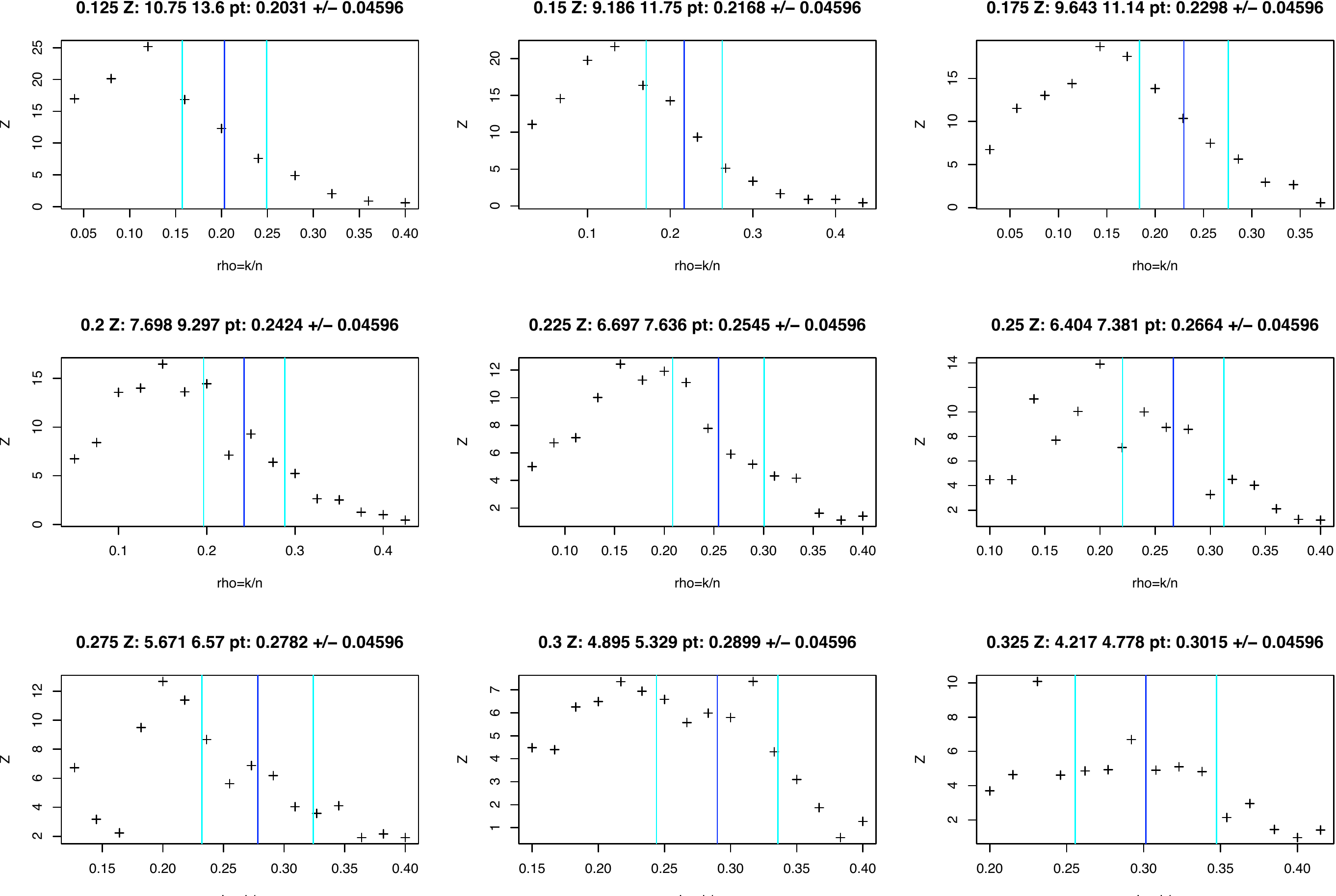}
\caption{$Z$-scores for suite 12, $N=200$ plotted versus $\rho = k/n$ at
each value of $\delta=n/N$. Vertical bars indicate asymptotic phase
transition and nominal width of the transition zone.}
\label{fig-ZscoresVsRhoWPTDeco-200}
\end{center}
\end{figure}

\begin{figure}
\begin{center}
\includegraphics[width=4in]{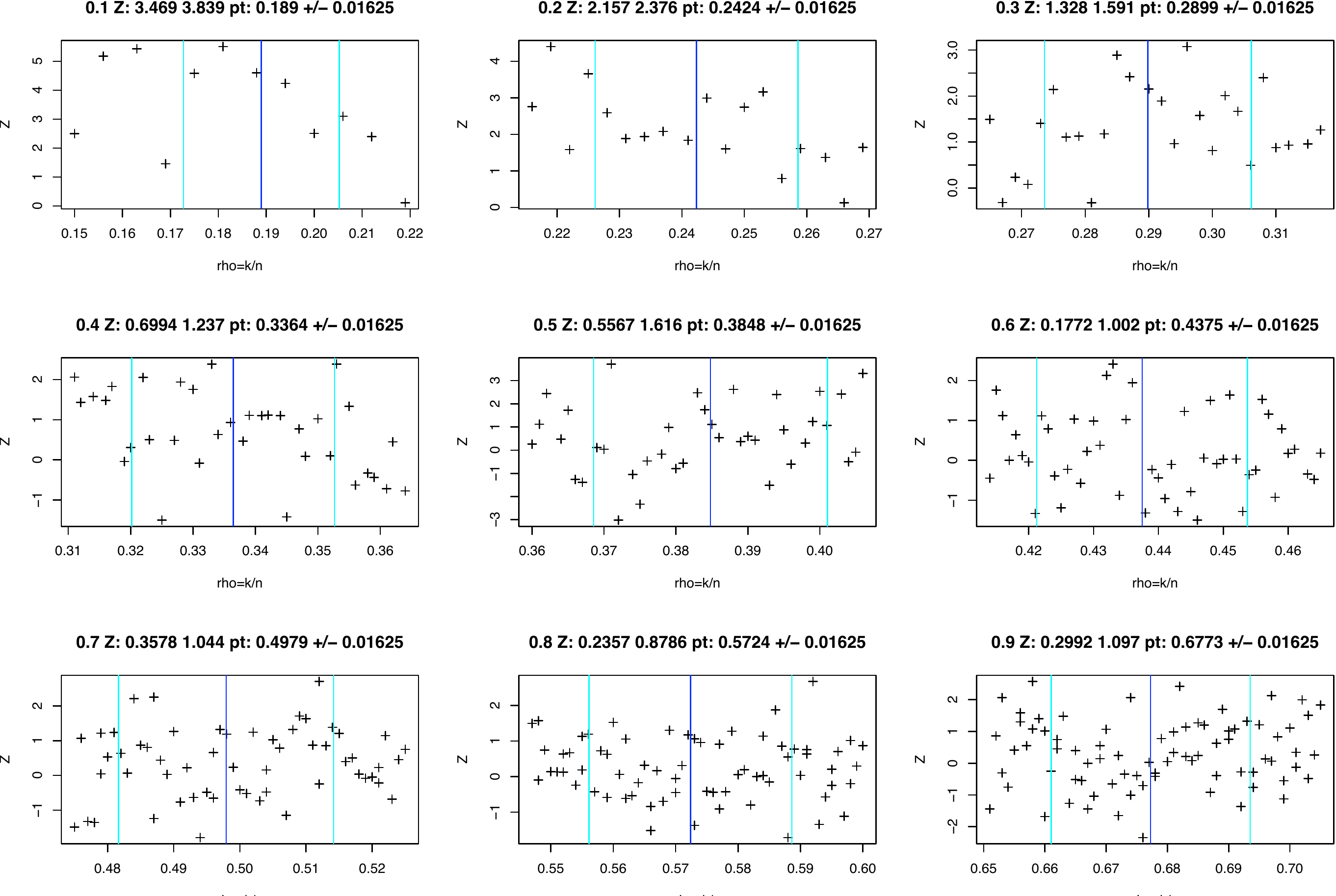}
\caption{$Z$-scores for suite 12, $N=1600$ plotted versus $\rho = k/n$ at
each value of $\delta=n/N$. Vertical bars indicate asymptotic phase
transition and nominal width of the transition zone.}
\label{fig-ZscoresVsRhoWPTDeco-1600}
\end{center}
\end{figure}

Suite 16 is adequately explained,
since the standard deviation is fairly close to 1.00
at each $N$. 

In all suites the residual standard deviation is decreasing 
with $N$, and in every case the residual standard
deviation is less than 1.10 at $N=1600$.
In our view this, together with the 
structure of the scaling model,
merits  the conclusion  that {\it these
random matrix ensembles agree with the Gaussian ensemble
for large $n$.}

However, there {\it is} noticeable structure at small $n$.
The key pattern visible in the table is the tendency of residuals
to be positive at small $n$. 


Some further structure  is visible 
 in the $\delta$ variable;
 figures \ref{fig-HingedResidualsVsDelta-200}
 and \ref{fig-HingedResidualsVsDelta-1600}
 show that at $N=200$ there is a dramatic 'blowup'
 in variance at small $\delta$, which has largely
 disappeared at the larger problem size $N=1600$.
 They also show that suite 16 is already well-behaved at $N=200$,
 and the "variance blowup" is largely a phenomenon of suites 11 and 12.
 Finally, at $N=1600$ one can see indications of curvilinear structure,
 so evidently the hinged model can only be regarded as
 an approximation.

As it turns out, the "variance blowup" is not a variance
phenomenon at all; it is caused by significant unmodeled
structure in the means as a function of $\rho=k/n$.
Figure \ref{fig-ZscoresVsRhoWPTDeco-200} shows the $Z$-scores at
suite 12, for problem size $N=200$.  The $Z$-scores are
quite large and nonrandom, indicating significant unmodelled
structure.

On the other hand, this unmodeled structure is largely a
phenomenon of small problem size.
Figure \ref{fig-ZscoresVsRhoWPTDeco-1600} shows the $Z$-scores at
suite 12, now for problem size $N=1600$.  The $Z$-scores are
not so large and much more random, indicating that unmodelled
structure is less of a problem, except at small $\delta$.

 The unmodelled structure can be largely accounted for
 as a {\it displacement of the $LD50$} by about 1 transition width at $N=200$.
 The asymptotic phase transition is at $\rho(\delta;Q)$;
 if we assume that in suite 12 the $LD50$ is not
 at  $\rho(\cdot;Q)$ but instead at $\rho^\dagger = \rho(\delta;Q) - w(\delta;N,Q)$,
 and we assume a probit link,
 we can expect the $Z$-scores to have a mean
 \[
     \mu \propto   \phi( (\rho - \rho^\dagger)/w)/w,
 \]
 where $\phi$ denotes the standard Gaussian density.

\begin{figure}
\begin{center}
 \includegraphics[width=4in]{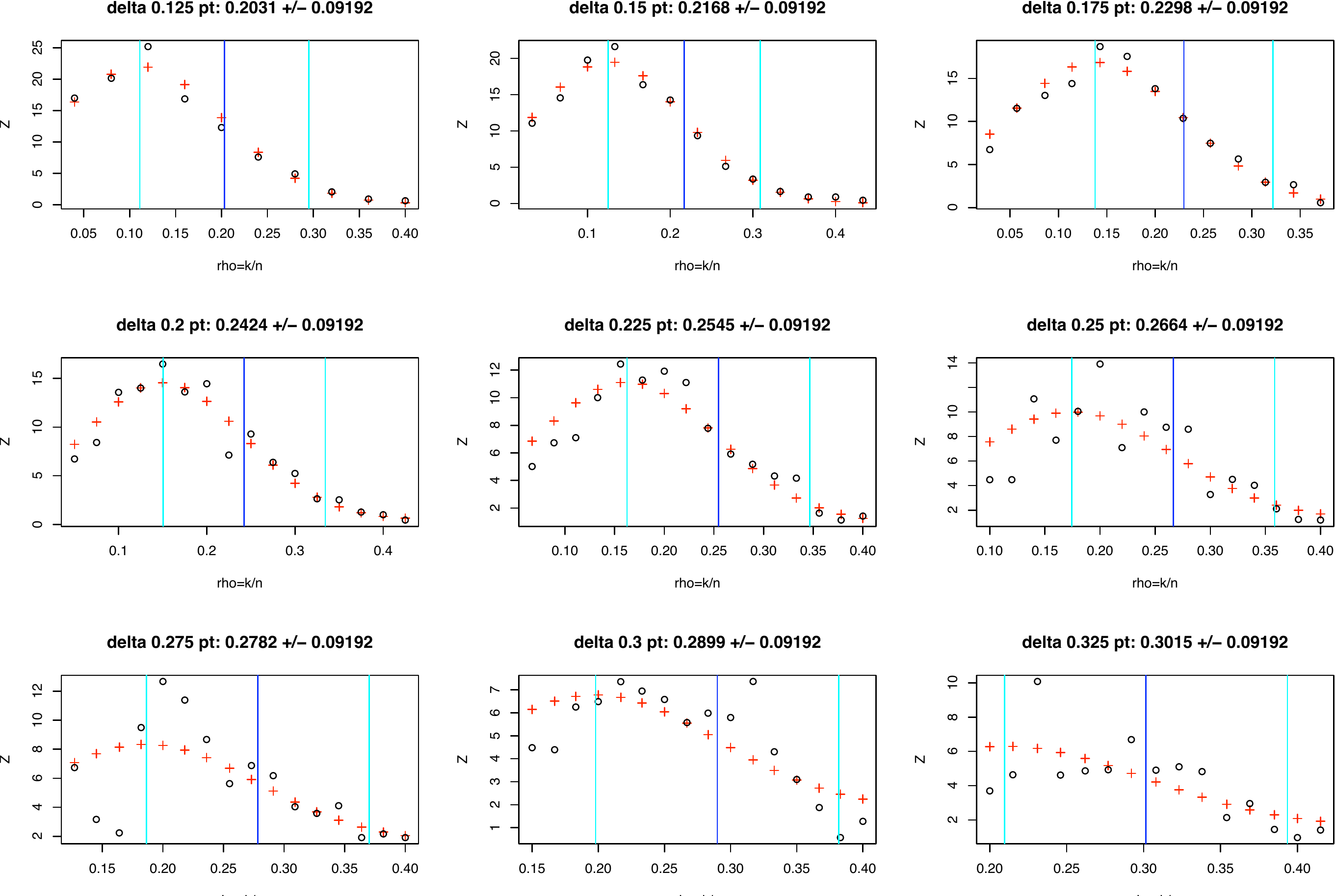}
 \caption{$Z$-scores for suite 12 plotted versus $\rho = k/n$ at
 9 different values of $\delta=n/N$, together with fitted mean shift model.
 Red Crosses portray the displaced $LD50$ model.
 }
 \label{fig-ZvsAdjustPT}
\end{center}
 \end{figure}

 Figure \ref{fig-ZvsAdjustPT} fits this model to the $Z$-scores at
 suite 12, for sample size $N=200$. It accounts for most of the
 structure in the $Z$-scores at these values of $\delta$.

\subsubsection{Understanding the exceptional ensembles}

When $N$ is
small and $\delta$ is small, the exceptional ensembles have
$LD50$'s  displaced substantially below the asymptotic phase transition,
It also seems that if we keep $\delta$ fixed and increase
$N$, such displacement eventually disappears.

It seems to us that one can understand this phenomenon
as a transient effect of random sampling in discrete structures.
To explain this, note that when "success" is declared in
our computations, this means that the image of a specific $k$-face $F$
of the polytope $Q$ `survives' projection under $A$, i.e. $AF$
is a face of $AQ$. (Here of course $Q$ is either $T^{N-1}$ or $C^N$,
see the main text; and the preceding statement
deserves and requires rigorous proof; it receives such proof in our other papers). 
In particular that there is no pair $(x_0, x_1)$ both in $F$
with $Ax_0 = Ax_1$. 
But then there is no vector $v$ supported on those $k+1$ sites
with $Av = 0$.

On the other hand, if there exist vectors supported on $k+1$ sites that
belong to this null space, then there definitely will be $k$-faces $F$
that get lost under projection. 

Let $E_{k,n,N}$ denote the
event that there exist $k+1$-sparse vectors in the null
space of $A$. Then when event $E_{k,n,N}$ occurs
there will be a $k$-face of $Q$ that does not survive projection.

For a Gaussian random matrix $A$, with probability one
the columns of $A$ are in general position,  and so $P(E_{k,n,N}) = 0$
provided $k < n$.  However, for the exceptional ensembles
$P(E_{k,n,N}) > 0$ even for quite small $k$.
 

The Ternary (1/10) and Expander (1/15) ensembles generate highly sparse
matrices $A$.   Particularly when $n$ is small, we have observed that these can 
have sparse solutions in their 
null space; that is, $E_{k,n,N}$ can have substantial probability.

It may help some readers to know we are 
effectively discussing  the sparsity of the {\it sparsest}
vector in the null space of $A$; this is often denoted $spark(A)$ (Donoho \& Elad, 2003).
For Gaussian random matrices $A$, $spark(A) = n$,
while for the Ternary and Expander ensembles, we find that
$spark(A)$ has a good chance of being much smaller than $n$.
 
Dossal {\it et al.} (2009) have proposed a greedy algorithm which 
can be used to find sparse vectors in the null space of a matrix.
It provides an heuristic upper bound $spark^+(A)$ on the
spark of a matrix .
Applying this algorithm to random draws of Ternary (1/10) and Expander (1/15) 
we observe sparse vectors satisfying $Az=0$; ranges of these observed sparsity 
levels are recorded in columns 2 and 3 of table \ref{tab:rip}.  

The sparsity levels recorded in columns 2 and 3 of table \ref{tab:rip} 
 indicate clearly that $E_{k,n,N}$
 can occur for quite small $k$, but they do not 
indicate the size of $P(E_{k,n,N})$
or give any other information about the
prevalence of these sparse null space vectors.  
Columns 4-9 of table \ref{tab:rip} indicate the empirical bound $spark^+(A)$ for 
suites 1-2,11-12,15-16 with $S/M$ closest to 50\%.  For Ternary (1/10) with 
$N=200$ and $n\le 60$ it is common to draw matrices which contain a column 
of all zeros, allowing for 1-sparse vectors in the null space.
In such cases, $spark(A)=1$, which is dramatically smaller than $n$,
what we would see with the Gaussian ensemble.

With such  easy counterexamples
to uniqueness it now seems remarkable that the
typical case, as measured by the observed $LD50$, is as high 
as we observed.  This may be explained by the fact that
typically, the bad  $k$-sets do not intersect
the particular $k$-set we are interested in 
(i.e. the one supporting the solution of interest).

For the Expander ensemble with $N=200$ and $n\le 40$ there are vectors in the 
null space with only 6 nonzeros, and null space vectors with more than 6 
nonzeros become common; as a consequence, it becomes increasingly common 
that of the $M=200$ problem instances presented to \p\ and \lp\ for each 
value of $k$ moderately larger than 6 there will be vectors $x_1$ which are 
sparser than $x_0$ and yet which satisfy $Ax_1=Ax_0=y$.  

We observe that 
the $LD50$ for suites 15-16 is significantly displaced below the asymptotic value for this range of 
$(N,n)$.  Suites 15-16 are observed to have similar values of $LD50$ for 
this range of $(N,n)$. Although these $LD50$s are similarly displaced, 
the effect is more noticeable in the $Z$-scores of suite 15 due to its 
being compared with suite 1 which has a markedly higher $LD50$ than does 
suite 2 (which suite 16 is compared against).

\begin{table}
\begin{center}
\begin{tabular}{|c||c|c||c|c|c|c|c|c|}
\hline
$(N,n)$ & Ternary (1/10) & Expander & 1 & 2 & 11 & 12 & 15 & 16 \\
\hline
(200,20) & 1 & 1 & 5 & 4 & 2 & 1 & 2 & 1 \\
\hline
(200,30) & 1 & 5-14 & 9 & 7 & 5 & 4 & 3 & 3 \\
\hline
(200,40) & 1 & 6-20 & 13 & 10 & 11 & 8 & 9 & 8 \\
\hline
(200,50) & 1 & 19-34 & 19 & 14 & 16 & 12 & 17 & 13 \\
\hline
(200,60) & 1-61 & 16-41 & 25 & 18 & 23 & 16 & 22 & 18 \\
\hline
(1600,160) & 1-61 & 1-61 & 39 & 30 & 38 & 29 & 40 & 31 \\
\hline
\end{tabular}
\end{center}
\caption{{\it Sparse vectors in $Nullspace(A)$}. 
Columns 2 and 3 of table \ref{tab:rip} indicate minimal 
sparsity levels of vectors $z$ obeying
$Az=0$ with $A$ drawn from ensembles Ternary (1/10) 
and Expander (1/15).  Columns 4-9 indicate the value 
of the empirical upper bound $spark^+(A)$ for 
suites 1-2,11-12,15-16 with $S/M$ closest to 50\%.}\label{tab:rip}
\end{table}

If we keep $n/N$ and $k/n$ constant, and $k/n$ small, but let $N$ increase
it seems that
the probability $P(E_{k,n,N})$ tends to zero rapidly. 
For $N=1600$ and
$\delta=0.1$, the observed $LD50$ values for these exceptional ensembles are 
within 1 of those observed for the Gaussian suites.  We also observe that 
for $N=1600$ and $n=160$ we are unable to find null space vectors with 
fewer than $161$ nonzeros.  This suggests the following interpretation of
what is happening: for each $(\delta,\rho)$ there  is a transient effect, such that
as $N$ increases, initially it is most likely that
 $spark(A)$ is less than $\rho \cdot n$, but for sufficiently
 large $N$, it becomes likely that $spark(A)  \approx n$.
 
 Under our interpretation, the driving effect
 behind the exceptional ensembles
is the phenomenon that for small $N$
the matrix $A$ has its columns not in general position;
and, for fixed $\delta = n/N$
the chance of observing this phenomenon decays rapidly with increasing
problem size $N$.

\subsection{Validation ensembles: Rademacher and Hadamard}\label{sec:reserve}

The computer runs reported here took place over a
substantial interval of calendar time.  As it happens,
two of the matrix ensembles -- Rademacher and partial Hadamard --
were completed long after the others and the data were not
available at the time of the analyses reported so far.

A basic principle in science is out-of-sample validation, namely
building a model using data available at a certain moment in time and then later using
data that were not available to the model construction to test the model construction.

The fresh data provided by the Rademacher and partial Hadamard
runs provide an excellent opportunity for validation of our fitted models.

Such validation could be helpful for the following reason.
A critic of our analysis could point out that we have used
{\it the same data} to choose the decay exponent $\gamma = 1/2$
in our scaling model $\mu \sim N^{-1/2}$, as well as to 
study lack of fit of that model.
This is not inferentially rigorous, as the exponent has specifically
been chosen to minimize lack of fit, and we then argue that
the lack of fit at such $\gamma$ is not significant.
The use of the same data for both tasks means
that the $t$-scores and $p$-values no longer have
the assumed distributions.

In our case, we have about 10,000 $Z$-scores, and we believe
this criticism is not as serious as in many other occurrences
of this practice, and we are
prepared to argue this point.
However, in this instance 
because we have an independent set of validation data,
we can actually perform direct model validation
and avoid lengthy rationalizations of our earlier analysis.

The validation data have one oddity: The Hadamard ensemble
is computable only for special choices of $N$, and we happen only
to have data for $N=256$ and $N=512$, while for the Rademacher
ensemble we have data for $N=200$, $400$, and $1600$ as
per usual.

We now use those data to validate the following points in our analysis.

\bitem
 \item Bulk distribution of $Z$-scores;
 \item Fit of $N^{-1/2}$ model for scaling of means;
 \item Lack of Fit at small $\delta$ caused by displacement of the LD50.
\eitem

It will turn out that the Rademacher ensemble
is not exceptional while the Hadamard is;
and for both Rademacher and Hadamard
the models developed using the other ensembles
fit rather well.
 
\subsubsection{Bulk distribution of $Z$-scores}
\begin{figure}
\begin{center}
\begin{tabular}{c}
\includegraphics[width=3in]{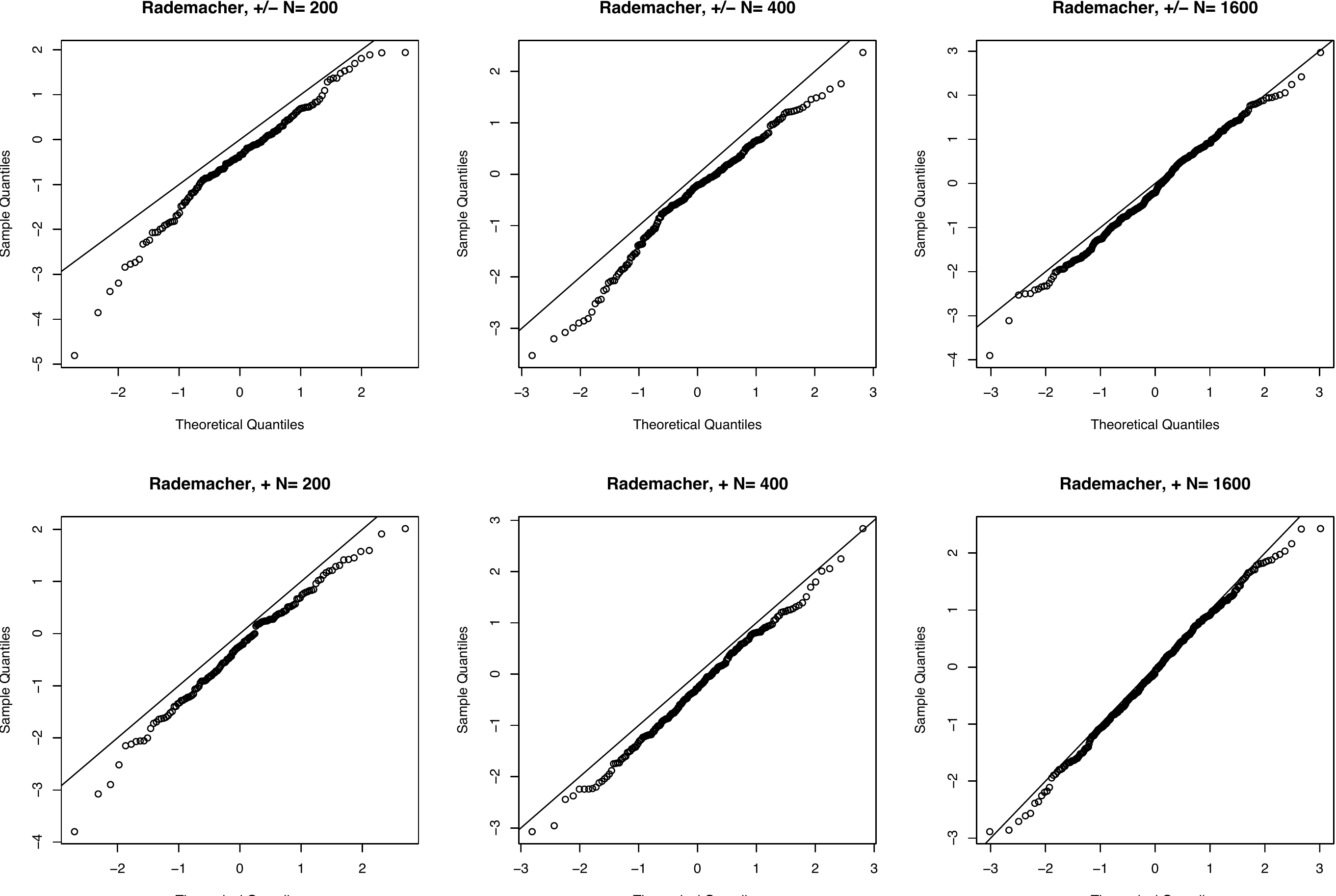} \\
\end{tabular}
\caption{$Z$-scores for suites 19-20, the Rademacher ensemble, at $N=200$, $400$ and $1600$}
\label{fig-Rademacher-Zscores}
\end{center}
\end{figure}

\begin{figure}
\begin{center}
\begin{tabular}{c}
\includegraphics[width=3in]{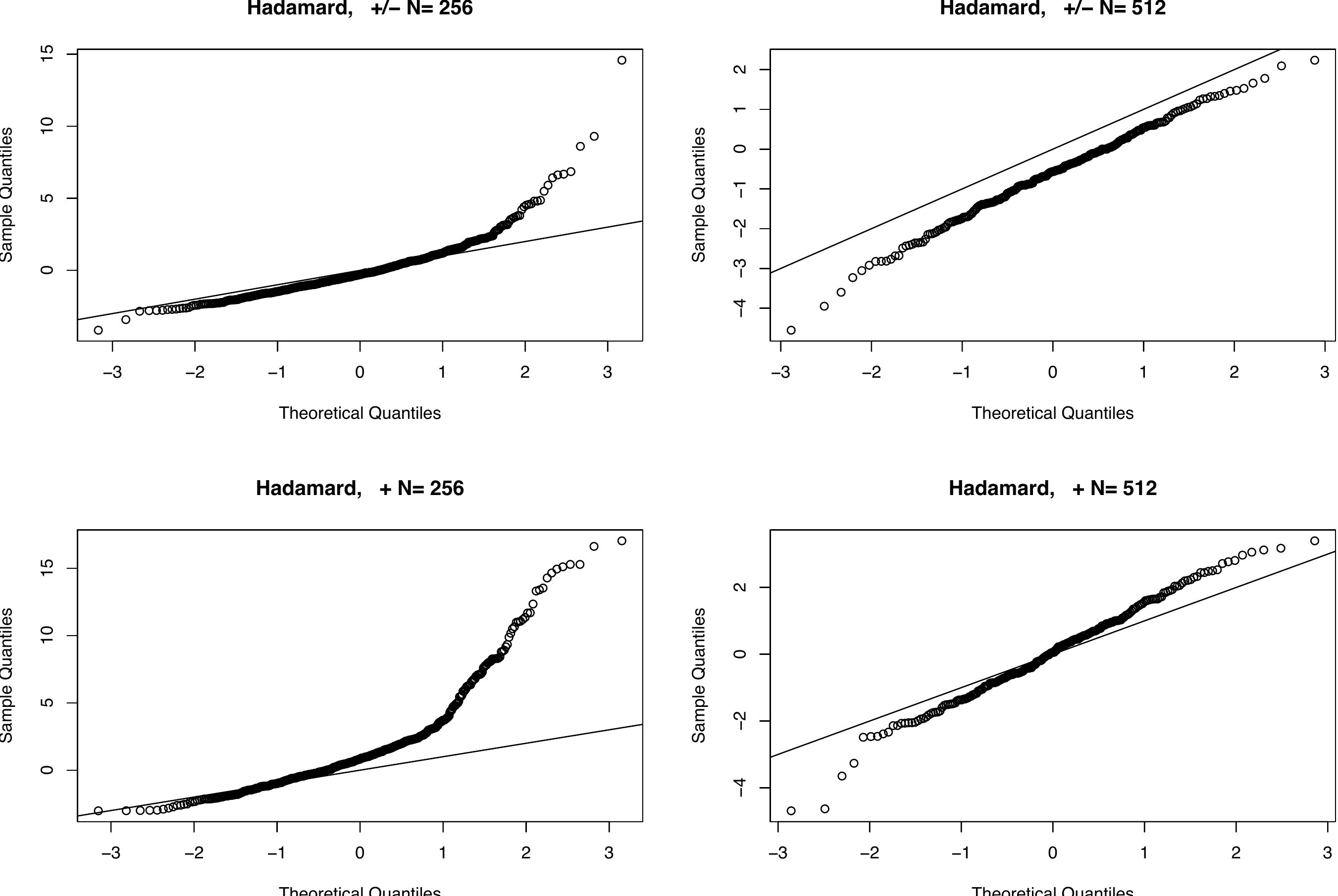} \\
\end{tabular}
\caption{$Z$-scores for suites 13-14, the Hadamard ensemble, at $N=256$ and $512$}
\label{fig-Hadamard-Zscores}
\end{center}
\end{figure}

In figures \ref{fig-Rademacher-Zscores} and \ref{fig-Hadamard-Zscores}
we present a sequence of  QQ plots showing the bulk distribution
of $Z$-scores for suites $19-20$ (at $N=200$, $400$ and $1600$) 
and for suites 13-14 (at $N=256$ and $512$).
In each plot the identity line is also displayed; if the $Z$-scores
had truly a standard normal distribution, they would oscillate around this line.

The Rademacher suite $Z$-scores are typical of what we have already seen in
suites 3-10.  Recall that $Z$-scores are expected to cluster around the line
$Y=X$ at the null hypothesis.
There is some noticeable bulk misfit at $N=200$, less
at $N=400$, and by $N=1600$ the visual misfit has mostly disappeared.
The misfit at $N=200$ is mostly in the lower tail,  and in a downward shift
away from the $Y=X$ line, meaning that 
there are situations where our realization of the Rademacher ensemble at $N=200$ 
gives noticeably {\it better} success probabilities than the Gaussian ensemble,
and it also gives slightly {\it better} success probabilities on average.

The Hadamard suite $Z$-scores are typical of what we have already seen
in suites 11,12, and 15.  There is some very noticeable bulk misfit at $N=256$,  and much less
at $N=512$. We do not have data at larger $N$.
The misfit at $N=256$ is mostly in the upper tail, meaning that 
there are some situations where the Hadamard ensemble at $N=256$ 
gives noticeably {\it worse} success probabilities than the Gaussian ensemble.
This upper tail effect has mostly disappeared by $N=512$,
however for the suite with positive coeficients
at $N=512$ we see that the slope of the $Z$-score graph is
higher than 1, so the bulk distribution of $Z$-scores
has seemingly higher standard deviation than 1.0.
This is a sign of {\it overdispersion}, i.e. the underlying
success probability $p_1$ might best be treated as a random variable
which varies from realization to realization, but which has the correct expectation
$p_0$.

\subsubsection{Linear models of the $Z$-scores}
Above we reported fits explaining the
observed $Z$-scores using two models of the form (\ref{satmodel})-(\ref{linmodel}).
We considered two such models, one with 
the restrictions $\alpha_0 = 0$ and $\beta_0 = 0$, and one without.

The display below shows the output of the $R$ linear model
command for the restricted model.  
The symbol {\tt de} denotes  $\delta - 1/2$
and the symbol {\tt probSize} denotes 
$1/\sqrt{N}$.   The output lists interaction
terms such as {\tt probSize:de}, which is what we have called
$\beta_1(N)$.

\begin{verbatim}
Call:
lm(formula = Z ~ probSize * de + probSize * dea2 - de - dea2 - 1)

Residuals:
     Min       1Q   Median       3Q      Max 
-3.08857 -0.65803  0.04351  0.71908  3.04688 

Coefficients:
              Estimate Std. Error t value Pr(>|t|)    
probSize        -3.357      1.684  -1.993   0.0466 *  
probSize:de      2.640      6.410   0.412   0.6805    
probSize:dea2  -71.876     14.293  -5.029 6.17e-07 ***
---
Residual standard error: 1.037 on 756 degrees of freedom
Multiple R-Squared: 0.1641,	Adjusted R-squared: 0.1608 
F-statistic: 49.49 on 3 and 756 DF,  p-value: < 2.2e-16 
\end{verbatim}

We have these key points:
\bitem
 \item Scaling with $N$:  the mean $Z$-score tends toward zero
 with increasing $N$.
 \item Standard deviation $\approx 1$. The standard devation of
 the residual $Z$-scores is close to 1.
 \item Hinge terms.  The linear trend term is not significant 
 but the hinge term is.
 \eitem

When we fit the unrestricted model, allowing non scaling terms
that do not decay with increasing $N$, we see that in fact
the only term approaching significance is again the scaling Hinge term.
\begin{verbatim}
Call:
lm(formula = Z ~ probSize * de + probSize * dea2)

Residuals:
     Min       1Q   Median       3Q      Max 
-3.08564 -0.66130  0.04149  0.71351  3.04902 

Coefficients:
               Estimate Std. Error t value Pr(>|t|)  
(Intercept)    -0.05629    0.18748  -0.300   0.7641  
probSize       -2.19063    4.15678  -0.527   0.5983  
de              0.51222    0.71232   0.719   0.4723  
dea2            0.03154    1.61315   0.020   0.9844  
probSize:de    -7.85130   15.83140  -0.496   0.6201  
probSize:dea2 -72.84768   35.45530  -2.055   0.0403 *
---

Residual standard error: 1.038 on 753 degrees of freedom
Multiple R-Squared: 0.1179,	Adjusted R-squared: 0.112 
F-statistic: 20.12 on 5 and 753 DF,  p-value: < 2.2e-16 
\end{verbatim}
The unrestricted model has a worse adjusted $R^2$
than the  scaling model and does not produce a meaningfully 
better residual standard deviation.
We conclude that the scaling model fits the Rademacher
ensemble adequately. Also, since the form of the model
and the particular analyses we are presenting were both
specified before the data became available,
the statistics can be considered inferentially rigorous.

The situation with the Hadamard ensemble is quite different.
It belongs with the exceptional ensembles, in the sense that
substantial lack of fit is evident, particularly at small $\delta$ and
and small $N$.  Also, since we have only two values of $N$ -- $256$ and $512$ --
it is unclear that fitting the linear model will in any way validate scaling.

\subsubsection{Lack of fit at the Hadamard ensembles}

\begin{figure}
\begin{center}
\includegraphics[width=4in]{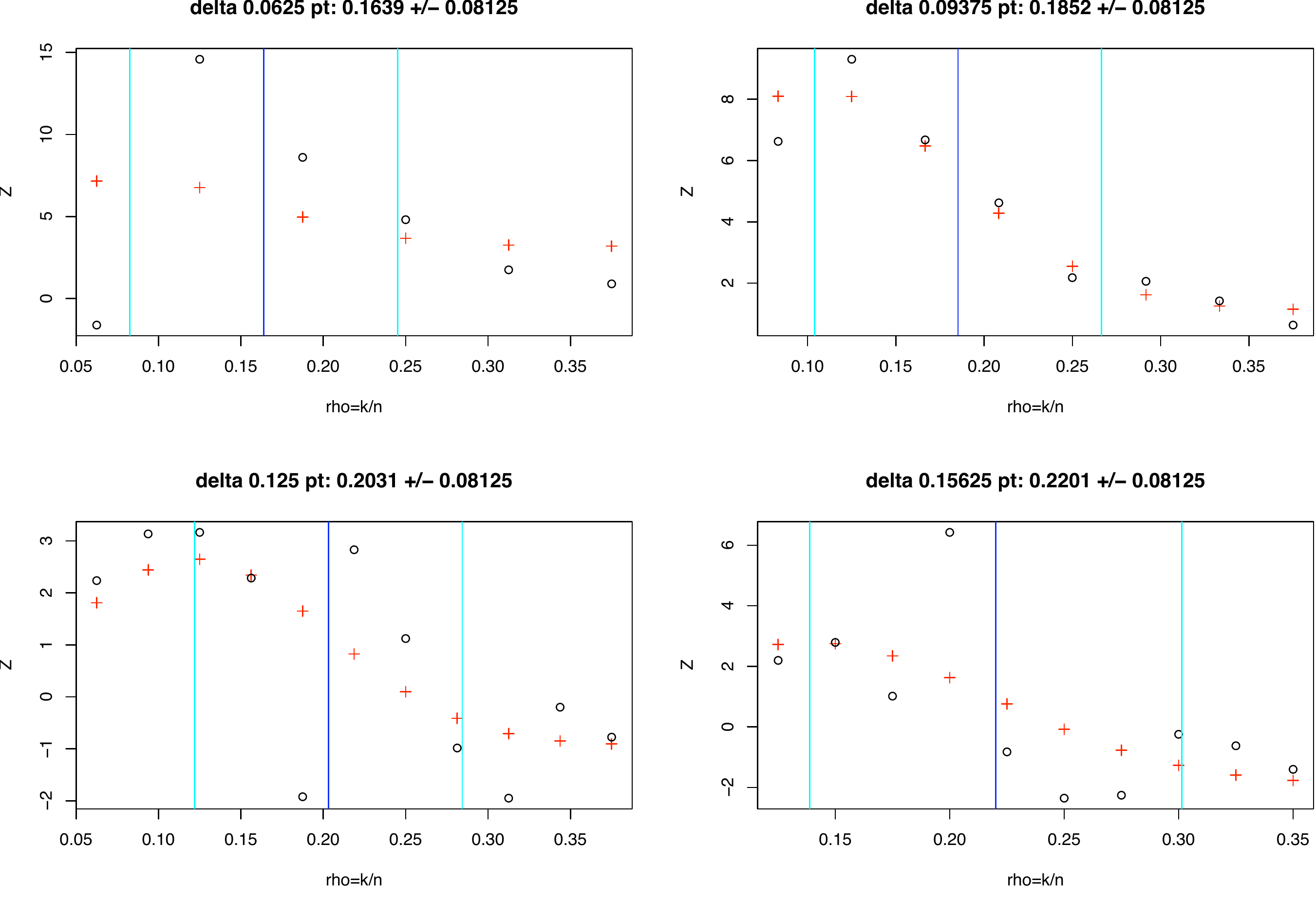}
\caption{$Z$-scores for suite 14, $N=256$ plotted versus $\rho = k/n$ at
each value of $\delta=n/N$. Vertical bars indicate asymptotic phase
transition and nominal width of the transition zone. Red Crosses
indicate fits using the displaced $LD50$ model developed earlier.
This model describes reasonably well the lack-of-fit.}
\label{fig-HadamardZvsAdjustPT}
\end{center}
\end{figure}

The Hadamard ensemble displays, at $N=256$, the lack of fit
which is familiar to us from the exceptional ensembles 11,12, and 15.
Figure \ref{fig-HadamardZvsAdjustPT} displays 
the behaviour of $Z$-scores as a function of $\rho =k/n$
within each constant-$\delta$ slice, together with
fits using the displaced $LD50$ model which we
showed described the other exceptional ensembles, 
with the same displacement of one transition width.
Note that this model was developed entirely before these
data became available. Hence this display
presents a validation of the earlier approach.

\begin{figure}
\begin{center}
\includegraphics[width=4in]{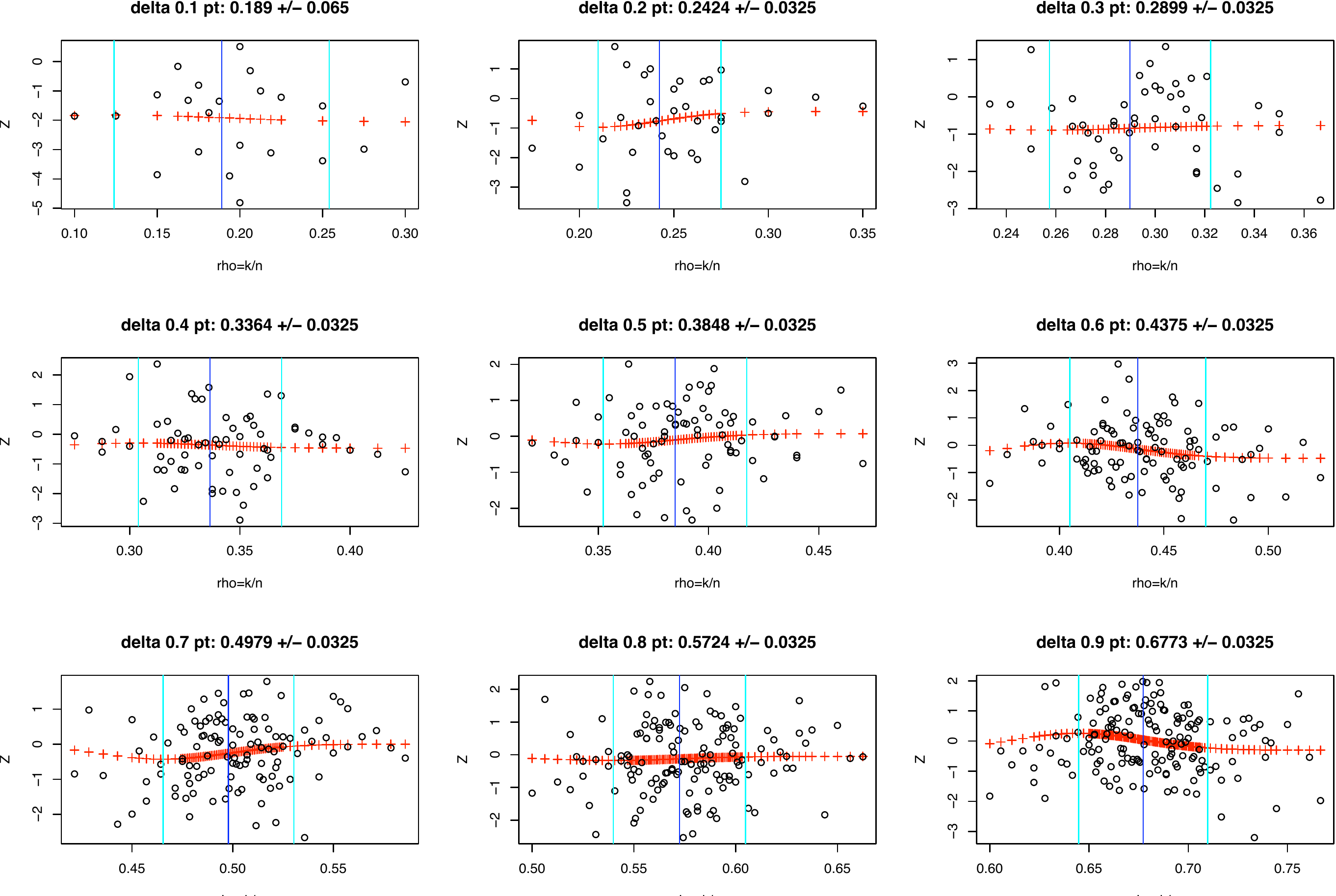}
\caption{$Z$-scores for suite 20, $N=200$ plotted versus $\rho = k/n$ at
each value of $\delta=n/N$. Vertical bars indicate asymptotic phase
transition and nominal width of the transition zone. Red Crosses
indicate fits using displaced $LD50$ model.
This lack of fit is not substantial and
the model does not capture any that might be evident.}
\label{fig-RademacherZvsAdjustPT}
\end{center}
\end{figure}

One can see in figure \ref{fig-HadamardZvsAdjustPT} 
that the displaced $LD50$ model indeed adequately describes the structure
at small $\delta$ in the Hadamard ensemble at $N=256$.  For comparison,
we show in figure \ref{fig-RademacherZvsAdjustPT} 
which demonstrates that no such model is needed in the 
Rademacher ensemble.

\subsubsection{Implications of the validation ensembles}

Our analysis of the validation ensembles follows the script
we arrived at in earlier analysis. Since the models
being fit -- which in earlier analysis {\it depended on the same data}
as those being explained -- 
now do not depend on the data being explained, we gain  additional confidence
about the model we developed.  In particular, the
partitioning into ordinary and exceptional ensembles, the modelling
of mean $Z$-scores with power law means having decay $N^{-1/2}$
in ordinary ensembles  and
the existence of transient effects in exceptional ensembles
at small $N$ and $\delta$ -- all seem confirmed in our validation
ensembles.
 
\subsection{Conclusions}

We made extensive computational
experiments, solving underdetermined systems of equations
$y = Ax$ using linear programming based algorithms
in cases where a known solution
$x_0$ is a $k$-sparse vector,
and $A$ is an $n$ by $N$ matrix with $n<N$.
We tracked the success of the linear
programming algorithms at
recovering this known sparse solution.

We considered  millions
of  such systems at a range of problem
sizes and covering a range of random matrix ensembles. 
Our work would have required more than 6 years
using a single modern desktop computer.

The Gaussian ensemble is well understood
in the large-problem limit based on prior
theoretical work.   It is known there is
a phase transition in the $(\delta,\rho)$
phase diagram at the curve $(\delta,\rho(\delta;Q))$.
For $(\delta, \rho)$ below this curve we
see success -- the probability of success tends to $1$ --
 and above this curve we
see failure -- the probability of success tends to $0$.
Here $\delta =n/N$ is a measure of matrix
shape and $\rho = k/n$ is a measure
of sparsity of the solution vector.

Empirical work at the Gaussian ensemble
shows that for finite $N$ the location of 
50\% success closely matches the asymptotic phase
transition $\rho(\delta;Q)$.  In fact the
success probability behaves with $\rho$
as $ \bar{\Phi}((\rho -\rho(\delta;Q))/w(\delta,N))$
where $\bar{\Phi}$
denotes the $N(0,1)$ survival function.
and $w(\delta,N)$ is a width parameter.
In picturesque terms there is a
transition zone: for $\rho < \rho(\cdot;Q) - 3w$
success is overwhelmingly likely;
for $\rho > \rho(\cdot;Q) + 3w$ success is
overwhelmingly unlikely, and in between
there is a transition zone of width $\propto w$.
$w$ tends to zero with increasing problem size $N$
as $O(N^{-1/2})$.

Broadly similar behaviour is evident at  all but two of the
non-Gaussian matrix ensembles considered here.
Hence, the evidence points to asymptotic phase
transitions at the conforming matrix ensembles 
which are all located at precisely the same place as in the Gaussian
case.  This is the hypothesis of weak universality
advanced and maintained in the main text.


The above conclusions emerged
from using standard two-sample statistical
inference tools to compare results from  non-Gaussian
ensembles with their Gaussian counterparts at
the same problem size and sparsity level. We
considered the behaviour of over 16,000 $Z$-scores
and  observed good bulk agreement of the two-sample $Z$-scores
with the $N(0,1)$ distribution, supporting
the null hypothesis of no difference; however,
fitting a linear model to an array of such $Z$-scores
we were able to identify statistically significant nonzero mean $Z$-scores
varying with problem size, with undersampling fraction $\delta=n/N$,
and, at small $\delta$, with $\rho$. 
The means exhibit trends  varying from ensemble
to ensemble, but they decay with problem
size $N$ 
as $O(N^{-1/2})$. 
The nonzero means  are consistent with weak, `asymptotic', universality
but not with strong,  finite-$N$, universality.  After accounting for
these means, the $Z$-scores exhibit standard deviations
close to one, with a small fraction of exceptions.

For the non-Gaussian ensembles studied here
the fact that  phase diagram behaviour matches the Gaussian case
goes far beyond current theory (e.g.  Adamczak {\em et al.} (2009) )

\label{lastpage}
\end{document}